\documentclass[10pt]{article}

\usepackage{amsmath}
\usepackage{amsmath, amsthm, amsfonts, amssymb, color}
\usepackage{mathrsfs}
\usepackage{latexsym,bm}
\usepackage{graphicx}
\usepackage[toc,page]{appendix}

\newtheorem{theorem}{Theorem}[section]
\newtheorem{proposition}[theorem]{Proposition}

\newtheorem{lemma}[theorem]{Lemma}
\newtheorem{corollary}[theorem]{Corollary}

\newtheorem{remark}[theorem]{Remark}

\definecolor{wco}{rgb}{0.5,0.2,0.3}

\newcommand{\p}{\mathrm{P}}

\newcommand{\pp}{\mathbb{P}}

\newcommand{\mc}{\mathcal{M}_{c}(\mathbb{R})}
\newcommand{\mf}{\mathcal{M}_{F}(\mathbb{R})}

\newcommand{\R}{\mathbb{R}}

\numberwithin{equation}{section} 

\begin{document}

\allowdisplaybreaks

\title{\bf The extremal process of super-Brownian motion: a probabilistic approach via skeletons
\footnote{The research of this project is supported by the National Key R\&D Program of China (No. 2020YFA0712900).}
     }
\author{{\bf Yan-Xia Ren}\thanks{The research of this author is supported by   NSFC (Grant Nos. 11731009 and 12071011) and LMEQF.}
\quad
{\bf Ting Yang }\thanks{The research of this  author is  supported by NNSFC 11731009.}
\quad
{\bf Rui Zhang}\thanks{The research of this  author is  supported by  NSFC (Grant No. 11601354), Beijing Municipal Natural Science Foundation(Grant No. 1202004), and Academy for Multidisciplinary Studies, Capital Normal University.}
}

\date{}
\maketitle

\begin{abstract}
Recently Ren et al.
[Stoch. Proc. Appl., 137 (2021)] have proved
that the extremal process of the super-Brownian motion converges in distribution in the limit of large times. Their techniques rely heavily on the study of the convergence of solutions to the Kolmogorov-Petrovsky-Piscounov equation along the lines of
[M. Bramson, Mem. Amer. Math. Soc., 44 (1983)].
In this paper we take a different approach. Our approach is based on the skeleton decomposition of super-Brownian motion.
The skeleton may be interpreted as immortal particles that determine the large time behaviour of the process. We exploit this fact and carry asymptotic properties from the skeleton over to the super-Brownian motion.
Some new results concerning the probabilistic representations of the limiting process are obtained, which cannot be directly obtained through the results of [Y.-X. Ren et al., Stoch. Proc. Appl., 137 (2021)]. Apart from the results,
our approach
offers insights into the driving force behind the limiting process for super-Brownian motions.

\end{abstract}

\medskip


\medskip


\section{Introduction}

\subsection{Super-Brownian motion and the skeleton space}\label{sec1}

Throughout this paper, we use
``$:=$"
as a way of definition.
 Suppose $\mf$ is the space of finite measures on $\mathbb{R}$ equipped with the topology of weak convergence.
 The set of finite and compactly supported measures on $\mathbb{R}$ is denoted by $\mc$.
 We use $\mathcal{B}_{b}(\R)$ (respectively,
$\mathcal{B}^{+}(\R)$) to denote
the space of bounded (respectively,
nonnegative) Borel functions on $\R$.
The space of continuous (and compactly supported) functions on $\R$ will be denoted as $C(\R)$ (and $C_{c}(\R)$ resp.).
 We use the notation $\langle f,\mu\rangle:=\int_{\R}f(x)\mu(dx)$ and $\|\mu\|:=\langle 1,\mu\rangle$.
 The main process of interest in this paper is an $\mf$-valued Markov process $X=\{X_{t}:t\ge 0\}$ with evolution depending on two quantities $P_{t}$ and $\psi$. Here $P_{t}$ is the semigroup of the standard Brownian motion
 $\{((B_{t})_{t\ge0},\Pi_x):x\in\R\}$
 and $\psi$ is the so-called
branching mechanism,
 which takes the form
 \begin{equation}
 \psi(\lambda)=-\alpha {\lambda}+\beta{\lambda}^{2}+\int_{(0,{\infty})}\left(e^{-{\lambda} y}-1+{\lambda} y\right)\pi(dy)\quad\mbox{ for }{\lambda}\ge 0,\label{2.1}
 \end{equation}
 with $\alpha>0$, $\beta\ge 0$ and  $\pi(dy)$ being a measure concentrated on $(0,{\infty})$ such that
 $
 \int_{(0,+\infty)}\left(y\wedge y^{2}\right)\pi(dy)<+\infty
 $.
 The distribution of $X$ is denoted by $\p_{\mu}$ if it is started
 at $\mu\in\mf$ at $t=0$.
 With abuse of notation, we also use $\p_{\mu}$ to denote the expectation with respect to $\p_{\mu}$.
 $X$ is called a
(supercritical) \textit{$(P_{t},\psi)$-superprocess} or \textit{super-Brownian motion with branching mechanism $\psi$}
if $X$ is
 an
 $\mf$-valued process such that
 for any
 $\mu\in \mf$, $f\in\mathcal{B}^{+}_{b}(\mathbb{R})$ and $t\ge 0$,
 \begin{equation}\label{eq0}
 \p_{\mu}\left[e^{-\langle f,X_{t}\rangle}\right]=e^{-\langle u_{f}(t, \cdot),\mu\rangle},
 \end{equation}
where
\begin{equation}\label{eq0'}
u_{f}(t,x):=-\log \p_{\delta_{x}}\left(e^{-\langle f,X_{t}\rangle}\right)
\end{equation}
is
 the unique nonnegative locally bounded solution
 to the following integral equation:
 \begin{equation}\nonumber
 u_{f}(t,x)=P_{t}f(x)-\int_{0}^{t}P_{s}\left(\psi(u_{f}(t-s, \cdot))\right)(x)ds \quad
 \mbox{ for any }x\in \mathbb{R} \mbox{ and }t\ge 0.
 \end{equation}
We note that $u_{f}(t,x)$ is also a solution to the partial differential equation
 \begin{equation}\label{eq1}
 \frac{\partial}{\partial t}u(t,x)=\frac{1}{2}\frac{\partial^{2}}{\partial x^{2}} u(t,x)-\psi(u(t,x))
 \end{equation}
 with initial condition $u(0,x)=f(x)$.
 Moreover, if $f$ is a nonnegative bounded continuous function on $\R$,  $\lim_{t\to 0}u_{f}(t,x)=f(x)$, for $x\in\R$.
We will refer to \eqref{eq1} as the Kolmogorov-Petrovsky-Piscounov (K-P-P) equation.
The existence of such a process $X$ is established by \cite{D93}. Moreover, such a super-Brownian motion $X$
has a Hunt realization in $\mf$ (see, for example, \cite[Theorem 5.12]{Li}) such that $t\mapsto \langle f,X_{t}\rangle$ is almost surely right continuous
for any bounded continuous functions $f$.
We shall always work with this version.

 The $(P_{t},\psi)$-superprocess can be constructed as the high density limit of a sequence of branching Markov processes. Another link between superprocesses and branching Markov processes is provided by the so-called skeleton decomposition, which is developed by
 \cite{BKM,CRY,EKW,KPR}.
 The skeleton decomposition provides a pathwise description of
a superprocesses in terms of immigrations along a branching Markov process called the skeleton.
The following condition is fundamental for the skeleton construction.
\medskip

\textbf{(A1).}\quad $\psi(+\infty)=+\infty$.
\medskip

\noindent Condition (A1) implies that there exists some $\lambda^{*}\in (0,+\infty)$ such that $\psi(\lambda^{*})=0$.
The key property of $\lambda^{*}$ used in the skeleton construction is that
$$\p_{\mu}\left(\lim_{t\to+\infty}\|X_{t}\|=0\right)=\mathrm{e}^{-\lambda^{*}\|\mu\|}\quad\forall \mu\in\mf,$$
and so $\lambda^{*}$ gives rise to the multiplicative $\p_{\mu}$-martingale
$(e^{-\lambda^{*}\|X_{t}\|})_{t\ge 0}$. A more detailed description of skeleton construction is given in Proposition \ref{prop1}.

\begin{proposition}[Skeleton space]\label{prop1}
Suppose (A1) holds.
For every $\mu\in \mf$ and
every finite point measure $\nu$ on $\R$,
there exists a probability space with probability measure $\pp_{\mu,\nu}$ that carries three processes: $(Z_{t})_{t\ge 0}$, $(I_{t})_{t\ge 0}$ and $(X^{*}_{t})_{t\ge0}$, where
$\left((Z_{t})_{t\ge 0};\pp_{\mu,\nu}\right)$ is a branching Brownian motion with  branching rate $q>0$ and offspring distribution $\{p_{k}:k\ge 2\}$
uniquely determined by
\begin{equation}\label{eq2}
q(F(s)-s)=\frac{1}{\lambda^{*}}\psi(\lambda^{*}(1-s))\quad\forall s\in [0,1],
\end{equation}
here $F(s):=\sum_{k=2}^{+\infty}p_{k}s^{k}$, and $\pp_{\mu,\nu}(Z_{0}=\nu)=1$;
$\left((X^{*}_{t})_{t\ge 0};\pp_{\mu,\nu}\right)$ is a subcritical
super-Brownian motion
with branching mechanism
$\psi^{*}(\lambda)=\psi(\lambda+\lambda^{*})$ and
$\pp_{\mu,\nu}(X^{*}_{0}=\mu)=1$; $((I_{t})_{t\ge 0},\pp_{\mu,\nu})$ is an $\mf$-valued process with $\pp_{\mu,\nu}\left(I_{0}=0\right)=1$, which denotes the immigration at time $t$ that occurred along the skeleton $Z$. Under $\pp_{\mu,\nu}$, both $Z$ and $I$ are independent of $X^{*}$.

If $\pp_{\mu}$ denotes the measure $\pp_{\mu,\nu}$ with $\nu$ replaced by a Poisson random measure with intensity $\lambda^{*}\mu(dx)$,
then $\left(
\widehat X :=
X^{*}+I;\pp_{\mu}\right)$ is Markovian and
has the same distribution as
$(X;\p_{\mu})$.
Moreover, under $\pp_{\mu}$, given $\widehat  X_{t}$,
the measure $Z_{t}$ is a Poisson random measure with intensity
$\lambda^{*}\widehat X_{t}(dx)$.
\end{proposition}

Since
$({\widehat X};\pp_{\mu})$
is equal in distribution to the $(P_{t},\psi)$-superprocess $(X;\p_{\mu})$, we may work on this skeleton space whenever it is convenient.
 For notational simplification, we will abuse the notation and denote $\widehat X$ by $X$.
We will refer to $(Z_{t})_{t\ge 0}$ as the skeleton branching Brownian motion (skeleton BBM) of $X$.
Since the distributions of $X^{*}$ (resp. $I$) under $\pp_{\mu,\nu}$ do not depend on $\nu$ (resp. $\mu$), we sometimes write $\pp_{\mu,\cdot}$ (resp. $\pp_{\cdot,\nu}$) for $\pp_{\mu,\nu}$.

The aim of this article is to
study the asymptotic behavior of the super-Brownian motion.
Since the solutions to the K-P-P equation fully capture the space-time behavior of the super-Brownian motion, it is natural to use this relationship to investigate the large time behaviour for super-Brownian motions.
K-P-P equation has been studied extensively using both analytic and probabilistic methods (see, for example,
\cite{Bramson,CR,Kyprianou,KLMR}).
Among these papers,
the seminal work of Bramson \cite{Bramson} established the convergence of $u(t,x+m(t))$ to travelling waves for some centering term $m(t)$.
We use $u\in Z_{t}$ to denote a particle of the skeleton BBM which is alive at time $t$ and $z_{u}(t)$ for its spatial location at $t$.
Based on the results of \cite{Bramson},
\cite{ABBS,ABK}
proved that for a branching Brownian motion, the extremal process, namely the random point measure
$$\mathcal{E}^{Z}_{t}:=\sum_{u\in Z_{t}}\delta_{z_{u}(t)-m(t)}$$
converges in distribution to a limiting process as $t\to+\infty$,
and gave an explicit construction of the limiting process.
In Section \ref{sec2.1} we will review some of the facts concerning on the convergence of
the solutions to the K-P-P equation
with applications to branching Brownian motions.

It is to be noted that Bramson \cite{Bramson} considers the solutions to \eqref{eq1} with initial condition taking values in $[0,1]$, and thus cannot be applied directly to super-Brownian motions.
A key step is to establish the convergence of solutions to \eqref{eq1} for more general initial conditions.
Methodologically, the convergence result established in \cite{Bramson} relies mainly on approximations of the solutions by Feynman-Kac formula. Bramson's reasoning can be applied, with modifications, to solutions of the K-P-P equation with initial conditions
taking values in $[0,\infty)$.
This method has been adopted by Ren et al. in their recent work \cite{RSZ}. In this paper we shall offer a different approach. We appeal to the skeleton techniques for superprocesses.
Intuitively, the super-Brownian motion may be interpreted as a cloud of subcritical diffusive mass immigrating off a supercritical branching Brownian motion, the skeleton, which governs the large time behaviour of the process. We exploit this fact and carry the long time behaviour from the skeleton over to the super-Brownian motion.
Our work is partly inspired by
\cite{CRY,EKW}, where the skeleton techniques have been used successfully to establish the laws of large numbers for superprocesses. Apart from the result itself, our approach provides structural insights into the driving force behind the limiting process for the super-Brownian motions.

We observe that up to a space-time scaling transform, the branching mechanism $\psi$ can be assumed to satisfy that
\begin{equation}\label{condi0}
\psi'(0)=-1\mbox{ and }\lambda^{*}=1.
\end{equation}
In fact, if $u(t,x)$ is a solution to \eqref{eq1} then $u(\alpha^{-1}t,\alpha^{-1/2}x)/\lambda^{*}$ is a solution to \eqref{eq1} with $\psi$ replaced by $\widetilde{\psi}$ where $\widetilde{\psi}(\lambda)=\psi(\lambda^{*}\lambda)/\alpha\lambda^{*}$ satisfies that $\widetilde{\psi}'(0)=-1$ and $\tilde{\psi}(1)=0$. This implies that
for a
$(P_{t},\psi)$-superprocess
$(X_{t})_{t\ge 0}$, if we define the random measures $\tilde{X}_{t}$ by
\begin{equation}\label{tildex}
\langle f,\widetilde{X}_{t}\rangle =\lambda^{*}\langle f(\alpha^{1/2}\cdot),X_{\alpha^{-1}t}\rangle \quad\forall t\ge 0,\ f\in\mathcal{B}_{b}^{+}(\mathbb{R}),
\end{equation}
then $(\widetilde{X}_{t})_{t\ge 0}$ is a
$(P_{t},\widetilde{\psi})$-superprocess.
It suffices to study
the long time behavior of $\widetilde{X}_{t}$.
Therefore, in the rest of this paper the branching mechanism $\psi$ is assumed to satisfy \eqref{condi0}, which will simplify computations and notations later on.

\subsection{The extremal process of the skeleton and facts}\label{sec2}

\subsubsection{Derivative martingales for the skeleton}\label{sec2.0}
In this and the next two subsections we assume (A1) and \eqref{condi0} hold. Additional conditions used are stated explicitly.
Recall that $u\in Z_{t}$ and $z_{u}(t)$ denote, respectively, a particle of the skeleton BBM which is alive at time $t$ and its spatial location at $t$.
Define for $t\ge 0$,
$$\partial M_{t}:=\sum_{u\in Z_{t}}\left(\sqrt{2}t-z_{u}(t)\right)\mathrm{e}^{\sqrt{2}\left(z_{u}(t)-\sqrt{2}t\right)}$$
	It is known that $((\partial M_{t})_{t\ge 0},\pp_{\cdot,\nu})$ is a signed martingale for every compactly supported finite point measure $\nu$ on $\R$, which is referred to as the derivative martingale of the skeleton BBM $(Z_{t})_{t\ge 0}$. This martingale is deeply related to the travelling wave solutions to the K-P-P equation and plays an important role in the limit theory of the skeleton BBM. \cite{Kyprianou} has proved that the martingale $((\partial M_{t})_{t\ge 0},\pp_{\cdot,\delta_{0}})$ has an almost sure nonnegative limit, and later \cite{RY} established the sufficient and necessary condition for the limit to be non-degenerate. We give the statement below which reproduces the same results in the setting of skeleton space.

\begin{proposition}\label{prop:llogl}
	Suppose $\mu\in\mc$. The limit $\partial M_{\infty}=\lim_{t\to+\infty}\partial M_{t}$ exists and is nonnegative $\pp_{\mu}$-a.s.
	Moreover, if
	\begin{equation}\label{llogl}
		\int_{(1,+\infty)}x(\log x)^{2}\pi(dx)<+\infty,
	\end{equation}
     then
    $$\pp_{\mu}\left(\partial M_{\infty}=0\right)=\mathrm{e}^{-\|\mu\|}.$$
	If \eqref{llogl} fails, then $\partial M_{\infty}=0$ $\pp_{\mu}$-a.s.
\end{proposition}
	
\proof
    By decomposing $\partial M_{t}$ into contributions derived from the population at time $s\in [0,t)$, one has
	$$\partial M_{t}=\sum_{u\in Z_{s}}\sum_{v\in Z_{t},u\prec v}\left(\sqrt{2}t-z_{v}(t)\right)\mathrm{e}^{\sqrt{2}\left(z_{v}(t)-\sqrt{2}t\right)}.$$
	Here $u\prec v$ means that $u$ is an ancestor of $v$.
	We use $z^{(u)}_{v}(t-s)$ to denote $z_{v}(t)-z_{u}(s)$. Then we have
	\begin{eqnarray*}
		\partial M_{t}&=&\sum_{u\in Z_{s}}\sum_{v\in Z_{t},u\prec v}\left(\sqrt{2}s-z_{u}(s)+\sqrt{2}(t-s)-z^{(u)}_{v}(t-s)\right)\mathrm{e}^{\sqrt{2}\left(z_{u}(s)-\sqrt{2}s+z^{(u)}_{v}(t-s)-\sqrt{2}(t-s)\right)}\\
		&=&\sum_{u\in Z_{s}}\mathrm{e}^{\sqrt{2}(z_{u}(s)-\sqrt{2}s)}\left[\sum_{v\in Z_{t},u\prec v}\left(\sqrt{2}(t-s)-z^{(u)}_{v}(t-s)\right)\mathrm{e}^{\sqrt{2}\left(z^{(u)}_{v}(t-s)-\sqrt{2}(t-s)\right)}\right]\\
		&&\quad+\sum_{u\in Z_{s}}\left(\sqrt{2}s-z_{u}(s)\right)\mathrm{e}^{\sqrt{2}(z_{u}(s)-\sqrt{2}s)}\left[\sum_{v\in Z_{t},u\prec v}\mathrm{e}^{\sqrt{2}\left(z^{(u)}_{v}(t-s)-\sqrt{2}(t-s)\right)}\right]\\
		&=:&\sum_{u\in Z_{s}}\mathrm{e}^{\sqrt{2}(z_{u}(s)-\sqrt{2}s)}\partial M^{(u)}_{t-s}+\sum_{u\in Z_{s}}(\sqrt{2}s-z_{u}(s))\mathrm{e}^{\sqrt{2}(z_{u}(s)-\sqrt{2}s)}M^{(u)}_{t-s}.
	\end{eqnarray*}
	In particular by setting $s=0$, we have
	\begin{equation}\label{3.3}
		\partial M_{t}=\sum_{u\in Z_{0}}\mathrm{e}^{\sqrt{2}z_{u}(0)}\partial M^{(u)}_{t}-\sum_{u\in Z_{0}}z_{u}(0)\mathrm{e}^{\sqrt{2}z_{u}(0)}M^{(u)}_{t}.
	\end{equation}
	Define for $t\ge 0$,
	$$M_{t}:=\sum_{u\in Z_{t}}\mathrm{e}^{\sqrt{2}(z_{u}(t)-\sqrt{2}t)}.$$
	It is easy to see that for $u\in Z_{0}$, $\partial M^{(u)}_{t}$ and $M^{(u)}_{t}$ are independent copies of $(\partial M_{t},\pp_{\cdot,\delta_{0}})$ and
	$(M_{t},\pp_{\cdot,\delta_{0}})$, respectively.
By \cite{Kyprianou} one has $\pp_{\cdot,\delta_{0}}\left(\lim_{t\to+\infty}\partial M_{t}\mbox{ exists and is nonnegative}\right)=1$ and $\pp_{\cdot,\delta_{0}}\left(\lim_{t\to+\infty}M_{t}=0\right)=1$.
We also note that $Z_{0}$ is a Poisson random measure with compactly supported intensity $\mu$.
So each of the sums in the right hand side of \eqref{3.3} contains
finite terms almost surely. These facts together with \eqref{3.3} imply that the limit $\partial M_{\infty}=\lim_{t\to+\infty}\partial M_{t}$ exists and is nonnegative
	$\pp_{\mu}$-a.s. for every $\mu\in\mc$.

By letting $t\to +\infty$ in \eqref{3.3}, we have
	\begin{equation}\label{n1}
		\partial M_{\infty}\stackrel{\mbox{d}}{=}\sum_{u\in Z_{0}}\mathrm{e}^{\sqrt{2}z_{u}(0)}\partial M^{(u)}_{\infty}.
	\end{equation}
	where for $u\in Z_{s}$, $\partial M^{(u)}_{\infty}$ are independent copies of $(\partial M_{\infty},\pp_{\cdot,\delta_{0}})$.
	Recall that $(Z_{0},\pp_{\mu})$ is a Poisson random measure with intensity $\mu(dx)$. Using the Poisson computations, we get by \eqref{n1} that
	\begin{equation}\label{prop2.2}
		\pp_{\mu}\left[\mathrm{e}^{-\lambda\partial M_{\infty}}\right]=\exp\left\{-\int_{\mathbb{R}}\left(1-\pp_{\cdot,\delta_{0}}\left[\mathrm{e}^{-\lambda \mathrm{e}^{\sqrt{2}x}\partial M_{\infty}}\right]\right)\mu(dx)\right\},\quad\forall \lambda>0.
	\end{equation}
	By letting $\lambda\to+\infty$ we have
	\begin{equation}\label{prop2.4}
		\pp_{\mu}\left(\partial M_{\infty}=0\right)=\mathrm{e}^{-\left(1-\pp_{\cdot,\delta_{0}}\left(\partial M_{\infty}=0\right)\right)\|\mu\|}.
	\end{equation}
	By \cite{RY} $\pp_{\cdot,\delta_{0}}(\partial M_{\infty}=0)=1$ if $\sum_{k}k(\log k)^{2}p_{k}=+\infty$,
and otherwise $\pp_{\cdot,\delta_{0}}(\partial M_{\infty}=0)=0$.
On the other hand by Lemma \ref{lemA1} $\sum_{k}k(\log k)^{2}p_{k}$ is finite if and only if so is $\int_{(1,+\infty)}x(\log x)^{2}\pi(dx)$. Thus it follows by \eqref{prop2.4} that $\pp_{\mu}\left(\partial M_{\infty}=0\right)=\mathrm{e}^{-\|\mu\|}$ if $\int_{(1,+\infty)}x(\log x)^{2}\pi(dx)<+\infty$, and $\pp_{\mu}\left(\partial M_{\infty}=0\right)=1$ if $\int_{(1,+\infty)}x(\log x)^{2}\pi(dx)=+\infty$. Hence we complete the proof.\qed
	
\bigskip

\subsubsection{The extremal process of the skeleton}\label{sec2.1}
Let $\mathcal{M}(\R)$ be the space of all  Radon measures on $\R$ equipped with the vague topology.
A sequence $\{\Xi_{n}:n\ge 1\}$ of random Radon measures on $\R$
is said to converge in distribution to $\Xi$ if and only if for all
$\phi\in C^{+}_{c}(\mathbb{R})$,
the random variables $\langle \phi,\Xi_{n}\rangle$ converges in distribution to $\langle \phi,\Xi\rangle$. For $x\in\mathbb{R}$ and a function $f$ on $\mathbb{R}$, we define the shift operator $\mathcal{T}_{x}$ by $\mathcal{T}_{x}f(y):=f(x+y)$ for all $y\in \mathbb{R}$.
For $\mu\in\mathcal{M}(\R)$,
we use $\mu+x$ and sometimes $\mathcal{T}_{x}\mu$ to denote the measure induced by $\mathcal{T}_{x}$, that is, $\int_{\R}f(y)\mathcal{T}_{x}\mu(dy)=\int_{\mathbb{R}}f(y)(\mu+x)(dy)=\int_{\mathbb{R}}\mathcal{T}_{x}f(y)\mu(dy)$ for all $f\in\mathcal{B}^{+}(\R)$.

Given \eqref{condi0}, \eqref{eq2} can be written as
\begin{equation}
\label{eq2'}
q(F(s)-s)=\psi(1-s),\quad\forall s\in [0,1],
\end{equation}
where $F(s):=\sum_{k=2}^{+\infty}p_{k}s^{k}$.
Let $f:\mathbb{R}\to [0,1]$ be a Borel function.
It is known that the function $(t,x)\mapsto \pp_{\cdot,\delta_{x}}\left[\prod_{u\in Z_{t}}f(z_{u}(t))\right]$ is a solution of the equation
\begin{equation*}
\frac{\partial}{\partial t}u(t,x)=\frac{1}{2}\frac{\partial^{2}}{\partial x^{2}}u(t,x)+q(F(u(t,x))-u(t,x)).
\end{equation*}
with initial condition $u(0,x)=f(x)$. Then $(t,x)\mapsto 1-\pp_{\cdot,\delta_{x}}\left[\prod_{u\in Z_{t}}(1-f(z_{u}(t)))\right]$
is a solution to the equation
\eqref{eq1} with initial condition $u(0,x)=f(x)$.

Recall the definition of $u_{f}(t,x)$ for $f\in\mathcal{B}_{b}^{+}(\R)$ given
in \eqref{eq0'}.
We note that $u_{f}(t,x)$ is the unique nonnegative solution to \eqref{eq1} with initial condition $f$.
In particular,
if $\|f\|_{\infty}\le 1$,
\begin{equation}\label{uf}
	u_f(t,x)=1-\pp_{\cdot,\delta_{x}}\left[\prod_{u\in Z_{t}}(1-f(z_{u}(t)))\right],
\end{equation}
where $(Z_{t})_{t\ge 0}$ is the skeleton BBM.
Let $\max Z_{t}:=\max\{z_{u}(t):u\in Z_{t}\}$ be the maximal displacement of the skeleton BBM.
In particular, $\pp_{\cdot,\delta_{x}}\left(\max Z_{t}>0\right)$
is  a solution to \eqref{eq1} with initial condition $1_{(0,+\infty)}(x)$.

Bramson \cite{Bramson} studied the asymptotic behavior of the solution to the
K-P-P equation \eqref{eq1}
with initial condition $u(0,x)$ taking values in $[0,1]$. Actually in \cite{Bramson}, the nonlinear function $-\psi$ can be any function on $[0,1]$ satisfying that
\begin{eqnarray*}
&&\psi\in C^{1}[0,1],
\quad
\psi(0)=\psi(1)=0,
\quad -\psi(u)>0\mbox{ for }u\in (0,1),\\
&&\psi'(0)=-1,\quad
-\psi'(u)\le 1\mbox{ for }0<u\le 1.
\end{eqnarray*}
Bramson \cite{Bramson} shows in particular that if $\psi$ also satisfies that
\begin{equation}\label{bramson1.6}
1+\psi'(u)=O(u^{\rho})\quad\mbox{ as }u\to 0
\end{equation}
for some $\rho>0$, then when the initial condition $u(0,x)$
satisfies a certain integrability condition,
it holds that
$$u(t,x-m(t))\to w(-x)\quad\mbox{ uniformly in }x\in\mathbb{R},\mbox{ as }t\to+\infty.$$
where
\begin{equation}\label{def-m(t)}m(t):=\sqrt{2}t-\frac{3}{2\sqrt{2}}\log t,\end{equation}
and $w(x)$ is a travelling wave solution with speed $\sqrt{2}$, that is, $w(x)$ is the unique (up to translations) solution to the ordinary differential equation
\begin{equation}
\frac{1}{2}w''(x)+\sqrt{2}w'(x)-\psi(w(x))=0,
\end{equation}
with $1-w(x)$ being a distribution function on $\R$.
The integral representation of $w(x)$ is established in \cite{LS} (see also, \cite{Kyprianou,RY}): When $\partial M_{\infty}$ is nondegenerate, one has
\begin{equation}\label{repforw(x)}
w(x)=1-\pp_{\cdot,\delta_{0}}\left[\exp\{-C\partial M_{\infty}\mathrm{e}^{-\sqrt{2}x}\}\right]
\end{equation}
for some constant $C>0$. Moreover, it holds that
\begin{equation}\label{asymforw}
\lim_{x\to+\infty}\frac{w(x)}{x\mathrm{e}^{-\sqrt{2}x}}=C.
\end{equation}
Later \cite{ABK} recovered the above representation of the form \eqref{repforw(x)} and provided an expression for the constant $C$ as a function of the initial condition. As a result \cite{ABK} established the convergence in distribution of the extremal process of the branching Brownian motion.
To apply \cite{ABK,Bramson}'s
results directly to the skeleton BBM, we assume the following condition holds.
\medskip

\textbf{(A2).}  There exists $\beta\in (0,1)$ such that
$$\int_{(1,+\infty)}y^{1+\beta}\pi(dy)<+\infty.$$

\noindent
It is easy to see from Proposition \ref{prop:llogl} that
(A2) is sufficient for $\partial M_{\infty}$ to be nondegenerate. Besides, (A2) is also sufficient for \eqref{bramson1.6}. This is because,
by Lemma \ref{lemA2}, (A2) holds if and only if $\int_{0}^{1}(1+\psi'(s))s^{-(1+\beta)}ds<+\infty$, and the latter implies that
$1+\psi'(s)=O(s^\beta)$.

Define
$$\mathcal{H}_{1}:=\left\{\phi\in\mathcal{B}^{+}_{b}(\mathbb{R}):\ \|\phi\|_{\infty}\le 1\mbox{ and }\int_{0}^{+\infty}y\mathrm{e}^{\sqrt{2}y}\phi(-y)dy<+\infty\right\}.$$

\begin{proposition}\label{lem3.5}
Assume that (A2) holds.
For every $\phi\in\mathcal{H}_{1}$, the limit
\begin{equation}\nonumber
C(\phi):=\lim_{r\to +\infty}\sqrt{\frac{2}{\pi}}\int_{0}^{+\infty}y\mathrm{e}^{\sqrt{2}y}u_{\phi}(r,-y-\sqrt{2}r)dy
\end{equation}
exists and is finite. Moreover,
$$u_{\phi}(t,x-m(t))\to 1-\pp_{\cdot,\delta_{0}}\left[\exp\{-C(\phi)\partial M_{\infty}\mathrm{e}^{\sqrt{2}x}\}\right]
\mbox{ locally uniformly in $x\in\R$, as }
t\to+\infty.$$
\end{proposition}

The proof of Proposition \ref{lem3.5} is given in Section \ref{sec3.1}. Proposition \ref{lem3.5} refines \cite[Lemma 4.10]{ABK}.  In fact, \cite[Lemma 4.10]{ABK} proves this result for
the functions $\phi\in\mathcal{H}_{1}$
which can be represented by $\phi(x)=1-\mathrm{e}^{-g(x)}$ for some $g\in C^{+}_{c}(\mathbb{R})$.

Taking $\phi=1-\mathrm{e}^{-g}$ for some $g\in C^{+}_{c}(\R)$,  one can rewrite $u_{\phi}(t,x)$ as $1-\pp_{\cdot,\delta_{x}}\left[\mathrm{e}^{-\langle g,Z_{t}\rangle}\right]$,
and the above result yields that
\begin{equation}\label{eq:ABK}
\pp_{\cdot,\delta_{0}}\left[\mathrm{e}^{-\langle g,Z_{t}-m(t)\rangle}\right]\to\pp_{\cdot,\delta_{0}}\left[\exp\{-C(\phi)\partial M_{\infty}\}\right]
\quad\mbox{ as }t\to+\infty.
\end{equation}
This implies the convergence in distribution of the extremal process of a branching Brownian motion.
To be more specific, \cite{ABK} proved that,
the extremal process of a branching Brownian motion defined by
$$\mathcal{E}^{Z}_{t}:=Z_{t}-m(t),\quad\forall t\ge 0,$$
converges in distribution, as $t\to\infty$, to a random point measure $\mathcal{E}^{Z}_{\infty}$, which is a decorated Poisson point process with intensity $c_{*}\partial M_{\infty}\sqrt{2}\mathrm{e}^{-\sqrt{2}z}dz$ and decoration law
$\triangle^{Z}$,
where
\begin{equation}\label{def-c*}
c_{*}:=C(1_{(0,+\infty)})=\lim_{r\to+\infty}\sqrt{\frac{2}{\pi}}\int_{0}^{+\infty}z\mathrm{e}^{\sqrt{2}z}\pp_{\cdot,\delta_{0}}\left(\max Z_{r}
>
\sqrt{2}r+z\right)dz\in (0,+\infty),
\end{equation}
and $\triangle^{Z}$ is a random point measure supported on $(-\infty,0]$, with an atom at $0$, which satisfies that
\begin{equation}\label{def-triangleZ}
\mathrm{E}\left[\mathrm{e}^{-\langle \phi,\triangle^{Z}\rangle}\right]=\lim_{t\to+\infty}\pp_{\cdot,\delta_{0}}\left[\mathrm{e}^{-\langle \phi,Z_{t}-\max Z_{t}\rangle}|\max Z_{t}
>
\sqrt{2}t\right],\quad\forall \phi\in\mathcal{C}^{+}_{c}(\mathbb{R}).
\end{equation}
We denote $\mathcal{E}^{Z}_{\infty}$ by DPPP($c_{*}\partial M_{\infty}\sqrt{2}\mathrm{e}^{-\sqrt{2}z}dz$,$\triangle^{Z}$).
\cite{ABBS} proved the same result by a totally different approach almost at the same time as \cite{ABK}.

The limiting extremal process $\mathcal{E}^{Z}_{\infty}$ can be constructed as follows. Given $\partial M_{\infty}$, let $\{e_{i}:i\ge 1\}$ be the atoms of a Poisson point process on $\mathbb{R}$ with intensity $c_{*}\partial M_{\infty}\sqrt{2}\mathrm{e}^{-\sqrt{2}z}dz$ and $\{\triangle^{Z}_{i}:j\ge 1\}$ be a sequence of i.i.d. point measures with the same law as $\triangle^{Z}$, then
$$\mathcal{E}^{Z}_{\infty}
\stackrel{\mbox{d}}{=}
\sum_{i\ge 1}\mathcal{T}_{e_{i}}\triangle^{Z}_{i}.$$

Recall that  $\max Z_{t}=\max\{z_{u}(t):u\in Z_{t}\}$ is the maximal displacement of the skeleton BBM.
 Then
$$u_{1_{(0,+\infty)}}(t,x)=\pp_{\cdot,\delta_{x}}\left(\max Z_{t}>0\right)=\pp_{\cdot,\delta_{0}}\left(\max Z_{t}+x>0\right)$$
is a solution to
\eqref{eq1}
with initial condition $1_{(0,+\infty)}(x)$. Proposition \ref{lem3.5} yields that
for any $x\in\R$,
$$\pp_{\cdot,\delta_{0}}\left(\max Z_{t}-m(t)\le x\right)\to \pp_{\cdot,\delta_{0}}\left[\exp\{-c_{*}\partial M_{\infty}\mathrm{e}^{-\sqrt{2}x}\}\right]\mbox{ as }t\to +\infty.$$
This implies that under the assumptions of Proposition \ref{lem3.5}, the maximal displacement of the skeleton BBM centered by $m(t)$ converges in distribution  to a randomly shifted Gumbel distribution. In fact, Proposition \ref{lem3.5} implies the joint convergence in distribution of $(\mathcal{E}^{Z}_{t},\max \mathcal{E}^{Z}_{t})$, see, for example, \cite[Lemma 4.4]{BBCM}.

\subsubsection{Relation between the limits of the derivative martingales of super-Brownian motion and its skeleton}

Define for $t\ge 0$,
$$\partial W_{t}:=\langle (\sqrt{2}t-\cdot)\mathrm{e}^{\sqrt{2}(\cdot-\sqrt{2}t)},X_{t}\rangle.$$
By \cite{KLMR} for every $\mu\in\mc$,  $((\partial W_{t})_{t\geq 0}, \pp_{\mu})$ is a martingale which is usually called the derivative martingale of the super-Brownian motion.
Obviously $(\partial W_{t})_{t\ge 0}$ is the counterpart of $(\partial M_{t})_{t\ge 0}$ in the setting of superprocess. Interest of this martingale is stimulated by its
close connection
with the travelling wave solutions to the K-P-P equation (see, for example, \cite{KLMR} and the references therein).
Note that $((\partial W_{t})_{t\geq 0}, \pp_{\mu})$ is a signed martingale which does not necessarily converge almost surely.
 To study its convergence, \cite{KLMR} imposed the following condition on $\psi$:
\begin{equation}\label{condi:extra}
\int^{+\infty}_{z}\frac{1}{\sqrt{\int_{1}^{y}\psi(u)du}}dy<+\infty\quad\forall z>1.
\end{equation}
The following result is from \cite[Theorem 2.4]{KLMR}.
\begin{lemma}\label{prop2}
	Assume \eqref{condi:extra} holds.
Then for every $\mu\in\mc$, the limit $\partial W_{\infty}=\lim_{t\to+\infty}\partial W_{t}$ exists and is nonnegative $\pp_{\mu}$-a.s. Moreover,
	$\partial W_{\infty}$ is non-degenerate if and only if \eqref{llogl} holds.
\end{lemma}

\begin{proposition}\label{prop4}
Assume \eqref{condi:extra} holds. Then one has
$(\partial W_{\infty},\pp_{\mu})\stackrel{\mbox{d}}{=}(\partial M_{\infty},\pp_{\mu})$
for any $\mu\in\mc$.
\end{proposition}

\proof Fix an arbitrary $\mu\in\mc$. If \eqref{llogl} fails, then by Proposition \ref{prop:llogl} and Lemma \ref{prop2}, one has $\partial W_{\infty}=\partial M_{\infty}=0$ $\pp_{\mu}$-a.s. Now we assume \eqref{llogl} holds. Let $w(x)$ be the travelling wave solution to the K-P-P equation with speed $\sqrt{2}$. It follows by \cite[Theorem 2.4 and Theorem 2.6]{KLMR} that under \eqref{condi:extra} and \eqref{llogl}, $w(x)$ is given by
$$w(x)=-\log \pp_{\delta_{0}}\left[\exp\{-C\partial W_{\infty}\mathrm{e}^{-\sqrt{2}x}\}\right]$$
for some constant $C>0$ satisfying that
$\lim_{x\to+\infty}\frac{w(x)}{x\mathrm{e}^{-\sqrt{2}x}}=C.$
Hence by \eqref{repforw(x)} one has
\begin{equation}\label{prop4.1}
-\log\pp_{\delta_{0}}\left[\exp\{-C\partial W_{\infty}\mathrm{e}^{-\sqrt{2}x}\}\right]=1-\pp_{\cdot,\delta_{0}}\left[\exp\{-C\partial M_{\infty}\mathrm{e}^{-\sqrt{2}x}\}\right]
\end{equation}
for all $C>0$ and $x\in\R$. We observe that for every $x\in\R$, the process $((X_{t})_{t\ge 0},\pp_{\delta_{x}})$ is equal in law to $((X_{t}+x)_{t\ge 0},\pp_{\delta_{0}})$. It follows that $(\partial W_{t},\pp_{\delta_{x}})\stackrel{\mbox{d}}{=}
(\mathrm{e}^{\sqrt{2}x}\partial W_{t}-x\mathrm{e}^{\sqrt{2}x}W_{t},\pp_{\delta_{0}})$
for any $t\ge 0$
where $W_{t}:=\langle \mathrm{e}^{\sqrt{2}(\cdot-\sqrt{2}t)},X_{t}\rangle$. Since by \cite[Theorem 2.4(i)]{KLMR} $\lim_{t\to+\infty}W_{t}=0$ $\pp_{\delta_{0}}$-a.s., one gets that
$$(\partial W_{\infty},\pp_{\delta_{x}})\stackrel{\mbox{d}}{=}(\mathrm{e}^{\sqrt{2}x}\partial W_{\infty},\pp_{\delta_{0}}).$$
Using this and the branching property of superprocesses,
one has for any $\lambda>0$,
\begin{eqnarray}
\pp_{\mu}\left[\mathrm{e}^{-\lambda \partial W_{\infty}}\right]&=&\exp\Big\{\int_{\R}\log\pp_{\delta_{x}}\left[\mathrm{e}^{-\lambda \partial W_{\infty}}\right]\mu(dx)\Big\}\nonumber\\
&=&\exp\Big\{\int_{\R}\log\pp_{\delta_{0}}\left[\mathrm{e}^{-\lambda \partial W_{\infty}\mathrm{e}^{\sqrt{2}x}}\right]\mu(dx)\Big\}.\nonumber
\end{eqnarray}
Hence by \eqref{prop4.1} and \eqref{prop2.2} one gets $\pp_{\mu}\left[\mathrm{e}^{-\lambda \partial W_{\infty}}\right]=\pp_{\mu}\left[\mathrm{e}^{-\lambda \partial M_{\infty}}\right]$ for all $\lambda>0$ and so $(\partial W_{\infty},\pp_{\mu})\stackrel{\mbox{d}}{=}(\partial M_{\infty},\pp_{\mu})$.\qed

\subsection{Statement of main results}
In what follows and for the remainder of the paper we assume (A1), (A2) and \eqref{condi0} hold. Additional conditions used are stated explicitly.

Define
\begin{equation}\label{def-H}
\mathcal{H}:=\left\{\phi\in\mathcal{B}^{+}_{b}(\R):\ \int_{0}^{+\infty}y\mathrm{e}^{\sqrt{2}y}\phi(-y)dy<+\infty\right\}.
\end{equation}

\begin{theorem}\label{them1}
Suppose $\phi\in\mathcal{H}$.
Then
\begin{equation}\label{conv-m(t)}
u_{\phi}(t,x-m(t))\to -\log\pp_{\delta_{x}}\left[\exp\{-C(\phi)\partial M_{\infty}\}\right]\mbox{ locally uniformly in $x\in\R$, as }t\to+\infty,
\end{equation}
where $C(\phi)$ is a nonnegative constant given by
\begin{equation}\label{them1.4}
C(\phi)=\lim_{r\to+\infty}\int_{0}^{+\infty}y\mathrm{e}^{\sqrt{2}y}u_{\phi}(r,-y-\sqrt{2}r)dy.
\end{equation}
\end{theorem}

\begin{remark}\rm
\begin{description}
\item{(i)} If additionally \eqref{condi:extra} holds,
then by Proposition \ref{prop4} and  \eqref{conv-m(t)} one has
\begin{equation}\label{conv-m(t)-2}
u_{\phi}(t,x-m(t))\to -\log\pp_{\delta_{x}}\left[\exp\{-C(\phi)\partial W_{\infty}\}\right]\mbox{ locally uniformly in $x\in\R$, as }t\to+\infty,
\end{equation}
This result is obtained independently in \cite[Proposition 1.3(1)]{RSZ}.
In fact, the constant  $C(\phi)$ given in
 \eqref{them1.4} is the same as the one given in
   \cite[Proposition 1.3]{RSZ}.
   This is because our  $u_{\phi}(t,x)$ is the solution to \eqref{eq1} with initial condition $u(0,x)=\phi(x)$, while $U_{\phi}(t,x)$ defined in \cite{RSZ} is the solution to \eqref{eq1} with initial condition $u(0,x)=\phi(-x)$. It holds that $U_{\phi}(t,x)=u_{\phi}(t,-x)$.
 Comparing
 the definitions of  $C(\phi)$ in  \eqref{them1.4}  and the one in
 \cite[Proposition 1.3]{RSZ},
 we see that they are the same.
We also note that
the ``locally uniform convergence"
is slightly stronger than \cite[Proposition 1.3(1)]{RSZ},
where the convergence of \eqref{conv-m(t)-2} is established for each fixed $x\in\R$.
\item{(ii)} Let $\max X_{t}$ denote the supremum of the support of $X_{t}$, i.e., $\max X_{t}:=\inf\{x\in \R:\ X_{t}(x,+\infty)=0\}$. Here we take the convention that $\inf\emptyset=+\infty$.
Unlike for the skeleton BBM, Theorem \ref{them1} does not imply the growth order of $\max X_{t}$ is $m(t)$. In fact, the asymptotic behavior of $\max X_{t}$ depends heavily on the branching mechanism $\psi(\lambda)$ and it may grow much faster than $m(t)$.
We give such examples in Remark \ref{rm4.13}.
\end{description}
\end{remark}

Theorem \ref{them1} yields the existence of the limiting process of the extremal process of super-Brownian motion.

\begin{theorem}\label{them2}
For $t\ge 0$, set
$$\mathcal{E}_{t}:=X_{t}-m(t)\quad \mbox{ and }\quad \mathcal{E}^{Z}_{t}:=Z_{t}-m(t).$$
Then for every $x\in\mathbb{R}$, the process $((\mathcal{E}_{t},\mathcal{E}^{Z}_{t})_{t\ge 0},\pp_{\delta_{x}})$ converges in distribution to a limit $(\mathcal{E}_{\infty},\mathcal{E}^{Z}_{\infty})$ as $t\to +\infty$, where $\mathcal{E}_{\infty}$ is a random Radon measure and $\mathcal{E}^{Z}_{\infty}$ is a random point measure satisfying that
\begin{equation}\label{them2.1}
\mathrm{E}\left[\mathrm{e}^{-\langle f,\mathcal{E}_{\infty}\rangle-\langle g,\mathcal{E}^{Z}_{\infty}\rangle}\right]=\pp_{\delta_{x}}\left[\exp\{-C\left(f+1-\mathrm{e}^{-g}\right)\partial M_{\infty}\}\right]
\end{equation}
for all $f\in\mathcal{H}$  and $g\in\mathcal{B}^{+}(\R)$ with $1-\mathrm{e}^{-g}\in \mathcal{H}$. Moreover, given $\mathcal{E}_{\infty}$, $\mathcal{E}^{Z}_{\infty}$ is a Poisson random measure with intensity $\mathcal{E}_{\infty}(dx)$.
\end{theorem}

In the above statement and for the remainder of this paper,
when we talk about the distributional limit, we do not
specify the probability space where the limit is defined, just  use $\p$ to denote the probability measure, and  $\mathrm{E}$ to denote the  corresponding  expectation.
We remark that
the distribution of $(\mathcal{E}_{\infty},\mathcal{E}^{Z}_{\infty})$
depends on $x$ since it is  the distributional limit of
 $((\mathcal{E}_{t},\mathcal{E}^{Z}_{t})_{t\ge 0},\pp_{\delta_{x}})$.
We will not remark this dependence in similar situation later in this paper.

In the following proposition, we establish a dichotomy on the finiteness of the supremum of the support for the limiting process.
\begin{proposition}\label{prop3}
Suppose $x\in\R$ and $\mathcal{E}_{\infty}$ is the limit of $\left((\mathcal{E}_{t})_{t\ge 0},\pp_{\delta_{x}}\right)$ in distribution.
Let $\max \mathcal{E}_{\infty}$ be the supremum of the support of $\mathcal{E}_{\infty}$. Then $\max\mathcal{E}_{\infty}$
is a.s. finite if and only if
\begin{equation}\nonumber
\sup_{\lambda}C\left(\lambda 1_{(0,+\infty)}\right)<+\infty.
\end{equation}
Otherwise $\mathrm{P}\left(\max\mathcal{E}_{\infty}<+\infty\right)=\pp_{\delta_{x}}\left(\partial M_{\infty}=0\right)
=\mathrm{e}^{-1}
$.
\end{proposition}

We also obtain some new results on the probabilistic representations for the limiting process.
For $u\in Z_{t}$, denote by $I^{(u)}_{s}$ the immigration at time $t+s$ that occurred along the subtree of the skeleton rooted at
$u$ with location $z_{u}(t)$.
Lemma \ref{lem4.4} below shows that under $\pp_{\delta_{x}}$,
for every $s>0$, conditioned on $\{\max Z_{t}-\sqrt{2}t>0\}$,
 $\left(\sum_{u\in Z_{t}}I^{(u)}_{s}-\sqrt{2}s-\max Z_{t},
 Z_{t}-\max Z_{t}\right)$ converges in distribution to a limit
 $\left(\triangle^{I,s},\triangle^{Z}\right)$ as $t\to +\infty$.
Moreover, the limit $\left(\triangle^{I,s},\triangle^{Z}\right)$
does not depend on $x$.

\begin{theorem}\label{them3}
Suppose $x\in\R$ and $(\mathcal{E}_{\infty},\mathcal{E}^{Z}_{\infty})$ is the limit of $((\mathcal{E}_{t},\mathcal{E}^{Z}_{t})_{t\ge 0},\pp_{\delta_{x}})$ in distribution. Then $\mathcal{E}^{Z}_{\infty}$ is a decorated Poisson point process with intensity $c_{*}\partial M_{\infty}\sqrt{2}\mathrm{e}^{-\sqrt{2}y}dy$ and decoration law $\triangle^{Z}$, here $c_{*}=C(1_{(0,+\infty)})$.
Let $\{d_{i}:i\ge 1\}$ be the atoms of $\mathcal{E}^{Z}_{\infty}$, and for every $s>0$, let $\{\triangle^{s}_{i}:i\ge 1\}$ be an independent sequence of i.i.d. random measures with the same law as $(I_{s}-\sqrt{2}s,\pp_{\cdot,\delta_{0}})$, then
\begin{equation}\label{them3.2}
\mathcal{E}_{\infty}\stackrel{\mbox{d}}{=}\lim_{s\to +\infty}\sum_{i\ge 1}\mathcal{T}_{d_{i}}\triangle^{s}_{i}.
\end{equation}
\end{theorem}

The limit in \eqref{them3.2} can not be put into the summation. In fact,
for each $i\ge 1$, $\triangle^{s}_{i}$
converges in distribution to the null measure as $s\to +\infty$,
see Remark \ref{rm5} below.
The following result gives an alternative
description of $\mathcal{E}_{\infty}$.

\begin{proposition}\label{prop5}
Suppose the assumptions of Theorem \ref{them3} hold.
Given $(\partial M_{\infty},\pp_{\delta_{x}})$, let $\{e_{i}:i\ge 1\}$ be the atoms of a Poisson point process with intensity $c_{*}\partial M_{\infty}\sqrt{2}\mathrm{e}^{-\sqrt{2}x}dx$, and
for every $s>0$, let $\{\triangle^{I,s}_{i}:i\ge 1\}$ be an independent sequence of i.i.d random measures with the same law as $\triangle^{I,s}$, then
$$\mathcal{E}_{\infty}\stackrel{\mbox{d}}{=}\lim_{s\to+\infty}\sum_{i\ge 1}\mathcal{T}_{e_{i}}\triangle^{I,s}_{i}.$$
\end{proposition}

For every $s>0$, $\sum_{i\ge 1}\mathcal{T}_{e_{i}}\triangle^{I,s}_{i}$ is a Poisson random measure with exponential intensity, in which each atom is decorated by an independent copy of an auxiliary measure. However, their distributional limit $\mathcal{E}_{\infty}$ may not inherit such a structure. This is revealed by the following theorem.

\begin{theorem}\label{them5}
Suppose $x\in\R$ and $\mathcal{E}_{\infty}$ is the limit of $((\mathcal{E}_{t})_{t\ge 0},\pp_{\delta_{x}})$ in distribution. There exist a constant $\iota\ge 0$ and a measure $\Lambda$ on $\mathcal{M}(\R)\setminus\{0\}$ satisfying that
$$\int_{-\infty}^{+\infty}\mathrm{e}^{-\sqrt{2}x}dx\int_{\mathcal{M}(\R)\setminus\{0\}}\left(1\wedge \mathcal{T}_{x}\mu(A)\right)\Lambda(d\mu)<+\infty,\quad\forall \mbox{ bounded Borel set }A\subset\R,$$
such that
$$\mathcal{E}_{\infty}\stackrel{\mbox{d}}{=}\iota\,\partial M_{\infty}\vartheta+\int_{\mathcal{M}(\R)\setminus\{0\}}\mu\eta(d\mu),$$
where $\vartheta(dx):=\mathrm{e}^{-\sqrt{2}x}dx$ is the (non-random) measure on $\R$ and given $(\partial M_{\infty},\pp_{\delta_{x}})$, $\eta$ is a Poisson random measure on $\mathcal{M}(\R)\setminus\{0\}$ with intensity $c_{*}\partial M_{\infty}\int_{-\infty}^{+\infty}\sqrt{2}\mathrm{e}^{-\sqrt{2}x}\mathcal{T}_{x}\Lambda(d\mu)dx$.
Moreover, $\iota$ and $\Lambda(d\mu)$
satisfy
\eqref{lem5.11.2} and \eqref{lem5.11.3} of Lemma \ref{lem5.11}, respectively.
\end{theorem}

The constant $\iota$ may not be $0$ in general.
The argument of Remark \ref{rm6} shows that $\iota=0$ if the following condition holds:
\medskip

\textbf{(A3)} There are constants $a,b>0$ and $0<\gamma\le 1$ such that
$$\psi(\lambda)\ge -a\lambda+b\lambda^{1+\gamma},\quad \forall \lambda>0.$$

\noindent
When (A3) is satisfied, $\mathcal{E}_{\infty}$ is equal in law to a Poisson random measure on $\mathcal{M}(\R)$ with intensity\\ $c_{*}\partial M_{\infty}\int_{-\infty}^{+\infty}\sqrt{2}\mathrm{e}^{-\sqrt{2}x}\mathcal{T}_{x}\Lambda(d\mu)dx$.
Theorem \ref{them4} below further shows that
$\Lambda(d\mu)=\frac{\tilde{c}_{0}}{c_{*}}\mathrm{P}\left(\widetilde{\triangle}^{X}\in d\mu\right)$, where
$\tilde{c}_{0}$ is a constant given by \eqref{defitildec} with $\phi=0$, and
$\widetilde{\triangle}^{X}$ is the limit of $X_{t}-\max X_{t}$ conditioned on $\{\max X_{t}-\sqrt{2}t>0\}$.
In fact, it is proved in Lemma \ref{lem6.2} that under (A3),
conditioned on $\{\max X_{t}-\sqrt{2}t>0\}$, the random measures
 $(X_{t}-\max X_{t},Z_{t}-\max X_{t})$
converges, as $t\to+\infty$, in distribution to a limit
$(\widetilde{\triangle}^{X},\widetilde{\triangle}^{Z})$, and that given $\widetilde{\triangle}^{X}$, $\widetilde{\triangle}^{Z}$ is a Poisson random measure with intensity $\widetilde{\triangle}^{X}$.

\begin{theorem}\label{them4}
Assume in addition that (A3) holds.
Suppose $x\in \R$ and $(\mathcal{E}_{\infty},\mathcal{E}^{Z}_{\infty})$ is the limit of $((\mathcal{E}_{t},\mathcal{E}^{Z}_{t})_{t\ge 0}$, $\pp_{\delta_{x}})$
in distribution.
Let $\tilde{c}_{0}$ be given by \eqref{defitildec} with $\phi=0$.
Given $(\partial M_{\infty},\pp_{\delta_{x}})$, let $\{\tilde{e}_{i}:i\ge 1\}$ be the atoms of a Poisson point process with intensity $\tilde{c}_{0}\partial M_{\infty}\sqrt{2}\mathrm{e}^{-\sqrt{2}y}dy$ and $\{(\widetilde{\triangle}^{X}_{i},\widetilde{\triangle}^{Z}_{i}):i\ge 1\}$ be an independent sequence of i.i.d. random measures with the same law as $(\widetilde{\triangle}^{X},\widetilde{\triangle}^{Z})$. Then we have
$$\left(\mathcal{E}_{\infty},\mathcal{E}^{Z}_{\infty}\right)\stackrel{\mbox{d}}{=}\left(\sum_{i\ge 1}\mathcal{T}_{\tilde{e}_{i}}\widetilde{\triangle}^{X}_{i},\sum_{i\ge 1}\mathcal{T}_{\tilde{e}_{i}}\widetilde{\triangle}^{Z}_{i}\right).$$
\end{theorem}

\begin{remark}\rm
It is easy to see that condition (A3) implies \eqref{condi:extra}. So by Proposition \ref{prop4} one can replace $\partial M_{\infty}$ by $\partial W_{\infty}$ in the statement of Theorem \ref{them4}.
\end{remark}

\begin{remark}\label{rm7}\rm
Assume that conditions (A1)-(A3) and  \eqref{condi0} hold.
Theorem \ref{them4} implies
that $\mathcal{E}^{Z}_{\infty}$
is a DPPP($\tilde{c}_{0}\partial M_{\infty}\sqrt{2}\mathrm{e}^{-\sqrt{2}y}dy$, $\widetilde{\triangle}^{Z}$),
while Theorem \ref{them3}
says that $\mathcal{E}^{Z}_{\infty}$
is a DPPP($c_*\partial M_{\infty}\sqrt{2}\mathrm{e}^{-\sqrt{2}y}dy$, ${\triangle}^{Z}$).
The two theorems
 give two interpretations of $\mathcal{E}^{Z}_{\infty}$ as a decorated Poisson point process.
Though the two interpretations are equal in law,  they have different intensities
and then different decoration laws.
To see this,
we only need to prove that $c_{*}<\tilde{c}_{0}$.
Using  $\mathbb{P}_{\delta_x}(\partial M_\infty>0)>0$ and Theorem \ref{them3},
one has that
$$\pp_{\delta_{x}}\left(\mathcal{E}_{\infty}(0,+\infty)=0\right)=\pp_{\delta_{x}}\left(\max \mathcal{E}_{\infty}\le 0\right)\le \pp_{\delta_{x}}\left(\max \mathcal{E}^{Z}_{\infty}\le 0\right)=\pp_{\delta_{x}}\left[\exp\{-c_{*}\partial M_{\infty}\}\right]<1.$$
Then we have,
for $\lambda>1$,
\begin{align*}
\mathbb{P}_{\delta_{x}}\left[e^{-C(\lambda 1_{(0,+\infty)})\partial M_\infty}\right]
&=\mathrm{E}\left[e^{-\lambda \mathcal{E}_\infty(0,\infty)}\right]\\
 &<\mathrm{E}\left[e^{- \mathcal{E}_\infty(0,\infty)} \right]
 =\mathbb{P}_{\delta_{x}}\left[e^{-C(1_{(0,+\infty)})\partial M_\infty}\right].
\end{align*}
Using $\mathbb{P}(\partial M_\infty>0)>0$ again,
we get
$c_*=C(1_{(0,+\infty)})<C(\lambda 1_{(0,+\infty)})\le \tilde{c}_0.$
\end{remark}

\section{Convergence of the extremal process of super-Brownian motion}\label{sec3}

\subsection{Proof of Proposition \ref{lem3.5}}\label{sec3.1}

Recall that $((B_{t})_{t\ge0},\Pi_x)$ a standard Brownian motion starting at $x$, and
that $u_{f}(t,x)$ is  the unique nonnegative solution to \eqref{eq1} with initial condition $f$.
\begin{lemma}\label{lem3.1}
Suppose $f,f_{1},f_{2}\in\mathcal{B}^{+}_{b}(\mathbb{R})$, $s,t\ge 0$ and $x,y\in\mathbb{R}$. Then
\begin{description}
\item{(1)} $u_{f}(t,x)\le \mathrm{e}^{t}P_{t}f(x)\le \mathrm{e}^{t}\|f\|_{\infty}.$
\item{(2)} For any $M\ge 1$, $u_{Mf}(t,x)\le M u_{f}(t,x).$
\item{(3)} $u_{f_{1}}(t,x)\vee u_{f_{2}}(t,x)\le u_{f_{1}+f_{2}}(t,x)\le u_{f_{1}}(t,x)+u_{f_{2}}(t,x)$.
\item{(4)} $u_{f}(t+s,x+y)=u_{\mathcal{T}_{y}u_{f}(s,\cdot)}(t,x)$. In particular, $u_{f}(t,x+y)=u_{\mathcal{T}_{y}f}(t,x)$.
\end{description}
\end{lemma}

\proof (1) and (2) follow from Jensen's inequality. In fact, one has
$$u_{f}(t,x)=-\log \pp_{\delta_{x}}\left[\mathrm{e}^{-\langle f,X_{t}\rangle}\right]\le -\log \mathrm{e}^{-\pp_{\delta_{x}}\left[\langle f,X_{t}\rangle\right]}=\pp_{\delta_{x}}\left[\langle f,X_{t}\rangle\right]=\mathrm{e}^{t}P_{t}f(x),$$
and
$$u_{Mf}(t,x)=-\log \pp_{\delta_{x}}\left[\mathrm{e}^{-M\langle f,X_{t}\rangle}\right]\le -\log \left(\pp_{\delta_{x}}\left[\mathrm{e}^{-\langle f,X_{t}\rangle}\right]\right)^{M}=M u_{f}(t,x).$$

(3) follows directly from \cite[Lemma 2.3(2)]{RSZ}.

(4) Fix $s\ge 0$ and $y\in\mathbb{R}$. Let $v(t,x):=u_{f}(t+s,x+y)$ for all $t\ge 0$ and $x\in\mathbb{R}$. It is easy to verify that $v$ is the unique nonnegative solution to \eqref{eq1} with initial condition $v(0,x)=u_{f}(s,x+y)=\mathcal{T}_{y}u_{f}(s,\cdot)(x)$. Hence we get $v(t,x)=u_{\mathcal{T}_{y}u_{f}(s,\cdot)}(t,x)$. In particular by setting $s=0$ we get that $u_{f}(t,x+y)=u_{\mathcal{T}_{y}u_{f}(0,\cdot)}(t,x)=u_{\mathcal{T}_{y}f}(t,x).$ \qed

\bigskip

Our proof of Proposition \ref{lem3.5} follows from two main steps:
the first step is
to establish the convergence for the Laplace functionals of the skeleton BBM which are truncated by a certain cutoff, and
the second step is
to show
the convergence continues to hold when the cutoff is lifted.
We need the following lemmas,
which are refinements of \cite[Proposition 4.4 and Lemma 4.9]{ABK}.
In fact, \cite{ABK} proves the same results
for $[0,1]$-valued functions with support bounded on the left. We extend
their results to all functions of $\mathcal{H}_{1}$.
Though the idea of our proofs is similar to \cite{ABK}, we give the details here for the reader's convenience.

\begin{lemma}\label{lem3.3}
Suppose $\phi\in \mathcal{H}_{1}$. Then for all $r>0$,
$$C_{r}(\phi):=\sqrt{\frac{2}{\pi}}\int_{0}^{+\infty}y\mathrm{e}^{\sqrt{2}y}u_{\phi}(r,-y-\sqrt{2}r)dy$$
exists and is finite. Moreover, the limit
$$C(\phi):=\lim_{r\to+\infty}C_{r}(\phi)$$
exists and is finite, and for every $x\in\mathbb{R}$,
$$\lim_{t\to+\infty}\frac{t^{3/2}}{\frac{3}{2\sqrt{2}}\log t}\mathrm{e}^{-\sqrt{2}x} u_{\phi}(t,x-\sqrt{2}t)=C(\phi).$$
\end{lemma}

\proof
By Lemma \ref{lem3.1}(1),
$$u_{\phi}(t,-x)\le \mathrm{e}^{t}P_{t}\phi(-x)\quad\forall t\ge 0,\ x\in\mathbb{R}.$$
Then
\begin{equation}\label{u-Cr}
C_{r}(\phi)
\le\sqrt{\frac{2}{\pi}}e^r \int_{-\infty}^{+\infty}|y|\mathrm{e}^{\sqrt{2}y}P_r\phi(-y-\sqrt{2}r)dy.
\end{equation}
By Fubini's theorem and change of variables we have
\begin{eqnarray}
\int_{-\infty}^{+\infty}|y|\mathrm{e}^{\sqrt{2}y}P_r\phi(-y-\sqrt{2}r)dy
&= &\frac{1}{\sqrt{2\pi r}}\int_{-\infty}^{+\infty}|y|\mathrm{e}^{\sqrt{2}y}dy\int_{-\infty}^{+\infty}\mathrm{e}^{-\frac{(z-y-\sqrt{2}r)^{2}}{2r}}\phi(-z)dz\nonumber\\
&=&e^{-r}\frac{1}{\sqrt{2\pi r}}\int_{-\infty}^{+\infty}\phi(-z)e^{\sqrt{2}z}dz
\int_{-\infty}^{+\infty}
|y|\mathrm{e}^{-\frac{(z-y)^2}{2r}}dy\nonumber\\
&=&e^{-r}\frac{1}{\sqrt{2\pi }}\int_{-\infty}^{+\infty}\phi(-z)e^{\sqrt{2}z}dz
\int_{-\infty}^{+\infty}
|\sqrt{r}x+z|\mathrm{e}^{-\frac{x^2}{2}}dx\nonumber\\
&\le &e^{-r}\frac{1}{\sqrt{2\pi }}\int_{-\infty}^{+\infty}\phi(-z)e^{\sqrt{2}z}
(\sqrt{r}\Pi_0(|B_1|)+|z|)dz.
\label{lem3.3.1}
\end{eqnarray}
We get by \eqref{u-Cr} and \eqref{lem3.3.1} that
$$C_{r}(\phi)\le c_{1}\int_{-\infty}^{+\infty}\phi(-z)\mathrm{e}^{\sqrt{2}z}\left(|z|+1\right)dz$$
for some constant $c_{1}=c_{1}(r)>0$. The fact that $\phi\in\mathcal{H}_{1}$ implies that the integral in the right hand side is finite. Thus we get that $C_{r}(\phi)<+\infty$.

Let $u(t,x):=u_{\phi}(t,-x)$ for $t\ge 0$ and $x\in\mathbb{R}$. Then $u$ is a solution to \eqref{eq1} with initial condition $u(0,x)=\phi(-x)$ satisfying that
\begin{equation}\nonumber
\int_{0}^{+\infty}y\mathrm{e}^{\sqrt{2}y}u(0,y)dy<+\infty.
\end{equation}
It then follows by \cite[Proposition 4.3]{ABK} that for $r$ large enough, $t\ge 8r$ and $x\ge 8r-\frac{3}{2\sqrt{2}}\log t$,
\begin{equation}\label{lem3.2}
\gamma^{-1}(r)\Psi(r,t,x+\sqrt{2}t)\le u(t,x+\sqrt{2}t)\le \gamma(r)\Psi(r,t,x+\sqrt{2}t),
\end{equation}
where $\gamma(r)\downarrow 1$ as $r\to+\infty$ and
$$\Psi(r,t,x+\sqrt{2}t)=
 \frac{\mathrm{e}^{-\sqrt{2}x}}{\sqrt{2\pi(t-r)}}
\int_{0}^{+\infty}u(r,y+\sqrt{2}r)\mathrm{e}^{\sqrt{2}y}\mathrm{e}^{-\frac{(y-x)^{2}}{2(t-r)}}
\left(1-\mathrm{e}^{-2y\frac{\left(x+\frac{3}{2\sqrt{2}}\log t\right)}{t-r}}\right)dy.$$
We may rewrite $\Psi(r,t,x+\sqrt{2}t)$ as follows:
$$\Psi(r,t,x+\sqrt{2}t)=\mathrm{e}^{-\sqrt{2}x}\frac{x+\frac{3}{2\sqrt{2}}\log t}{(t-r)^{3/2}}\sqrt{\frac{2}{\pi}}\int_{0}^{+\infty}y\mathrm{e}^{\sqrt{2}y}u(r,y+\sqrt{2}r)\mathrm{e}^{-\frac{(y-x)^{2}}{2(t-r)}}G\left(\frac{2y\left(x+\frac{3}{2\sqrt{2}}\log t\right)}{t-r}\right)dy,$$
where $G(z):=\left(1-\mathrm{e}^{-z}\right)/z$. Using the fact that $G(z)\in [0,1]$ for all $z>0$ and $G(z)\sim 1$ as $z\to 0$, we get by the bounded convergence theorem that
\begin{equation*}
\lim_{t\to+\infty}\frac{t^{3/2}}{\frac{3}{2\sqrt{2}}\log t}\mathrm{e}^{\sqrt{2}x}\Psi(r,t,x+\sqrt{2}t)=C_{r}(\phi),\quad\forall x\in\mathbb{R}.
\end{equation*}
Consequently by letting $t\to+\infty$ in \eqref{lem3.2}, we have
\begin{eqnarray*}
0\le \gamma^{-1}(r)C_{r}(\phi)&\le&\liminf_{t\to+\infty}\frac{t^{3/2}}{\frac{3}{2\sqrt{2}}\log t}\mathrm{e}^{\sqrt{2}x}u(t,x+\sqrt{2}t)\\
&\le&\limsup_{t\to+\infty}\frac{t^{3/2}}{\frac{3}{2\sqrt{2}}\log t}\mathrm{e}^{\sqrt{2}x}u(t,x+\sqrt{2}t)\le \gamma(r)C_{r}(\phi)<+\infty.
\end{eqnarray*}
Then by letting $r\to +\infty$, we have
\begin{eqnarray*}
0\le \limsup_{r\to+\infty}C_{r}(\phi)&\le& \liminf_{t\to+\infty}\frac{t^{3/2}}{\frac{3}{2\sqrt{2}}\log t}\mathrm{e}^{\sqrt{2}x}u(t,x+\sqrt{2}t)\\
&\le&\limsup_{t\to+\infty}\frac{t^{3/2}}{\frac{3}{2\sqrt{2}}\log t}\mathrm{e}^{\sqrt{2}x}u(t,x+\sqrt{2}t)\le \liminf_{r\to+\infty}C_{r}(\phi)<+\infty.
\end{eqnarray*}
This implies that the limit $\lim_{r\to +\infty}C_{r}(\phi)$ exists and is finite, and is equal to
$$\lim_{t\to+\infty}\frac{t^{3/2}}{\frac{3}{2\sqrt{2}}\log t}\mathrm{e}^{\sqrt{2}x}u(t,x+\sqrt{2}t).$$
Therefore we complete the proof.\qed

\bigskip

\begin{corollary}\label{cor2}
For all $\phi\in\mathcal{H}_{1}$ and $x\in\mathbb{R}$,
$$C(\mathcal{T}_{x}\phi)=\mathrm{e}^{\sqrt{2}x}C(\phi).$$
Moreover, $\lim_{\delta\to+\infty}C(1_{[\delta,+\infty)})=0$.
\end{corollary}
\proof It follows by Lemma \ref{lem3.3} and Lemma \ref{lem3.1}(4) that
\begin{eqnarray}
C(\mathcal{T}_{x}\phi)&=&\lim_{t\to+\infty}\frac{t^{3/2}}{\frac{3}{2\sqrt{2}}\log t}\mathrm{e}^{-\sqrt{2}y}u_{\mathcal{T}_{x}\phi}(t,y-\sqrt{2}t)\nonumber\\
&=&\lim_{t\to+\infty}\frac{t^{3/2}}{\frac{3}{2\sqrt{2}}\log t}\mathrm{e}^{-\sqrt{2}y}u_{\phi}(t,x+y-\sqrt{2}t)\nonumber\\
&=&\mathrm{e}^{\sqrt{2}x}\lim_{t\to+\infty}\frac{t^{3/2}}{\frac{3}{2\sqrt{2}}\log t}\mathrm{e}^{-\sqrt{2}(x+y)}u_{\phi}(t,x+y-\sqrt{2}t)
=
\mathrm{e}^{\sqrt{2}x}C(\phi).\nonumber
\end{eqnarray}
We observe that $1_{[\delta,+\infty)}(x)=\mathcal{T}_{-\delta}1_{[0,+\infty)}(x)$ for all $\delta,x\in\mathbb{R}$. Thus
$$C(1_{[\delta,+\infty)})=C(\mathcal{T}_{-\delta}1_{[0,+\infty)})=\mathrm{e}^{-\sqrt{2}\delta}C(1_{[0,+\infty)})\to 0\mbox{ as }\delta\to +\infty.$$\qed

For $\phi\in\mathcal{H}_{1}$ and $\delta\in\R$,
put
\begin{equation}\label{def-phidelta}
\phi^\delta(x):=\phi(x)1_{(-\infty,\delta)}(x)+1_{[\delta,+\infty)}(x).
\end{equation}
Note that $\phi^\delta\in\mathcal{H}_1$ and
$$u_{\phi^\delta}(t,x)=1-\pp_{\cdot,\delta_{x}}\left[\prod_{u\in Z_{t}}\left(1-\phi(z_{u}(t))\right)1_{\{z_{u}(t)<\delta\}}\right],\quad\forall t\ge0,\ x\in\mathbb{R}.$$
\begin{lemma}\label{lem3.4}
Suppose $\phi\in\mathcal{H}_{1}$ and $\delta\in\R$.
Then
$$u_{\phi^\delta}(t,x-m(t))\to 1-\pp_{\cdot,\delta_{0}}\left[\exp\{-C(\phi^\delta)\partial M_{\infty}\mathrm{e}^{\sqrt{2}x}\}\right]\mbox{ uniformly in }x\in\mathbb{R},\mbox{ as }t\to+\infty,$$
where $m(t)$ is defined by \eqref{def-m(t)}.
\end{lemma}

\proof
Let $v^{\delta}_{\phi}(t,x):=u_{\phi^\delta}(t,-x)$ for $t\ge 0$ and $x\in\mathbb{R}$. Then $v^{\delta}_{\phi}(t,x)$ is a solution to \eqref{eq1} with initial condition
$v^{\delta}_{\phi}(0,x)=\phi(-x)1_{(-\delta,+\infty)}(x)+1_{(-\infty,-\delta]}(x)$.
Using the fact that $\phi\in\mathcal{H}_{1}$, one can easily verify that $v^{\delta}_{\phi}(0,x)$ satisfies   conditions (8.1) and (1.17) of \cite{Bramson}. Hence by \cite[Theorem 8.3]{Bramson}, one has
$$v^{\delta}_{\phi}(t,x+m(t))\to w(x)\mbox{ uniformly in }x\in\mathbb{R}\mbox{ as }t\to +\infty,$$
where $w$ is the unique (up to translations) travelling wave solution with speed $\sqrt{2}$. It is established in \cite{Kyprianou} that
$$
w(x)=1-\pp_{\cdot,\delta_{0}}\left[\exp\{-C\partial M_{\infty}\mathrm{e}^{-\sqrt{2}x}\}\right]
$$
for some constant $C>0$ which is determined by
$C=\lim_{x\to+\infty}\frac{w(x)}{x\mathrm{e}^{-\sqrt{2}x}}.$
Hence to prove this lemma, it suffices to show that $C=C(\phi^\delta)$, or equivalently,
\begin{equation}\label{lem3.4.1}
\lim_{x\to+\infty}\lim_{t\to+\infty}\frac{v^{\delta}_{\phi}(t,x+m(t))}{x\mathrm{e}^{-\sqrt{2}x}}=C(\phi^\delta).
\end{equation}
By \eqref{lem3.2}, for $r$ large enough and $t,x\ge 8r$, we have the bounds
\begin{equation}\label{lem3.4.2}
\gamma(r)^{-1}\Psi(r,t,x+m(t))\le v^{\delta}_{\phi}(t,x+m(t))\le \gamma(r)\Psi(r,t,x+m(t)),
\end{equation}
where $\gamma(r)\downarrow 1$ as $r\to+\infty$ and
$$\Psi(r,t,x+m(t))=\sqrt{\frac{2}{\pi}}\frac{t^{3/2}}{(t-r)^{3/2}}x\mathrm{e}^{-\sqrt{2}x}\int_{0}^{+\infty}y\mathrm{e}^{\sqrt{2}y}v^{\delta}_{\phi}(r,y+\sqrt{2}r)
\mathrm{e}^{-\frac{\left(y-x+\frac{3}{2\sqrt{2}}\log t\right)^{2}}{2(t-r)}}\frac{1-\mathrm{e}^{-\frac{2xy}{t-r}}}{2xy/(t-r)}dy.$$
It follows by the bounded convergence theorem and Lemma \ref{lem3.3} that
\begin{equation}
\lim_{t\to +\infty}\Psi(r,t,x+m(t))=x\mathrm{e}^{-\sqrt{2}x}C_{r}(\phi^\delta).\nonumber
\end{equation}
This together with \eqref{lem3.4.2} yields that for $r$ large enough and $x\ge 8r$,
\begin{equation}\label{lem3.4.3}
\gamma(r)^{-1}C_{r}(\phi^\delta)\le \lim_{t\to+\infty}\frac{v^{\delta}_{\phi}(t,x+m(t))}{x\mathrm{e}^{-\sqrt{2}x}}\le \gamma(r)C_{r}(\phi^\delta).
\end{equation}
Since $\gamma(r)\to 1$ and $C_{r}(\phi^\delta)\to C(\phi^\delta)$ as $r\to+\infty$, we get \eqref{lem3.4.1} by letting $r\to+\infty$ in \eqref{lem3.4.3}.\qed

\bigskip

\begin{corollary}\label{cor3}
Suppose $\phi\in\mathcal{H}_{1}$. Then
$C(\phi)=\lim_{\delta\to+\infty}C(\phi^\delta)$.
\end{corollary}

\proof We note that for every $\delta\in\mathbb{R}$,
$$\phi(x)\le
\phi^{\delta}(x)
\le \phi(x)+1_{[\delta,+\infty)}(x),\quad\forall x\in\mathbb{R}.$$
Thus by Lemma \ref{lem3.1}(3), we have
\begin{equation}\label{u-phi}
u_{\phi}(t,x)\le u_{\phi^\delta}(t,x)\le u_{\phi}(t,x)+
u_{1_{[\delta,+\infty)}}(t,x),
\quad\forall t\ge0,\ x\in\mathbb{R},
\end{equation}
which implies that $C(\phi)\le C(\phi^\delta)\le C(\phi)+C(1_{[\delta,+\infty)})$. Since $\lim_{\delta\to+\infty}C(1_{[\delta,+\infty)})=0$ by Corollary \ref{cor2}, it follows that
$\lim_{\delta\to+\infty}C(\phi^\delta)=C(\phi)$ .\qed

\bigskip

\noindent\textbf{Proof of Proposition \ref{lem3.5}:} The first part of this proposition follows from Lemma \ref{lem3.3}. We only need to show the second part.
For $c>0$ and $x\in\R$, let
$$w_{c}(x):=1-\pp_{\cdot,\delta_{0}}\left[\exp\{-c\partial M_{\infty}\mathrm{e}^{-\sqrt{2}x}\}\right].$$
By the uniqueness (up to translations) of the travelling wave solution, one has $w_{c}(x)=w_{1}(x-\ln c/\sqrt{2})$ for all $x\in\mathbb{R}$. We need to show that
\begin{equation}\label{lem3.5.1}
u_{\phi}(t,x-m(t))\to
w_{C(\phi)}(-x)\quad
\mbox{locally uniformly in $x\in\R$, as }
t\to+\infty.
\end{equation}
For $\delta\ge 0$,
by \eqref{u-phi} we have that
\begin{equation}\label{lem3.5.2}
u_{\phi}(t,x-m(t))-w_{C(\phi)}(-x)\le \left(u_{\phi^\delta}(t,x-m(t))-w_{C(\phi^\delta)}(-x)\right)+\left(w_{C(\phi^\delta)}(-x)-w_{C(\phi)}(-x)\right)
\end{equation}
and
\begin{eqnarray}\label{lem3.5.3}
u_{\phi}(t,x-m(t))-w_{C(\phi)}(-x)&\ge& \left(u_{\phi^\delta}(t,x-m(t))-w_{C(\phi^\delta)}(-x)\right)-\left(u_{1_{[\delta,+\infty)}}(t,x)-w_{C(1_{[\delta,+\infty)})}(-x)\right)\nonumber\\
&&+\left(w_{C(\phi^\delta)}(-x)-w_{C(\phi)}(-x)\right)-w_{C(1_{[\delta,+\infty)}}(-x).
\end{eqnarray}
We note that
\begin{equation}\label{lem3.5.4}
w_{C(\phi^\delta)}(-x)-w_{C(\phi)}(-x)=w_{1}(-x-\frac{\ln C(\phi^\delta)}{\sqrt{2}})-w_{1}(-x-\frac{\ln C(\phi)}{\sqrt{2}}),
\end{equation}
and
\begin{equation}
w_{C(1_{[\delta,+\infty)})}(-x)=w_{1}(-x-\frac{\ln C(1_{[\delta,+\infty)})}{\sqrt{2}}).\label{lem3.5.4'}
\end{equation}
By Corollary \ref{cor3}, $C(\phi^\delta)\to C(\phi)$ and $C(1_{[\delta,+\infty)})\to 0$ as $\delta\to +\infty$. Then by the continuity of $w_{1}$ we get from \eqref{lem3.5.4} and  \eqref{lem3.5.4'} that
$$w_{C(\phi^\delta)}(-x)-w_{C(\phi)}(-x)\to 0,\mbox{ and }w_{C(1_{[\delta,+\infty)})}(-x)\to 0
\mbox{ locally uniformly in $x\in\R$, as }
\delta\to +\infty.$$
On the other hand, by Lemma \ref{lem3.4} we have for $\delta\ge 0$,
$$u_{\phi^\delta}(t,x-m(t))-w_{C(\phi^\delta)}(-x)\to 0,\mbox{ and }u_{1_{[\delta,+\infty)}}(t,x-m(t))-w_{C(1_{[\delta,+\infty)})}(-x)\to 0 \mbox{ uniformly in $x\in\mathbb{R}$},$$
as $t\to+\infty$. Hence we get \eqref{lem3.5.1} by letting first $t\to +\infty$ and then $\delta\to +\infty$ in both \eqref{lem3.5.2} and \eqref{lem3.5.3}.\qed

\subsection{Proof of Theorem \ref{them1}}

First we introduce notation to refer to the different parts of the skeleton decomposition which will be used later in the computation.
For $t\ge 0$, let $\mathcal{F}_{t}$ denote the $\sigma$-filed generated by $Z$, $X^{*}$ and $I$ up to time $t$. Denote by $I^{*,t}_{s}$ the immigration at time $t+s$ that occurred along the skeleton before time $t$.
For $u\in Z_{t}$,
denote by $I^{(u)}_{s}$ the immigration at time $t+s$ that occurred along the subtree of the skeleton rooted at
$u$ with location $z_{u}(t)$.
We have
\begin{equation}
	X_{s+t}=X^{*}_{s+t}+I^{*,t}_{s}+\sum_{u\in Z_{t}}I^{(u)}_{s}\quad\mbox{ for all }s,t\ge 0.\label{skeletondecomp}
\end{equation}
 It is known (see, e.g., \cite{CRY}) that given $\mathcal{F}_{t}$, $(X^{*}_{s+t}+I^{*,t}_{s})_{s\ge 0}$ is equal in distribution to
	$((X^{*}_{s})_{s\ge 0};\pp_{X_{t}})$
and $I^{(u)}:=(I^{(u)}_{s})_{s\ge 0}$ is equal in distribution to $(I;\pp_{\cdot,\delta_{z_{u}(t)}})$. Moreover, the processes $\{I^{(u)}:\ u\in Z_{t}\}$ are mutually independent and are independent of $X^{*}$.
For $f\in\mathcal{B}^{+}_{b}(\mathbb{R})$, $t\ge 0$ and $x\in\R$,
define
$$u^{*}_{f}(t,x):=-\log \pp_{\delta_{x},\cdot}\left[\mathrm{e}^{-\langle f,X^{*}_{t}\rangle}\right],$$
and
\begin{equation}\label{def-V}
V_{f}(t,x):=\pp_{\cdot,\delta_{x}}\left[\mathrm{e}^{-\langle f,I_{t}\rangle}\right].
\end{equation}
Since $X_{t}=X^{*}_{t}+\sum_{u\in Z_{0}}I^{(u)}_{t}$, we have
\begin{eqnarray}
u_{f}(t,x)&=&-\log\pp_{\delta_{x}}\left[\mathrm{e}^{-\langle f,X_{t}\rangle}\right]
=
-\log\pp_{\delta_{x}}\left[\mathrm{e}^{-\langle f,X^{*}_{t}\rangle}\right]-\log\pp_{\delta_{x}}\left[\prod_{u\in Z_{0}}\mathrm{e}^{-\langle f,I^{(u)}_{t}\rangle}\right]\nonumber\\
&=&u^{*}_{f}(t,x)-\log\pp_{\delta_{x}}\left[\prod_{u\in Z_{0}}V_{f}(t,z_{u}(0))\right]
=
u^{*}_{f}(t,x)-\log\pp_{\delta_{x}}\left[\mathrm{e}^{\langle \ln V_{f}(t,\cdot),Z_{0}\rangle}\right].\nonumber
\end{eqnarray}
Using the fact that $(Z_{0},\pp_{\delta_{x}})$ is a Poisson random measure with intensity $\delta_{x}(dy)$, one has
\begin{equation}
u_{f}(t,x)=u^{*}_{f}(t,x)+1-V_{f}(t,x).\label{decomforu1}
\end{equation}
In this section we will make extensive use of \eqref{decomforu1}, mostly when we deal with $u_{f}(t,x-\sqrt{2}t)$ for large $t$,
in which case, $u^{*}_{f}(t,x-\sqrt{2}t)$ becomes easy to handle.

Recall the definition of $\phi^{\delta}$ given in \eqref{def-phidelta}.
The following lemma gives an upperbound for the constant $C(\phi^{0})$ which will be used later.
\begin{lemma}\label{lem3.60}
There exists a constant $c>0$ such that
for any $\phi\in\mathcal{H}_{1}$,
\begin{equation}\label{lem3.60.1}
C\left(\phi^0\right)\le c\left(\int_{1}^{+\infty}x\mathrm{e}^{\sqrt{2}x}\phi(-x)dx+1\right).
\end{equation}
\end{lemma}

\proof
Let $v(t,x):=u_{\phi^0}(t,-x)$ for $t\ge 0$ and $x\in\R$.
It follows by \eqref{repforw(x)}, \eqref{asymforw} and Proposition \ref{lem3.5}
that
\begin{equation}
C\left(\phi^0\right)=\lim_{x\to+\infty}\lim_{t\to+\infty}\frac{v(t,x+m(t))}{x\mathrm{e}^{-\sqrt{2}x}}.\label{lem3.60.3}
\end{equation}
Let $k(s,y):=-\psi(v(s,y))/v(s,y)$ for $s\ge 0$ and $y\in\R$. By the Feyman-Kac formula
\begin{eqnarray}
v(t,x)&=&
\Pi_x\left[\mathrm{e}^{\int_{0}^{t}k(t-s,B_{s})ds}v(0,B_{t})\right]\nonumber\\
&=&\int_{-\infty}^{+\infty}v(0,y)\frac{1}{\sqrt{2\pi t}}\mathrm{e}^{-\frac{(y-x)^{2}}{2t}}\mathrm{E}\left[\mathrm{e}^{\int_{0}^{t}k(t-s,\zeta^{t}_{x,y}(s))ds}\right]dy.\label{lem3.60.2}
\end{eqnarray}
Here
$\{\zeta^{(t)}_{x,y}(s):0\le s\le t\}$ denotes a Brownian bridge of length $t$ starting at $x$ and ending at $y$.
Let $v^{H}(t,x)$ be the solution to
\eqref{eq1}
with heaviside initial condition $v^{H}(0,x)=1_{\{x\le 0\}}$ and let $k^{H}(s,y):=-\psi(v^{H}(s,y))/v^{H}(s,y)$ for $s\ge 0$ and $y\in\R$. Define $m^{H}_{1/2}(t):=\sup\{x\in \R:\ v^{H}(t,x)\ge 1/2\}$ for $t\ge 0$. By \cite[Proposition 8.1 and Proposition 8.2]{Bramson}, there are constants $C^{H}_{1}$ and $C^{H}_{2}$ ($C^{H}_{1}<C^{H}_{2}$) such that for $t$ large enough,
\begin{equation}\label{lem3.60.4}
m(t)+C^{H}_{1}\le m^{H}_{1/2}(t)\le m(t)+C^{H}_{2}.
\end{equation}
Moreover, it is established in \cite[equation (8.21)]{Bramson} that there are constants $C^{H}_{3}>0$ and $r>>1$ such that for all $t\ge 3r$, $x\ge m^{H}_{1/2}(t)+1$ and all $y\in \R$,
$$\mathrm{E}\left[\mathrm{e}^{\int_{0}^{t}k^{H}(t-s,\zeta^{t}_{x,y}(s))ds}\right]\le 2C^{H}_{3}r\mathrm{e}^{t}\left(1-\mathrm{e}^{-\frac{2\hat{y}\bar{z}}{t}}\right)$$
where $\hat{y}:=y\vee 1$ and $\bar{z}:=x-(m(t)+C^{H}_{1})$. Since $v(0,y)=1$ for $y\in (-\infty,0]$ and so $v(s,y)\ge v^{H}(s,y)$ for all $s\ge 0$ and $y\in\R$, it follows by the convexity of $\psi$ that $k(s,y)\le k^{H}(s,y)$ for all $s\ge 0$ and $y\in\R$. Hence one has
$$\mathrm{E}\left[\mathrm{e}^{\int_{0}^{t}k(t-s,\zeta^{t}_{x,y}(s))ds}\right]\le C_{4}\mathrm{e}^{t}\left(1-\mathrm{e}^{-\frac{2\hat{y}\bar{z}}{t}}\right)$$
for $C_{4}=2C^{H}_{3}r$.
Putting this back in \eqref{lem3.60.2} one gets that for $t\ge 3r$ and $x\ge m^{H}_{1/2}(t)+1$,
\begin{eqnarray}
v(t,x)&\le&C_{4}\mathrm{e}^{t}\int_{-\infty}^{+\infty}v(0,y)\frac{1}{\sqrt{2\pi t}}\mathrm{e}^{-\frac{(x-y)^{2}}{2t}}\left(1-\mathrm{e}^{-\frac{2\hat{y}\bar{z}}{t}}\right)dy\nonumber\\
&\le&\frac{C_{4}\mathrm{e}^{t}}{\sqrt{2\pi t}}\left[\int_{1}^{+\infty}\phi(-y)\mathrm{e}^{-\frac{(x-y)^{2}}{2t}}\left(1-\mathrm{e}^{-\frac{2y\bar{z}}{t}}\right)dy
+\left(1-\mathrm{e}^{-\frac{2\bar{z}}{t}}\right)\int_{-\infty}^{1}\mathrm{e}^{-\frac{(x-y)^{2}}{2t}}dy\right]\nonumber\\
&\le&C_{4}\sqrt{\frac{2}{\pi}}\frac{\mathrm{e}^{t}}{t^{3/2}}\bar{z}\left[\int_{1}^{+\infty}\phi(-y)y\mathrm{e}^{-\frac{(x-y)^{2}}{2t}}dy+\int_{-\infty}^{1}\mathrm{e}^{-\frac{(x-y)^{2}}{2t}}dy\right].\nonumber
\end{eqnarray}
The last inequality is from the fact that $1-\mathrm{e}^{-x}\le x$ for all $x\ge 0$. This together with \eqref{lem3.60.4} yields that for $t\ge 3r$ and $x\ge C^{H}_{2}+1$,
\begin{eqnarray}
v(t,x+m(t))&\le&C_{4}\sqrt{\frac{2}{\pi}}\frac{\mathrm{e}^{t}}{t^{3/2}}(x-C^{H}_{1})\left[\int_{1}^{+\infty}\phi(-y)y\mathrm{e}^{-\frac{(x+m(t)-y)^{2}}{2t}}dy
+\int_{-\infty}^{1}\mathrm{e}^{-\frac{(x+m(t)-y)^{2}}{2t}}dy\right]\nonumber\\
&=&C_{4}\sqrt{\frac{2}{\pi}}(x-C^{H}_{1})\mathrm{e}^{-\left(\sqrt{2}-\frac{3}{2\sqrt{2}}\frac{\log t}{t}\right)x-\frac{9}{16}\frac{\log^{2} t}{t}}
\Big[\int_{1}^{+\infty}\phi(-y)y\mathrm{e}^{-\frac{(x-y)^{2}}{2t}}\,\mathrm{e}^{\left(\sqrt{2}-\frac{3}{2\sqrt{2}}\frac{\log t}{t}\right)y}dy\nonumber\\
&&\quad+\int_{-\infty}^{1}\mathrm{e}^{-\frac{(x-y)^{2}}{2t}}\,\mathrm{e}^{\left(\sqrt{2}-\frac{3}{2\sqrt{2}}\frac{\log t}{t}\right)y}dy\Big].\nonumber
\end{eqnarray}
By letting $t\to+\infty$, we have
$$\lim_{t\to+\infty}v(t,x+m(t))\le C_{4}\sqrt{\frac{2}{\pi}}(x-C^{H}_{1})\mathrm{e}^{-\sqrt{2}x}\left[\int_{1}^{+\infty}y\mathrm{e}^{\sqrt{2}y}\phi(-y)dy+\frac{1}{\sqrt{2}}\mathrm{e}^{\sqrt{2}}\right]$$
for $x\ge C^{H}_{2}+1$.
Putting this back in \eqref{lem3.6.3}, one gets that
$$C\left(\phi^0\right)\le C_{4}\sqrt{\frac{2}{\pi}}\left[\int_{1}^{+\infty}y\mathrm{e}^{\sqrt{2}y}\phi(-y)dy+\frac{1}{\sqrt{2}}\mathrm{e}^{\sqrt{2}}\right].$$
Hence we complete the proof.

\qed

\begin{lemma}\label{lem3.6}
Suppose $\{\phi_{s}(x):s\ge 0\}$ is a sequence of functions in $\mathcal{H}_{1}$. If $\phi_{s}(x)\to 0$ as $s\to +\infty$ for all $x\in\R$, and $\int_{-\infty}^{+\infty}|x|\mathrm{e}^{\sqrt{2}x}\phi_{s}(-x)dx\to 0$ as $s\to +\infty$, then $\lim_{s\to+\infty}C(\phi_{s})=0$.
\end{lemma}

\proof Suppose $\alpha(s)\ge 0$ for all $s>0$. (The explicit value of $\alpha(s)$ will be given later.) By Lemma \ref{lem3.1}(3), one has
$$C(\phi_{s})\le C\left(\phi_{s} 1_{(-\infty,\alpha(s))}+1_{[\alpha(s),+\infty)}\right)
,\ \forall s\ge 0.
$$
So it suffices to prove that
\begin{equation}
\lim_{s\to+\infty}C\left(\phi_{s} 1_{(-\infty,\alpha(s))}+1_{[\alpha(s),+\infty)}\right)=0.\label{lem3.6.1}
\end{equation}
We note that $\phi_{s}(x)1_{(-\infty,\alpha(s))}(x)+1_{[\alpha(s),+\infty)}(x)=\mathcal{T}_{-\alpha(s)}\left(\mathcal{T}_{\alpha(s)}\phi_{s}\cdot 1_{(-\infty,0)}+1_{[0,+\infty)}\right)(x)$ for all $x\in\R$. By Corollary \ref{cor2} and Lemma \ref{lem3.60}, we have
\begin{eqnarray}
C\left(\phi_{s}1_{(-\infty,\alpha(s))}+1_{[\alpha(s),+\infty)}\right)&=&\mathrm{e}^{-\sqrt{2}\alpha(s)}C\left(\mathcal{T}_{\alpha_{s}}\phi_{s}\cdot 1_{(-\infty,0)}+1_{[0,+\infty)}\right)\nonumber\\
&\le&\mathrm{e}^{-\sqrt{2}\alpha(s)}c\left[\int_{1}^{+\infty}x\mathrm{e}^{\sqrt{2}x}\mathcal{T}_{\alpha(s)}\phi_{s}(-x)dx+1\right]\nonumber\\
&=&\mathrm{e}^{-\sqrt{2}\alpha(s)}c\left[\int_{1-\alpha(s)}^{+\infty}(y+\alpha(s))\mathrm{e}^{\sqrt{2}(y+\alpha(s))}\phi_{s}(-y)dy+1\right]\nonumber\\
&\le&c\int_{-\infty}^{+\infty}(|y|+\alpha(s))\mathrm{e}^{\sqrt{2}y}\phi_{s}(-y)dy+c\mathrm{e}^{-\sqrt{2}\alpha(s)}.\label{lem3.6.2}
\end{eqnarray}
Let $\alpha(s):=s\wedge\left(\int_{-\infty}^{+\infty}\phi_{s}(-y)\mathrm{e}^{\sqrt{2}y}dy\right)^{-1/2}$.
Then the right hand side of \eqref{lem3.6.2} is no larger than
\begin{equation}\label{lem3.6.3}
c\left[\int_{-\infty}^{+\infty}|y|\mathrm{e}^{\sqrt{2}y}\phi_{s}(-y)dy+\alpha(s)^{-1}+\mathrm{e}^{-\sqrt{2}\alpha(s)}\right].
\end{equation}
Note that
\begin{eqnarray*}
\int_{-\infty}^{+\infty}\phi_{s}(-y)\mathrm{e}^{\sqrt{2}y}dy
&\le&
\int_{|y|\le 1}\phi_{s}(-y)\mathrm{e}^{\sqrt{2}y}dy+\int_{|y|>1}\phi_{s}(-y)|y|\mathrm{e}^{\sqrt{2}y}dy
<+\infty.
\end{eqnarray*}
 Both integrals in the right hand side converge to $0$ as $s\to +\infty$ given that $\phi_{s}(x)\to 0$ as $s\to\infty$ for all $x\in \R$ and $\int_{-\infty}^{+\infty}\phi_{s}(-y)|y|\mathrm{e}^{\sqrt{2}y}dy\to 0$ as $s\to\infty$.
 This implies  $\alpha(s)\to +\infty$ as $s\to +\infty$.
 Thus \eqref{lem3.6.3} converges to $0$ as $s\to +\infty$ and so we prove \eqref{lem3.6.1}.\qed

\bigskip

We recall the definition of $\mathcal{H}$ given in \eqref{def-H}.
\begin{lemma}\label{cor1}
For
all $\phi\in\mathcal{H}$ and $s\ge 0$, $u_{\phi}(s,\cdot-\sqrt{2}s)$,  $u^{*}_{\phi}(s,\cdot-\sqrt{2}s)$ and $1-V_{\phi}(s,\cdot-\sqrt{2}s)\in \mathcal{H}$.
\end{lemma}
\proof Fix an arbitrary $\phi\in\mathcal{H}$ and $s\ge 0$. Since $0\le u^{*}_{\phi}(s,x-\sqrt{2}s),\ 1-V_{\phi}(s,x-\sqrt{2}s)\le u_{\phi}(s,x-\sqrt{2}s)$ for all $x\in\mathbb{R}$, it suffices to prove that
\begin{equation}\label{cor1.1}
u_{\phi}(s,\cdot-\sqrt{2}s)\in\mathcal{H}.
\end{equation}
Let $M:=\|\phi\|_{\infty}$. If $M\le 1$, then $\phi\in\mathcal{H}_{1}$ and \eqref{cor1.1} follows directly from the first conclusion of Lemma \ref{lem3.3}. Now suppose $M>1$. Let $\phi_{1}=\phi/M$. Then by Lemma \ref{lem3.1}(2), one has
$u_{\phi}(s,x-\sqrt{2}s)\le M u_{\phi_{1}}(s,x-\sqrt{2}s)$ for all $x\in\mathbb{R}$, where $u_{\phi_{1}}(s,\cdot-\sqrt{2}s)\in \mathcal{H}$ by Lemma \ref{lem3.3}. This implies that $u_{\phi}(s,\cdot-\sqrt{2}s)\in\mathcal{H}$.\qed

\bigskip

\begin{lemma}\label{lem3.7}
Suppose $\phi\in\mathcal{H}$. Then
$u^{*}_{\phi}(s,\cdot-\sqrt{2}s)\in \mathcal{H}_{1}$
for $s$ large enough. Moreover,
$$\lim_{s\to+\infty}C(u^{*}_{\phi}(s,\cdot-\sqrt{2}s))=0.$$
\end{lemma}

\proof It follows by Jensen's inequality that
$$u^{*}_{\phi}(s,x)=-\log\pp_{\delta_{x}}\left[\mathrm{e}^{-\langle \phi,X^{*}_{t}\rangle}\right]\le \pp_{\delta_{x}}\left[\langle \phi,X^{*}_{t}\rangle\right]=\mathrm{e}^{\alpha^{*}s}P_{s}\phi(x)\le \mathrm{e}^{\alpha^{*}s}\|\phi\|_{\infty}.$$
Since $\alpha^{*}=-\psi'(1)<0$, we have $u^{*}_{\phi}(s,x)\to 0$ as $s\to +\infty$ for all $x\in\R$ and $\|u^{*}_{\phi}(s,\cdot-\sqrt{2}s)\|_{\infty}\le 1$ for $s$ large enough.
Hence $u^{*}_{\phi}(s,\cdot-\sqrt{2}s)\in \mathcal{H}_{1}$
by Lemma \ref{cor1}.
 We note that
$$u^{*}_{\phi}(s,-x-\sqrt{2}s)\le \mathrm{e}^{\alpha^{*}s}P_{s}\phi(-x-\sqrt{2}s).$$
Thus by \eqref{lem3.3.1} we have
\begin{eqnarray}
\int_{-\infty}^{+\infty}\mathrm{e}^{\sqrt{2}x}|x|u^{*}_{\phi}(s,-x-\sqrt{2}s)dx
&\le&\mathrm{e}^{\alpha^{*}s}\int_{-\infty}^{+\infty}\mathrm{e}^{\sqrt{2}x}|x|P_{s}\phi(-x-\sqrt{2}s)dx\nonumber\\
&\le&\mathrm{e}^{(\alpha^{*}-1)s}\int_{-\infty}^{+\infty}\phi(-y)\mathrm{e}^{\sqrt{2}y}\left( |y|+\sqrt{s}
\Pi_0\left[|B_{1}|\right]\right)dy.\nonumber
\end{eqnarray}
The assumption that $\phi\in\mathcal{H}$ implies that
$\int_{-\infty}^{+\infty}\phi(-y)\mathrm{e}^{\sqrt{2}y}|y|dy<\infty.$
Since $\alpha^{*}<0$, we get by the above inequality that
$$\int_{-\infty}^{+\infty}\mathrm{e}^{\sqrt{2}x}|x|u^{*}_{\phi}(s,-x-\sqrt{2}s)dx\to 0 \mbox{ as } s\to +\infty,$$ and thus by Lemma \ref{lem3.6} $C(u^{*}_{\phi}(s,\cdot-\sqrt{2}s))\to 0$ as $s\to +\infty$.\qed

\bigskip

The following lemma extends the result of Lemma \ref{lem3.3} to all functions of $\mathcal{H}$.

\begin{lemma}\label{lem3.8}
Suppose $\phi\in\mathcal{H}$. Then for any $r>0$,
$$C_{r}(\phi):=\sqrt{\frac{2}{\pi}}\int_{0}^{+\infty}y\mathrm{e}^{\sqrt{2}y}u_{\phi}(r,-y-\sqrt{2}r)dy$$
exists and is finite. The limit
$$C(\phi):=\lim_{r\to+\infty}C_{r}(\phi)$$
exists and is finite.
Moreover, for every $x\in\mathbb{R}$,
\begin{equation}\label{lem3.8.1}
C(\phi)=\lim_{t\to+\infty}\frac{t^{3/2}}{\frac{3}{2\sqrt{2}}\log t}\mathrm{e}^{-\sqrt{2}x}u_{\phi}(t,x-\sqrt{2}t)=\lim_{s\to+\infty}C\left(1-V_{\phi}(s,\cdot-\sqrt{2}s)\right).
\end{equation}
\end{lemma}

\proof Without loss of generality we assume $\phi\in \mathcal{H}\setminus \mathcal{H}_{1}$. Let $M:=\|\phi\|_{\infty}$ and $\phi_{1}=\phi/M$. Then $\phi_{1}\in\mathcal{H}_{1}$. The finiteness of $C_{r}(\phi)$ is immediate since by Lemma \ref{lem3.1}(2) $u_{\phi}(r,-y-\sqrt{2}r)\le M u_{\phi_{1}}(r,-y-\sqrt{2}r)$ for all $r\ge 0$ and $y\in\mathbb{R}$ and so $C_{r}(\phi)\le M C_{r}(\phi_{1})<+\infty$. Since
\begin{equation}\nonumber
u_{\phi}(t,x)=u^{*}_{\phi}(t,x)+1-V_{\phi}(t,x),\quad\forall t\ge 0,\ x\in\mathbb{R},
\end{equation}
we get by Lemma \ref{lem3.1}(3)(4) that
\begin{eqnarray}
u_{1-V_{\phi}(s,\cdot-\sqrt{2}s)}(r,x-\sqrt{2}r)&\le&u_{u_{\phi}(s,\cdot-\sqrt{2}s)}(r,x-\sqrt{2}r)\nonumber\\
&=&u_{\phi}(s+r,x-\sqrt{2}(s+r))\nonumber\\
&\le&u_{u^{*}_{\phi}(s,\cdot-\sqrt{2}s)}(r,x-\sqrt{2}r)+u_{1-V_{\phi}(s,\cdot-\sqrt{2}s)}(r,x-\sqrt{2}r).\label{lem3.8.5}
\end{eqnarray}
This implies that
\begin{equation}\label{lem3.8.3}
C_{r}\left(1-V_{\phi}(s,\cdot-\sqrt{2}s)\right)\le C_{r+s}\left(\phi\right)\le C_{r}\left(u^{*}_{\phi}(s,\cdot-\sqrt{2}s)\right)+C_{r}\left(1-V_{\phi}(s,\cdot-\sqrt{2}s)\right)
\end{equation}
for all $r>0$. Since $1-V_{\phi}(s,\cdot-\sqrt{2}s)\in \mathcal{H}_{1}$ and $u^{*}_{\phi}(s,\cdot-\sqrt{2}s)\in \mathcal{H}_{1}$ for $s$ large enough, we get by \eqref{lem3.8.3} and  Proposition \ref{lem3.5} that
\begin{equation}\nonumber
C\left(1-V_{\phi}(s,\cdot-\sqrt{2}s)\right)\le \liminf_{r\to+\infty}C_{r}(\phi)\le \limsup_{r\to+\infty}C_{r}(\phi)\le C\left(u^{*}_{\phi}(s,\cdot-\sqrt{2}s)\right)+C\left(1-V_{\phi}(s,\cdot-\sqrt{2}s)\right).
\end{equation}
Since $\lim_{s\to+\infty}C(u^{*}_{\phi}(s,\cdot-\sqrt{2}s))=0$, we have
\begin{equation}\label{lem3.8.2}
\limsup_{s\to+\infty}C\left(1-V_{\phi}(s,\cdot-\sqrt{2}s)\right)\le \liminf_{r\to+\infty}C_{r}(\phi)\le \limsup_{r\to+\infty}C_{r}(\phi)\le \liminf_{s\to+\infty}C\left(1-V_{\phi}(s,\cdot-\sqrt{2}s)\right)
\end{equation}
Since the $C_{r}(\phi)\le MC_{r}(\phi_{1})$ for all $r>0$ and the latter is bounded in $r$, \eqref{lem3.8.2} implies that the limit $C(\phi)=\lim_{r\to+\infty}C_{r}(\phi)$ exists and is finite, and satisfies that
\begin{equation}\label{lem3.8.4}
C(\phi)=\lim_{s\to+\infty}C\left(1-V_{\phi}(s,\cdot-\sqrt{2}s)\right).
\end{equation}
On the other hand,
it follows from Lemma \ref{lem3.3} and \eqref{lem3.8.5} that
\begin{eqnarray}
C\left(1-V_{\phi}(s,\cdot-\sqrt{2}s)\right)&=&\lim_{t\to+\infty}\frac{t^{3/2}}{\frac{3}{2\sqrt{2}}\log t}\mathrm{e}^{-\sqrt{2}x}u_{1-V_{\phi}(s,\cdot-\sqrt{2}s)}(t,x-\sqrt{2}t)\nonumber\\
&\le&\liminf_{t\to+\infty}\frac{t^{3/2}}{\frac{3}{2\sqrt{2}}\log t}\mathrm{e}^{-\sqrt{2}x}u_{\phi}(t,x-\sqrt{2}t)\nonumber\\
&\le&\limsup_{t\to+\infty}\frac{t^{3/2}}{\frac{3}{2\sqrt{2}}\log t}\mathrm{e}^{-\sqrt{2}x}u_{\phi}(t,x-\sqrt{2}t)\nonumber\\
&\le&\lim_{t\to+\infty}\frac{t^{3/2}}{\frac{3}{2\sqrt{2}}\log t}\mathrm{e}^{-\sqrt{2}x}u_{u^{*}_{\phi}(s,\cdot-\sqrt{2}s)}(t,x-\sqrt{2}t)\nonumber\\
&&\quad+\lim_{t\to+\infty}\frac{t^{3/2}}{\frac{3}{2\sqrt{2}}\log t}\mathrm{e}^{-\sqrt{2}x}u_{1-V_{\phi}(s,\cdot-\sqrt{2}s)}(t,x-\sqrt{2}t)\nonumber\\
&=&C\left(u^{*}_{\phi}(s,\cdot-\sqrt{2}s)\right)+C\left(1-V_{\phi}(s,\cdot-\sqrt{2}s)\right).\nonumber
\end{eqnarray}
Letting $s\to +\infty$, we get by \eqref{lem3.8.4} and Lemma \ref{lem3.7} that
$$C(\phi)=\lim_{t\to+\infty}\frac{t^{3/2}}{\frac{3}{2\sqrt{2}}\log t}\mathrm{e}^{-\sqrt{2}x}u_{\phi}(t,x-\sqrt{2}t).$$\qed

\begin{corollary}\label{cor4}
For $f,f_{1},f_{2}\in \mathcal{H}$ and $M\ge 1$,
\begin{eqnarray*}
&&C(f)=C(u_{f}(s,\cdot-\sqrt{2}s))
,\ \forall s\ge 0,\\
&&C(Mf)\le MC(f),\\
&&C(f_{1})\vee C(f_{2})\le C(f_{1}+f_{2})\le C(f_{1})+C(f_{2}).
\end{eqnarray*}
\end{corollary}
\proof This result follows directly from Lemma \ref{lem3.8} and Lemma \ref{lem3.1}.\qed

\bigskip

\noindent\textbf{Proof of Theorem \ref{them1}:}
We first suppose that $\phi\in\mathcal{H}_{1}$. Then by Proposition \ref{lem3.5} we have
\begin{equation}\label{them1.1}
u_{\phi}(t,x-m(t))\to 1-\pp_{\cdot,\delta_{0}}\left[\exp\{-C(\phi)\mathrm{e}^{\sqrt{2}x}\partial M_{\infty}\}\right]\mbox{ locally uniformly in $x\in\R$, as }t\to+\infty.
\end{equation}
By \eqref{prop2.2} we have
\begin{equation}
\pp_{\delta_{x}}\left[\exp\{-C(\phi)\partial M_{\infty}\}\right]
=\exp\left\{-\left(1-\pp_{\cdot,\delta_{0}}\left[\exp\{-C(\phi)\mathrm{e}^{\sqrt{2}x}\partial M_{\infty}\}\right]\right)\right\}.\label{them1.0}
\end{equation}
Putting this back in \eqref{them1.1}, we get that
$$u_{\phi}(t,x-m(t))\to -\log\pp_{\delta_{x}}\left[\exp\{-C(\phi)\partial M_{\infty}\}\right]\mbox{ locally uniformly in $x\in\R$, as }t\to+\infty.$$

Now we suppose $\phi\in\mathcal{H}\setminus\mathcal{H}_{1}$. For any $r,s>0$ and $x\in\mathbb{R}$, one can rewrite $u_{\phi}(r+s,x-m(r+s))$ as
\begin{equation}
u_{u_{\phi}(s,\cdot-\sqrt{2}s)}(r,x-m(r)+O_{s}(r)),\nonumber
\end{equation}
where $O_{s}(r):=\frac{3}{2\sqrt{2}}\log \frac{r+s}{r}$.  Note that  for $s>0$, $O_{s}(r)\to 0$ as $r\to+\infty$. In view of this and \eqref{decomforu1} we get by Lemma \ref{lem3.1}(4) that
\begin{eqnarray}
&&u_{1-V_{\phi}(s,\cdot-\sqrt{2}s)}(r,x-m(r)+O_{s}(r))\nonumber\\
&\le& u_{\phi}(r+s,x-m(r+s))\nonumber\\
&\le& u_{1-V_{\phi}(s,\cdot-\sqrt{2}s)}(r,x-m(r)+O_{s}(r))+u_{u^{*}_{\phi}(s,\cdot-\sqrt{2}s)}(r,x-m(r)+O_{s}(r)).\label{them1.2}
\end{eqnarray}
It follows that
\begin{eqnarray}
&&\left|u_{\phi}(r+s,x-m(r+s))-\left(-\log \pp_{\delta_{x}}\left[\exp\{-C(\phi)\partial M_{\infty}\}\right]\right)\right|\nonumber\\
&\le&\left|u_{1-V_{\phi}(s,\cdot-\sqrt{2}s)}(r,x-m(r)+O_{s}(r))-\left(-\log \pp_{\delta_{x}}\left[\exp\{-C(1-V_{\phi}(s,\cdot-\sqrt{2}s))\partial M_{\infty}\}\right]\right)\right|\nonumber\\
&&+\left|u_{u^{*}_{\phi}(s,\cdot-\sqrt{2}s)}(r,x-m(r)+O_{s}(r))-\left(-\log \pp_{\delta_{x}}\left[\exp\{-C(u^{*}_{\phi}(s,\cdot-\sqrt{2}s))\partial M_{\infty}\}\right]\right)\right|\nonumber\\
&&+\left|\log \pp_{\delta_{x}}\left[\exp\{-C(\phi)\partial M_{\infty}\}\right]-\log \pp_{\delta_{x}}\left[\exp\{-C(1-V_{\phi}(s,\cdot-\sqrt{2}s))\partial M_{\infty}\}\right]\right|\nonumber\\
&&+\left|\log \pp_{\delta_{x}}\left[\exp\{-C(u^{*}_{\phi}(s,\cdot-\sqrt{2}s))\partial M_{\infty}\}\right]\right|.\label{them1.3}
\end{eqnarray}
We have proved in the first part that the first two terms of \eqref{them1.3} converge to $0$ locally uniformly in $x\in\R$ as $r\to+\infty$. On the other hand we have by \eqref{them1.0}
\begin{eqnarray}
&&\left|\log \pp_{\delta_{x}}\left[\exp\{-C(\phi)\partial M_{\infty}\}\right]-\log \pp_{\delta_{x}}\left[\exp\{-C(1-V_{\phi}(s,\cdot-\sqrt{2}s))\partial M_{\infty}\}\right]\right|\nonumber\\
&=&\left|\pp_{\cdot,\delta_{0}}\left[\exp\{-C(\phi)\mathrm{e}^{\sqrt{2}x}\partial M_{\infty}\}\right]-\pp_{\cdot,\delta_{0}}\left[\exp\{-C(1-V_{\phi}(s,\cdot-\sqrt{2}s))\mathrm{e}^{\sqrt{2}x}\partial M_{\infty}\}\right]\right|\nonumber\\
&\le&\left|C(\phi)-C(1-V_{\phi}(s,\cdot-\sqrt{2}s))\right|\mathrm{e}^{\sqrt{2}x}\pp_{\cdot,\delta_{0}}\left[\partial M_{\infty}\right],\label{4.27}
\end{eqnarray}
where in the lase inequality we use the fact that $|e^{-x_1}-e^{-x_2}|\le |x_1-x_2|$, for all $x_1,x_2\ge0$.
Since $C(\phi)-C(1-V_{\phi}(s,\cdot-\sqrt{2}s))\to 0$ as $s\to+\infty$,
we get by \eqref{4.27} that
$$\log \pp_{\delta_{x}}\left[\exp\{-C(\phi)\partial M_{\infty}\}\right]-\log \pp_{\delta_{x}}\left[\exp\{-C(1-V_{\phi}(s,\cdot-\sqrt{2}s))\partial M_{\infty}\}\right]\to 0\mbox{ locally uniformly in $x\in\R$,}$$
as $s\to +\infty$. Similarly one can prove that
$$\log \pp_{\delta_{x}}\left[\exp\{-C(u^{*}_{\phi}(s,\cdot-\sqrt{2}s))\partial M_{\infty}\}\right]\to 0\mbox{ locally uniformly in $x\in\R$, as }s\to+\infty.$$
Therefore we complete the proof.\qed

\bigskip

\subsection{Proofs of Theorem \ref{them2} and Proposition \ref{prop3}}

\noindent\textbf{Proof of Theorem \ref{them2}:}
Fix an arbitrary $x\in\mathbb{R}$ and functions $f,g$ satisfying our assumptions. We have
\begin{eqnarray}
\pp_{\delta_{x}}\left[\mathrm{e}^{-\langle f,\mathcal{E}_{t}\rangle -\langle g,\mathcal{E}^{Z}_{t}\rangle}\right]
&=&\pp_{\delta_{x}}\left[\mathrm{e}^{\langle \mathcal{T}_{-m(t)}f,X_{t}\rangle}\pp_{\delta_{x}}\left[\mathrm{e}^{-\langle \mathcal{T}_{-m(t)}g,Z_{t}\rangle}|X_{t}\right]\right]\nonumber\\
&=&\pp_{\delta_{x}}\left[\mathrm{e}^{-\langle\mathcal{T}_{-m(t)}f+1-\mathrm{e}^{-\mathcal{T}_{-m(t)}g},X_{t}\rangle}\right]\nonumber\\
&=&\mathrm{e}^{-u_{\mathcal{T}_{-m(t)}(f+1-\mathrm{e}^{-g})}(t,x)}
=
\mathrm{e}^{-u_{f+1-\mathrm{e}^{-g}}(t,x-m(t))}.\nonumber
\end{eqnarray}
The second equality follows from the fact that given $X_{t}$, $Z_{t}$ is a Poisson random measure with intensity $X_{t}(dx)$. Note that $f+1-\mathrm{e}^{-g}\in \mathcal{H}$ by our assumptions.
We get by Theorem \ref{them1} that
$$\lim_{t\to+\infty}\pp_{\delta_{x}}\left[\mathrm{e}^{-\langle f,\mathcal{E}_{t}\rangle -\langle g,\mathcal{E}^{Z}_{t}\rangle}\right]=\pp_{\delta_{x}}\left[\exp\{-C(f+1-\mathrm{e}^{-g})\partial M_{\infty}\}\right].$$
The above identity holds in particular for all $f,g\in C^{+}_{c}(\R)$.
On the other hand, it is easy to verify that the assumptions of Lemma \ref{lem3.6} are satisfied by $\phi_{s}(x):=\phi(x)/s$ for every $\phi\in C^{+}_{c}(\R)$. This implies that $\lim_{\lambda\to 0}C(\lambda\phi)=0$ for all $\phi\in C^{+}_{c}(\R)$. Since $C(\lambda_{1}f+1-\mathrm{e}^{-\lambda_{2}g})\le C(\lambda_{1}f)+C(\lambda_{2}g)$ for $\lambda_{1},\lambda_{2}>0$, we have $\lim_{\lambda_{1},\lambda_{2}\to 0+}C(\lambda_{1}f+1-\mathrm{e}^{-\lambda_{2}g})=0$.
Hence by \cite[Chapter 4]{Kallenberg}
we get the weak convergence of $(\mathcal{E}_{t},\mathcal{E}^{Z}_{t})$,
and consequently \eqref{them2.1} holds for all $f,g\in C^{+}_{c}(\R)$.
Using the monotone convergence theorem and Lemma \ref{lem3.6}, one can show by standard approximation that \eqref{them2.1} holds for all $f\in\mathcal{H}$ and $g\in\mathcal{B}^{+}(\R)$ with $1-\mathrm{e}^{-g}\in\mathcal{H}$.
Also, \eqref{them2.1}  yields that
$$\mathrm{E}\left[\mathrm{e}^{-\langle f,\mathcal{E}_{\infty}\rangle}\mathrm{E}\left[\mathrm{e}^{-\langle g,\mathcal{E}^{Z}_{\infty}\rangle}\,|\,\mathcal{E}_{\infty}\right]\right]=\pp_{\delta_{x}}\left[\exp\{-C\left(f+1-\mathrm{e}^{-g}\right)\partial M_{\infty}\}\right]=\mathrm{E}\left[\mathrm{e}^{-\langle f,\mathcal{E}_{\infty}\rangle}\mathrm{e}^{-\langle 1-\mathrm{e}^{-g},\mathcal{E}_{\infty}\rangle}\right].$$
Thus we get $\mathrm{E}\left[\mathrm{e}^{-\langle g,\mathcal{E}^{Z}_{\infty}\rangle}\,|\,\mathcal{E}_{\infty}\right]=\mathrm{e}^{-\langle 1-\mathrm{e}^{-g},\mathcal{E}_{\infty}\rangle}$ a.s., and the second conclusion follows immediately.\qed

\bigskip

\noindent\textbf{Proof of Proposition \ref{prop3}:}
Suppose $M\in\R$. We have
\begin{eqnarray}
\p\left(\max\mathcal{E}_{\infty}\le M\right)&=&\p\left(\langle1_{(M,+\infty)},\mathcal{E}_{\infty}\rangle=0\right)
=
\lim_{\lambda\to+\infty}\mathrm{E}\left[\mathrm{e}^{-\lambda\langle1_{(M,+\infty)},\mathcal{E}_{\infty}\rangle}\right]\nonumber\\
&=&\lim_{\lambda\to+\infty}\pp_{\delta_{x}}\left[\exp\{-\partial M_{\infty} C\left(\lambda 1_{(M,+\infty)}\right)\}\right]\nonumber\\
&=&\pp_{\delta_{x}}\left[\exp\{-\partial M_{\infty}\lim_{\lambda\to+\infty}C\left(\lambda 1_{(M,+\infty)}\right)\}\right].\nonumber
\end{eqnarray}
The third equality follows from Theorem \ref{them2}. Thus we have
\begin{equation}\label{prop3.1}
\p\left(\max\mathcal{E}_{\infty}<+\infty\right)=\lim_{M\to+\infty}\p\left(\max \mathcal{E}_{\infty}\le M\right)=\pp_{\delta_{x}}\left[\exp\{-\partial M_{\infty} \lim_{M\to+\infty}\lim_{\lambda\to+\infty}C\left(\lambda 1_{(M,+\infty)}\right)\}\right].
\end{equation}
Note that
$$C\left(\lambda 1_{(M,+\infty)}\right)=C\left(\mathcal{T}_{-M}\lambda 1_{(0,+\infty)}\right)=\mathrm{e}^{-\sqrt{2}M}C\left(\lambda 1_{(0,+\infty)}\right).$$
It follows that $\lim_{M\to+\infty}\lim_{\lambda\to+\infty}C\left(\lambda 1_{(M,+\infty)}\right)=\lim_{M\to+\infty}\mathrm{e}^{-\sqrt{2}M}\left(\lim_{\lambda\to+\infty}C\left(\lambda 1_{(0,+\infty)}\right)\right)=0$ or $+\infty$ corresponding to $\sup_{\lambda}C\left(\lambda 1_{(0,+\infty)}\right)<+\infty$ or $+\infty$. Hence we get by \eqref{prop3.1} that $\p\left(\max\mathcal{E}_{\infty}<+\infty\right)=1$ if and only if $\sup_{\lambda}C\left(\lambda 1_{(0,+\infty)}\right)<+\infty$, and prove the result.\qed

\bigskip

\begin{remark}\label{rm4.13}\rm
	
	Recall that $\max X_{t}$ denotes the supremum of the support of $X_{t}$.
	Unlike for the skeleton BBM, Theorem \ref{them1} does not imply the growth order of $\max X_{t}$ is $m(t)$.
We make a short discussion here.

	For $\phi\in\mathcal{B}^{+}_{b}(\R)$, define
	\begin{equation}\label{r0}
		\tilde{u}_{\phi}(t,x):=-\log\pp_{\delta_{x}}\left[\mathrm{e}^{-\langle\phi,X_{t}\rangle};\max X_{t}\le 0\right],\quad\forall t\ge 0,\ x\in\R.
	\end{equation}
	Then one has
	\begin{equation}
		\tilde{u}_{\phi}(t,x)
		=\lim_{\lambda\to +\infty}u_{\phi 1_{(-\infty,0]}+\lambda 1_{(0,+\infty)}}(t,x),\quad\forall t\ge 0,\ x\in\R.\label{r1}
	\end{equation}
	In fact,
	\begin{eqnarray*}
		\mathrm{e}^{-\tilde{u}_{\phi}(t,x)}
		&=&\pp_{\delta_{x}}\left[\mathrm{e}^{-\langle \phi 1_{(-\infty,0]},X_{t}\rangle};\langle 1_{(0,+\infty)},X_{t}\rangle=0\right]\\
		&=&\lim_{\lambda\to+\infty}\pp_{\delta_{x}}\left[\mathrm{e}^{-\langle\phi 1_{(-\infty,0]}+\lambda 1_{(0,+\infty)},X_{t}\rangle}\right]
=
\mathrm{e}^{-\lim_{\lambda\to+\infty}u_{\phi 1_{(-\infty,0]}+\lambda 1_{(0,+\infty)}}(t,x)}.
	\end{eqnarray*}
	By \eqref{r1} and the monotone convergence theorem, one has for $t,r\ge 0$ and $x\in \R$,
	\begin{equation}\label{4.36}
		\tilde{u}_{\phi}(t+r,x)=\lim_{\lambda\to+\infty}u_{\phi 1_{(-\infty,0]}+\lambda 1_{(0,+\infty)}}(t+r,x)=\lim_{\lambda\to+\infty}u_{u_{\phi 1_{(-\infty,0]}+\lambda 1_{(0,+\infty)}}(r,\cdot)}(t,x)=u_{\tilde{u}_{\phi}(r,\cdot)}(t,x).
	\end{equation}
	If we assume in addition that (A3) holds,
	then by \cite[Corollary 3.2]{RSZ} for all $\phi\in\mathcal{H}$ and $r>0$, $\tilde{u}_{\phi}(r,\cdot)\in\mathcal{B}^{+}_{b}(\R)$. So by \eqref{4.36}, for every $\phi\in\mathcal{H}$, $\tilde{u}_{\phi}(t,x)$ can be viewed as a solution to \eqref{eq1} with initial condition $u(0,x)=\phi(x)1_{(-\infty,0]}(x)+\infty 1_{(0,+\infty)}(x)$.
	Moreover, by \cite[Lemma 2.1 and Corollary 3.2]{RSZ},
	for all $\phi\in\mathcal{H}$ and $r>0$, $\tilde{u}_{\phi}(r,\cdot-\sqrt{2}r)\in \mathcal{H}$.
	Hence applying Theorem \ref{them1} to the function $\tilde{u}_{\phi}(t+r,x-\sqrt{2}r)=u_{\tilde{u}_{\phi}(r,\cdot-\sqrt{2}r)}(t,x)$, one gets that for all $\phi\in\mathcal{H}$, the limit
	\begin{equation}\label{defitildec}
		\widetilde{C}(\phi):=\lim_{t\to+\infty}\sqrt{\frac{2}{\pi}}\int_{0}^{+\infty}y\mathrm{e}^{\sqrt{2}y}\tilde{u}_{\phi}(t,-y-\sqrt{2}t)dy
	\end{equation}
	exists and is finite, and for all $x\in \R$,
	\begin{equation}\nonumber
		\pp_{\delta_{0}}\left[\mathrm{e}^{-\langle \phi,X_{t}-m(t)-x\rangle};\max X_{t}-m(t)\le x\right]\to \pp_{\delta_{0}}\left[\mathrm{e}^{-\widetilde{C}(\phi)\partial M_{\infty}\mathrm{e}^{-\sqrt{2}x}}\right]\mbox{ as }t\to +\infty.
	\end{equation}
	Taking $\phi=0$, one gets that
	\begin{equation}\label{r4}
		\pp_{\delta_{0}}\left(\max X_{t}-m(t)\le x\right)\to \pp_{\delta_{0}}\left[\exp\{-\tilde{c}_{0}\partial M_{\infty}\mathrm{e}^{-\sqrt{2}x}\}\right]\mbox{ as }t\to +\infty,
	\end{equation}
	where $\tilde{c}_{0}$ is the constant $\widetilde{C}(\phi)$ with $\phi=0$.

	On the other hand, if (A3) fails, \eqref{r4} may not be true as it is. For instance, we take the triplet $(\alpha,\beta,\pi(dy))$ of $\psi(\lambda)$ to be $(1,0,\delta_{y_{0}}(dy))$ where $y_{0}>1$ satisfies that $\mathrm{e}^{-y_{0}}=2-y_{0}$. In this case
\begin{equation}
	\psi(\lambda)=(y_{0}-1)\lambda-\left(1-\mathrm{e}^{-\lambda y_{0}}\right),\quad \forall \lambda>0.\label{eg}
\end{equation}
	Let $k_{\lambda}(s,y):=-\psi(u_{\lambda 1_{(0,+\infty)}}(s,y))/u_{\lambda 1_{(0,+\infty)}}(s,y)$ for $s\ge 0$ and $y\in \R$. Then by Feynman-Kac formula, one has
\begin{eqnarray}
		u_{\lambda 1_{(0,+\infty)}}(t,x)&=&
\Pi_x\left[\mathrm{e}^{\int_{0}^{t}k_{\lambda}(t-s,B_{s})ds}u_{\lambda 1_{(0,+\infty)}}(0,B_{t})\right]\nonumber\\
		&=&\lambda\mathrm{e}^{-(y_{0}-1)t}
\Pi_x\left[\exp\left\{\int_{0}^{t}\frac{1-\mathrm{e}^{-y_{0}u_{\lambda 1_{(0,+\infty)}}(t-s,B_{s})}}{u_{\lambda 1_{(0,+\infty)}}(t-s,B_{s})}ds\right\};\ B_{t}>0\right]\nonumber\\
&\ge&\lambda \mathrm{e}^{-(y_{0}-1)t}\Pi_x\left(B_{t}>0\right).\nonumber
	\end{eqnarray}
Since
$\Pi_x\left(B_{t}>0\right)>0$
	for all $t>0$ and $x\in \R$, one has $\tilde{u}_{0}(t,x)=\lim_{\lambda\to+\infty}u_{\lambda 1_{(0,+\infty)}}(t,x)=+\infty$. This implies that
	$\pp_{\delta_{0}}\left(\max X_{t}\le x\right)=\pp_{\delta_{-x}}\left(\max X_{t}\le 0\right)=\mathrm{e}^{-\tilde{u}_{0}(t,-x)}=0$
	for all $t>0$ and $x\in\R$, and so by letting $x\to +\infty$, $\pp_{\delta_{0}}\left(\max X_{t}=+\infty\right)=1-\pp_{\delta_{0}}\left(\max X_{t}<+\infty\right)=1$ for all $t>0$.
The branching mechanism given by \eqref{eg} satisfies in particular that
\begin{equation}\label{sheu}
\int_{z}^{+\infty}\frac{1}{\psi(y)}dy=+\infty\quad\forall z>1\mbox{ and }\int_{0}^{1}y\pi(dy)<+\infty.
\end{equation}
In fact, \cite[Theorem 4.4]{Sheu} shows that for a super-Brownian motion with branching mechanism $\psi$ satisfying \eqref{sheu}, it holds that $\pp_{\mu}\left(\mbox{supp}X_{t}=\R\right)=1$ for all $t>0$.
\end{remark}

\section{Probabilistic representation of the limiting process}\label{sec4}

\subsection{Laws of decorations}
For the proofs of Theorem \ref{them3} and Proposition \ref{prop5} we need to show the existence of the limit for $(Z_{t}-\max Z_{t},\sum_{u\in Z_{t}}I^{(u)}_{s}-\sqrt{2}s-\max Z_{t})$ conditioned on $\{\max Z_{t}-\sqrt{2}t>0\}$. This is completed by the following lemmas.
\begin{lemma}\label{lem4.1}
For any $f,g\in\mathcal{B}^{+}_{b}(\mathbb{R})$, $x,z\in\mathbb{R}$ and $t,y\ge 0$, we have
\begin{eqnarray}\label{lem4.1.0}
&&\pp_{\delta_{x}}\left[\exp\{-\langle f,X_{t}-\sqrt{2}t-z\rangle -\langle g,Z_{t}-\sqrt{2}t-z\rangle\}1_{\{\max Z_{t}-\sqrt{2}t-z>y\}}\,|\,\max Z_{t}-\sqrt{2}t-z>0\right]\nonumber\\
&=&\frac{\exp\{-u_{f+\left(1-\mathrm{e}^{-g}\right)}(t,x-\sqrt{2}t-z)\}-\exp\{-u_{f+\left(1-\mathrm{e}^{-g}\right)1_{(-\infty,y]}+1_{(y,+\infty)}}(t,x-\sqrt{2}t-z)\}}{1-\exp\{-u_{1_{(0,+\infty)}}(t,x-\sqrt{2}t-z)\}}.
\end{eqnarray}
\end{lemma}

\proof  We have
\begin{eqnarray}
&&\pp_{\delta_{x}}\left[\exp\{-\langle f,X_{t}-\sqrt{2}t-z\rangle -\langle g,Z_{t}-\sqrt{2}t-z\rangle\}1_{\{\max Z_{t}-\sqrt{2}t-z>y\}}\right]\nonumber\\
&=&\pp_{\delta_{x}}\left[\exp\{-\langle f,X_{t}-\sqrt{2}t-z\rangle\}\pp_{\delta_{x}}\left[\exp\{ -\langle g,Z_{t}-\sqrt{2}t-z\rangle\}1_{\{\max Z_{t}-\sqrt{2}t-z>y\}}|X_{t}\right]\right].\label{lem4.1.1}
\end{eqnarray}
Recall that given $X_{t}$, $Z_{t}$ is a Poisson random measure with intensity $X_{t}(dx)$. Using Poisson computations, we have
\begin{eqnarray*}
&&\pp_{\delta_{x}}\left[\exp\{ -\langle g,Z_{t}-\sqrt{2}t-z\rangle\}1_{\{\max Z_{t}-\sqrt{2}t-z>y\}}|X_{t}\right]\nonumber\\
&=&\pp_{\delta_{x}}\left[\exp\{ -\langle g,Z_{t}-\sqrt{2}t-z\rangle\}\left(1-1_{\{\max Z_{t}-\sqrt{2}t-z\le y\}}\right)|X_{t}\right]\\
&=&\pp_{\delta_{x}}\left[\exp\{ -\langle g,Z_{t}-\sqrt{2}t-z\rangle\}|X_{t}\right]-\pp_{\delta_{x}}\left[\exp\{ -\langle g,Z_{t}-\sqrt{2}t-z\rangle\}1_{\{\langle 1_{(y,+\infty)}, Z_{t}-\sqrt{2}t-z\rangle =0\}}|X_{t}\right]\\
&=&\exp\{-\langle 1-\mathrm{e}^{-g},X_{t}-\sqrt{2}t-z\rangle\}-\exp\{-\langle(1-\mathrm{e}^{-g})1_{(-\infty,y]}+1_{(y,+\infty)},X_{t}-\sqrt{2}t-z\rangle\}.
\end{eqnarray*}
Putting this back in \eqref{lem4.1.1} we get that
\begin{eqnarray}
&&\pp_{\delta_{x}}\left[\exp\{-\langle f,X_{t}-\sqrt{2}t-z\rangle -\langle g,Z_{t}-\sqrt{2}t-z\rangle\}1_{\{\max Z_{t}-\sqrt{2}t-z>y\}}\right]\nonumber\\
&=&
\pp_{\delta_{x}}\left[\exp\{-\langle f+(1-\mathrm{e}^{-g}),X_{t}-\sqrt{2}t-z\rangle\}\right]\nonumber\\
&&\quad-\pp_{\delta_{x}}\left[\exp\{-\langle f+(1-\mathrm{e}^{-g})1_{(-\infty,y]}+1_{(y,+\infty)},X_{t}-\sqrt{2}t-z\rangle\}\right]\nonumber\\
&=&
\exp\{-u_{f+\left(1-\mathrm{e}^{-g}\right)}(t,x-\sqrt{2}t-z)\}-\exp\{-u_{f+\left(1-\mathrm{e}^{-g}\right)1_{(-\infty,y]}+1_{(y,+\infty)}}(t,x-\sqrt{2}t-z)\}.\label{lem4.1.2}
\end{eqnarray}
In particular by setting $f=g=0$ and $y=0$ in the above formula, we get that
\begin{equation}
\pp_{\delta_{x}}\left(\max Z_{t}-\sqrt{2}t-z>0\right)=1-\exp\{-u_{1_{(0,+\infty)}}(t,x-\sqrt{2}t-z)\}.\label{lem4.1.3}
\end{equation}
Hence \eqref{lem4.1.0} follows by $\eqref{lem4.1.2}/\eqref{lem4.1.3}$.

\qed

\begin{lemma}\label{lem4.2}
For any $f,g\in\mathcal{B}^{+}_{b}(\mathbb{R})$, $x,z\in\mathbb{R}$, $t\ge 0$ and $\lambda>0$,
\begin{eqnarray}\label{lem4.2.0}
&&\pp_{\delta_{x}}\left[\exp\left\{-\langle f,X_{t}-\sqrt{2}t-z\rangle -\langle g,Z_{t}-\sqrt{2}t-z\rangle-\lambda\left(\max Z_{t}-\sqrt{2}t-z\right)\right\}\,\left|\,\max Z_{t}-\sqrt{2}t-z>0\right.\right]\nonumber\\
&=&
\Big[\int_{0}^{+\infty}\mathrm{e}^{-y}\exp\left\{-u_{f+\left(1-\mathrm{e}^{-g}\right)1_{(-\infty,\frac{y}{\lambda}]}
+1_{(\frac{y}{\lambda},+\infty)}}(t,x-\sqrt{2}t-z)\right\}dy\nonumber\\
&& -\exp\left\{-u_{f+\left(1-\mathrm{e}^{-g}\right)
1_{(-\infty,0]}+1_{(0,+\infty)}}(t,x-\sqrt{2}t-z)\right\}\Big]
\left(1-\exp\left\{-u_{1_{(0,+\infty)}}(t,x-\sqrt{2}t-z)\right\}\right)^{-1}.
\end{eqnarray}
\end{lemma}

\proof We rewrite the left hand side of \eqref{lem4.2.0} as $I-II$ where
$$I:=\pp_{\delta_{x}}\left[\exp\{-\langle f,X_{t}-\sqrt{2}t-z\rangle -\langle g,Z_{t}-\sqrt{2}t-z\rangle\,|\,\max Z_{t}-\sqrt{2}t-z>0\right]$$
and
\begin{eqnarray*}
II&:=&\pp_{\delta_{x}}\Big[\exp\left\{-\langle f,X_{t}-\sqrt{2}t-z\rangle -\langle g,Z_{t}-\sqrt{2}t-z\rangle\right\}\\
&&\qquad\cdot\left(1-\exp\left\{-\lambda\left(\max Z_{t}-\sqrt{2}t-z\right)\right\}\right)\,\big|\,\max Z_{t}-\sqrt{2}t-z>0\Big].
\end{eqnarray*}
By Lemma \ref{lem4.1},
\begin{equation}\label{lem4.2.1}
I=\frac{\exp\{-u_{f+\left(1-\mathrm{e}^{-g}\right)}(t,x-\sqrt{2}t-z)\}-\exp\{-u_{f+\left(1-\mathrm{e}^{-g}\right)1_{(-\infty,0]}+1_{(0,+\infty)}}(t,x-\sqrt{2}t-z)\}}{1-\exp\{-u_{1_{(0,+\infty)}}(t,x-\sqrt{2}t-z)\}}.
\end{equation}
On the other hand, by Fubini's theorem and Lemma \ref{lem4.1} we have
\begin{eqnarray}
&&II\nonumber\\
&=&\pp_{\delta_{x}}\left[\exp\{-\langle f,X_{t}-\sqrt{2}t-z\rangle -\langle g,Z_{t}-\sqrt{2}t-z\rangle\}\int_{0}^{\lambda(\max Z_{t}-\sqrt{2}t-z)}\mathrm{e}^{-y}dy\,\big|\,\max Z_{t}-\sqrt{2}t-z>0\right]\nonumber\\
&=&\int_{0}^{+\infty}\mathrm{e}^{-y}\pp_{\delta_{x}}\left[
\mathrm{e}^{-\langle f,X_{t}-\sqrt{2}t-z\rangle -\langle g,Z_{t}-\sqrt{2}t-z\rangle}
1_{\{\max Z_{t}-\sqrt{2}t-z>\frac{y}{\lambda}\}}\,\big|\,\max Z_{t}-\sqrt{2}t-z>0\right]dy\nonumber\\
&=&\left(1-\exp\left\{-u_{1_{(0,+\infty)}}(t,x-\sqrt{2}t-z)\right\}\right)^{-1}
\Big[\exp\left\{-u_{f+\left(1-\mathrm{e}^{-g}\right)}(t,x-\sqrt{2}t-z)\right\}\nonumber\\
&&-\int_{0}^{+\infty}\mathrm{e}^{-y}\exp\left\{-u_{f+\left(1-\mathrm{e}^{-g}\right)1_{(-\infty,\frac{y}{\lambda}]}
+1_{(\frac{y}{\lambda},+\infty)}}(t,x-\sqrt{2}t-z)\right\}dy\Big].\label{lem4.2.2}
\end{eqnarray}
Hence we get \eqref{lem4.2.0} by letting
$\eqref{lem4.2.1}-\eqref{lem4.2.2}$.
\qed

\bigskip

\begin{lemma}\label{lem4.3}
Suppose $x,z\in\mathbb{R}$ and $s>0$. Under $\pp_{\delta_{x}}$, conditioned on $\{\max Z_{t}-\sqrt{2}t-z>0\}$, the random elements $\left(X_{t}-\sqrt{2}t-z,\sum_{u\in Z_{t}}I^{(u)}_{s}-\sqrt{2}(s+t)-z,
Z_{t}-\sqrt{2}t-z,\max Z_{t}-\sqrt{2}t-z\right)$ converges, as $t\to +\infty$, in distribution to a limit
$\left(\bar{\mathcal{E}}^{X}_{\infty},\bar{\mathcal{E}}^{I,s}_{\infty},\bar{\mathcal{E}}^{Z}_{\infty},Y\right)$,
where the limit is independent of $x$ and $z$, and $Y$ is an exponential random variable with mean $1/\sqrt{2}$. Moreover, we have for any $f,h\in\mathcal{H}$, $g\in\mathcal{B}^{+}(\mathbb{R})$ with $1-\mathrm{e}^{-g}\in\mathcal{H}_{1}$ and $\lambda\ge 0$,
\begin{eqnarray}
&&\mathrm{E}\left[\mathrm{e}^{-\langle f,\bar{\mathcal{E}}^{X}_{\infty}\rangle-\langle h,\bar{\mathcal{E}}^{I,s}_{\infty}\rangle-\langle g,\bar{\mathcal{E}}^{Z}_{\infty}\rangle-\lambda Y}\right]\nonumber\\
&=&\lim_{t\to+\infty}\pp_{\delta_{x}}\big[\exp\{-\langle f,X_{t}-\sqrt{2}t-z\rangle-\langle h,\sum_{u\in Z_{t}}I^{(u)}_{s}-\sqrt{2}(s+t)-z\rangle-\langle g,Z_{t}-\sqrt{2}t-z\rangle\nonumber\\
&&\quad\quad\quad\quad-\lambda(\max Z_{t}-\sqrt{2}t-z)\}\,|\,\max Z_{t}-\sqrt{2}t-z>0\big]\nonumber\\
&=&
\frac{1}{c_{*}}
\Big[ C\left(f+\left(1-\mathrm{e}^{-g}V_{h}(s,\cdot-\sqrt{2}s)\right)1_{(-\infty,0]}+1_{(0,+\infty)}\right)\nonumber\\
&&\quad -\int_{0}^{+\infty}\mathrm{e}^{-y}C\left(f+\left(1-\mathrm{e}^{-g}V_{h}(s,\cdot-\sqrt{2}s)\right)
1_{(-\infty,\frac{y}{\lambda}]}+1_{(\frac{y}{\lambda},+\infty)}\right)dy\Big],
\end{eqnarray}
where $c_{*}=C(1_{(0,+\infty)})$.
\end{lemma}

\proof In view of \eqref{lem4.1.3} we have for any $x,z\in\mathbb{R}$ and $y\ge 0$,
\begin{eqnarray}
&&
\lim_{t\to+\infty}\pp_{\delta_{x}}\left(\max Z_{t}-\sqrt{2}t-z>y\,|\,\max Z_{t}-\sqrt{2}t-z>0\right)
\nonumber\\
&=&\lim_{t\to +\infty}\frac{\pp_{\delta_{x}}\left(\max Z_{t}-\sqrt{2}t-z>y\right)}{\pp_{\delta_{x}}\left(\max Z_{t}-\sqrt{2}t-z>0\right)}\nonumber\\
&=&\lim_{t\to+\infty}\frac{1-\exp\{-u_{1_{(0,+\infty)}}(t,x-\sqrt{2}t-z-y)\}}{1-\exp\{-u_{1_{(0,+\infty)}}(t,x-\sqrt{2}t-z)\}}\nonumber\\
&=&\lim_{t\to +\infty}\frac{u_{1_{(0,+\infty)}}(t,x-\sqrt{2}t-z-y)}{u_{1_{(0,+\infty)}}(t,x-\sqrt{2}t-z)}
=
\mathrm{e}^{-\sqrt{2}y}.\nonumber
\end{eqnarray}
The final equality follows from Lemma \ref{lem3.3}. This implies that,
conditioned on
$\{\max Z_{t}-\sqrt{2}t-z>0\}$, $\max Z_{t}-\sqrt{2}t-z$
converges in distribution to an exponentially distributed random variable with mean $1/\sqrt{2}$.

Suppose $f,h,g$ are functions satisfying our assumptions.
Recall that
$\mathcal{F}_{t}$ is the $\sigma$-filed generated by $Z$, $X^{*}$ and $I$ up to time $t$, and
given $\mathcal{F}_{t}$, $I^{(u)}_{s}\stackrel{\mbox{d}}{=}(I_{s},\pp_{\cdot,\delta_{z_{u}(t)}})$ for $u\in Z_{t}$. We have
\begin{eqnarray}
\pp_{\delta_{x}}\left[\exp\left\{-\langle h,\sum_{u\in Z_{t}}I^{(u)}_{s}-\sqrt{2}(t+s)-z\rangle \right\}\,\Big|\,\mathcal{F}_{t}\right]
&=&\prod_{u\in Z_{t}}\pp_{\cdot,\delta_{z_{u}(t)}}\left[\exp\{-\langle h,I_{s}-\sqrt{2}(s+t)-z\rangle\}\right]\nonumber\\
&=&\prod_{u\in Z_{t}}V_{h}(s,z_{u}(t)-\sqrt{2}(s+t)-z)\nonumber\\
&=&\exp\{\langle \ln V_{h}(s,\cdot-\sqrt{2}s),Z_{t}-\sqrt{2}t-z\rangle\},
\end{eqnarray}
where $V_h$ is defined by \eqref{def-V} with $f$ replaced by $h$.
Thus we have
\begin{eqnarray}
&&\pp_{\delta_{x}}\left[\mathrm{e}^{-\langle f,X_{t}-\sqrt{2}t-z\rangle -\langle h,\sum_{u\in Z_{t}}I^{(u)}_{s}-\sqrt{2}(s+t)-z\rangle-
\langle g,Z_{t}-\sqrt{2}t-z\rangle-\lambda\left(\max Z_{t}-\sqrt{2}t-z\right)}\,|\,\max Z_{t}-\sqrt{2}t-z>0\right]\nonumber\\
 &=&\pp_{\delta_{x}}\Big[\mathrm{e}^{-\langle f,X_{t}-\sqrt{2}t-z\rangle -
\langle g,Z_{t}-\sqrt{2}t-z\rangle-\lambda\left(\max Z_{t}-\sqrt{2}t-z\right)}\nonumber\\
 &&\qquad\cdot\pp_{\delta_{x}}\left[\mathrm{e}^{-\langle h,\sum_{u\in Z_{t}}I^{(u)}_{s}-\sqrt{2}(t+s)-z\rangle }\,|\,\mathcal{F}_{t}\Big]
\,\big|\,\max Z_{t}-\sqrt{2}t-z>0\right]\nonumber\\
&=&\pp_{\delta_{x}}\left[\mathrm{e}^{-\langle f,X_{t}-\sqrt{2}t-z\rangle -
\langle g-\ln V_{h}(s,\cdot-\sqrt{2}s),Z_{t}-\sqrt{2}t-z\rangle-\lambda\left(\max Z_{t}-\sqrt{2}t-z\right)}\,\big|\,\max Z_{t}-\sqrt{2}t-z>0\right].\label{5.10}
\end{eqnarray}

By Lemma \ref{lem4.2} the right hand side of \eqref{5.10} equals
\begin{eqnarray}
&&\Big[\int_{0}^{+\infty}\mathrm{e}^{-y}\exp\left\{-u_{f+\left(1-\mathrm{e}^{-g}V_{h}(s,\cdot-\sqrt{2}s)\right)1_{(-\infty,\frac{y}{\lambda}]}
+1_{(\frac{y}{\lambda},+\infty)}}(t,x-\sqrt{2}t-z)\right\}dy\nonumber\\
&&-\exp\left\{-u_{f+\left(1-\mathrm{e}^{-g}V_{h}(s,\cdot-\sqrt{2}s)\right)
1_{(-\infty,0]}+1_{(0,+\infty)}}(t,x-\sqrt{2}t-z)\right\}\Big]\left(1-\mathrm{e}^{-u_{1_{(0,+\infty)}}(t,x-\sqrt{2}t-z)}\right)^{-1}\nonumber\\
&=&\Big[-\int_{0}^{+\infty}\mathrm{e}^{-y}\left(1-\mathrm{e}^{-u_{f+\left(1-\mathrm{e}^{-g}V_{h}(s,\cdot-\sqrt{2}s)\right)1_{(-\infty,\frac{y}{\lambda}]}
+1_{(\frac{y}{\lambda},+\infty)}}(t,x-\sqrt{2}t-z)}\right)dy\nonumber\\
&&\quad+\left(1 -\mathrm{e}^{-u_{f+\left(1-\mathrm{e}^{-g}V_{h}(s,\cdot-\sqrt{2}s)\right)
1_{(-\infty,0]}+1_{(0,+\infty)}}(t,x-\sqrt{2}t-z)}\right)\Big]\left(1-\mathrm{e}^{-u_{1_{(0,+\infty)}}(t,x-\sqrt{2}t-z)}\right)^{-1}.\label{lem4.3.1}
\end{eqnarray}
We observe that $1-\mathrm{e}^{-g}V_{h}(s,\cdot-\sqrt{2}s)=\mathrm{e}^{-g}\left(1-V_{h}(s,\cdot-\sqrt{2}s)\right)+(1-\mathrm{e}^{-g})\in\mathcal{H}_{1}$ since $1-V_{h}(s,\cdot-\sqrt{2}s),1-\mathrm{e}^{-g}\in \mathcal{H}_{1}$ by the assumptions.
Using the facts that $1-\mathrm{e}^{-x}\sim x$ as $x\to 0$ and that
$$\lim_{t\to+\infty}\frac{t^{3/2}}{\frac{3}{2\sqrt{2}}\log t}\mathrm{e}^{-\sqrt{2}x}u_{\phi}(t,x-\sqrt{2}t)=C(\phi)\quad\forall \phi\in\mathcal{H},$$
one can show by the bounded convergence theorem that the right hand side of \eqref{lem4.3.1} converges to
\begin{eqnarray*}
&&
\frac{1}{c_{*}}
\Big[ C\left(f+\left(1-\mathrm{e}^{-g}V_{h}(s,\cdot-\sqrt{2}s)\right)1_{(-\infty,0]}+1_{(0,+\infty)}\right)\nonumber\\
&&\quad -\int_{0}^{+\infty}\mathrm{e}^{-y}C\left(f+\left(1-\mathrm{e}^{-g}V_{h}(s,\cdot-\sqrt{2}s)\right)
1_{(-\infty,\frac{y}{\lambda}]}+1_{(\frac{y}{\lambda},+\infty)}\right)dy\Big]
\end{eqnarray*}
as $t\to+\infty$.
In particular, for any $\lambda_{i}>0$, $i=1,2,3,4$ and $f,h,g\in C^{+}_{c}(\mathbb{R})$, one has
\begin{eqnarray}\label{lem4.3.6}
&&\lim_{t\to +\infty}\pp_{\delta_{x}}\Big[\mathrm{e}^{-\lambda_{1}\langle f,X_{t}-\sqrt{2}t-z\rangle -\lambda_{2}\langle h,\sum_{u\in Z_{t}}I^{(u)}_{s}-\sqrt{2}(s+t)-z\rangle-
\lambda_{3}\langle g,Z_{t}-\sqrt{2}t-z\rangle-\lambda_{4}\left(\max Z_{t}-\sqrt{2}t-z\right)}\nonumber\\
&&\hspace{1.8cm}\big|\,\max Z_{t}-\sqrt{2}t-z>0\Big]\nonumber\\
&=&
\frac{1}{c_{*}}
\Big[ C\left(\lambda_{1}f+\left(1-\mathrm{e}^{-\lambda_{3}g}V_{\lambda_{2}h}(s,\cdot-\sqrt{2}s)\right)1_{(-\infty,0]}+1_{(0,+\infty)}\right)\nonumber\\
&&\quad -\int_{0}^{+\infty}\mathrm{e}^{-y}C\left(\lambda_{1}f+\left(1-\mathrm{e}^{-\lambda_{3}g}V_{\lambda_{2}h}(s,\cdot-\sqrt{2}s)\right)
1_{(-\infty,\frac{y}{\lambda_{4}}]}+1_{(\frac{y}{\lambda_{4}},+\infty)}\right)dy\Big]\label{lem4.3.2}
\end{eqnarray}
To show the convergence in distribution, it suffices to show the right hand side of \eqref{lem4.3.2} converges to $1$ as $\lambda_{i}\to 0$, $i=1,2,3,4$.
By Corollary \ref{cor4} and Lemma \ref{lem3.6}
we have
\begin{equation*}
C(1-V_{\lambda_{2}h}(s,\cdot-\sqrt{2}s))
\le
C(u_{\lambda_{2}h}(s,\cdot-\sqrt{2}s))
=C(\lambda_{2}h)
\to 0\mbox{ as }\lambda_{2}\to 0,
\end{equation*}
and
$$C(1-\mathrm{e}^{-\lambda_{3}g})\le C(\lambda_{3}g)
\to 0,\quad \mbox{ as }\lambda_{3}\to 0.
$$
Thus one has
\begin{eqnarray}\label{lem4.3.3}
C(1-\mathrm{e}^{-\lambda_{3}g}V_{\lambda_{2}h}(s,\cdot-\sqrt{2}s))&\le& C\left(\mathrm{e}^{-\lambda_{3}g}\left(1-V_{\lambda_{2}h}(s,\cdot-\sqrt{2}s)\right)\right)+C(1-\mathrm{e}^{-\lambda_{3}g})\nonumber\\
&\le &C\left(1-V_{\lambda_{2}h}(s,\cdot-\sqrt{2}s)\right)+C(1-\mathrm{e}^{-\lambda_{3}g})\nonumber\\
&\to &0,\quad\mbox{ as }\lambda_{2},\lambda_{3}\to 0.
\end{eqnarray}
On the other hand by Corollary \ref{cor4} we have for any $\delta\in\mathbb{R}$
\begin{eqnarray}
C(1_{(\delta,+\infty)})&\le& C(\lambda_{1}f+\left(1-\mathrm{e}^{-\lambda_{3}g}V_{\lambda_{2}h}(s,\cdot-\sqrt{2}s)\right)1_{(-\infty,\delta]}+1_{(\delta,+\infty)})\nonumber\\
&\le&C(\lambda_{1}f)+C\left(1-\mathrm{e}^{-\lambda_{3}g}V_{\lambda_{2}h}(s,\cdot-\sqrt{2}s)\right)+C(1_{(\delta,+\infty)}),\label{5.15}
\end{eqnarray}
Using \eqref{lem4.3.3} and the fact that $\lim_{\lambda_{1}\to 0}C(\lambda_{1}f)=0$ we get by \eqref{5.15} that
\begin{equation}\nonumber
\lim_{\lambda_{1},\lambda_{2},\lambda_{3}\to 0}C(\lambda_{1}f+\left(1-\mathrm{e}^{-\lambda_{3}g}V_{\lambda_{2}h}(s,\cdot-\sqrt{2}s)\right)1_{(-\infty,\delta]}+1_{(\delta,+\infty)})=C(1_{(\delta,+\infty)}).
\end{equation}
This implies that
\begin{equation}\label{lem4.3.4}
\lim_{\lambda_{1},\lambda_{2},\lambda_{3}\to 0}C(\lambda_{1}f+\left(1-\mathrm{e}^{-\lambda_{3}g}V_{\lambda_{2}h}(s,\cdot-\sqrt{2}s)\right)1_{(-\infty,0]}+1_{(0,+\infty)})=C(1_{(0,+\infty)}),
\end{equation}
and that for every $y\in\mathbb{R}$,
\begin{equation}\label{lem4.3.5}
\lim_{\lambda_{1},\lambda_{2},\lambda_{3},\lambda_{4}\to 0}C(\lambda_{1}f+\left(1-\mathrm{e}^{-\lambda_{3}g}V_{\lambda_{2}h}(s,\cdot-\sqrt{2}s)\right)1_{(-\infty,\frac{y}{\lambda_{4}}]}+1_{(\frac{y}{\lambda_{4}},+\infty)})=\lim_{\lambda_{4}\to 0}C(1_{(\frac{y}{\lambda_{4}},+\infty)})=0.
\end{equation}
In view of \eqref{lem4.3.4} and \eqref{lem4.3.5}, one can use the bounded convergence theorem to show that the right hand side of \eqref{lem4.3.6} converges to $1$ as $\lambda_{i}\to 0$, $i=1,2,3,4$. Hence we complete the proof.\qed

\bigskip

\begin{lemma}\label{lem4.4}
Suppose $x,z\in\mathbb{R}$ and $s>0$. Under $\pp_{\delta_{x}}$, conditioned on $\{\max Z_{t}-\sqrt{2}t-z>0\}$, the random elements $\left(X_{t}-\max Z_{t},\sum_{u\in Z_{t}}I^{(u)}_{s}-\sqrt{2}s-\max Z_{t},
Z_{t}-\max Z_{t},\max Z_{t}-\sqrt{2}t-z\right)$ converges, as $t\to +\infty$, in distribution to a limit
$\left(\triangle^{X},\triangle^{I,s},\triangle^{Z},Y\right):=
\left(\bar{\mathcal{E}}^{X}_{\infty}-Y,\bar{\mathcal{E}}^{I,s}_{\infty}-Y,\bar{\mathcal{E}}^{Z}_{\infty}-Y, Y\right)$,
which is independent of $x$ and $z$.
Moreover $\left(\triangle^{X},\triangle^{I,s},\triangle^{Z}\right)$ is independent of $Y$.
\end{lemma}

\proof The first conclusion is a direct result of Lemma \ref{lem4.3} and \cite[Lemma 4.13]{ABK}. We only need to show the independence. Suppose $f,g,h\in C^{+}_{c}(\mathbb{R})$ and $y\ge 0$. We have
\begin{eqnarray}
&&\mathrm{E}\left[\mathrm{e}^{-\langle f,\triangle^{X}\rangle-\langle g,\triangle^{I,s}\rangle-\langle h,\triangle^{Z}\rangle}1_{\{Y>y\}}\right]\nonumber\\
&=&\lim_{t\to +\infty}\pp_{\delta_{0}}\big[\mathrm{e}^{-\langle f,X_{t}-\max Z_{t}\rangle-\langle g,\sum_{u\in Z_{t}}I^{(u)}_{s}-\sqrt{2}s-\max Z_{t}\rangle
-\langle h, Z_{t}-\max Z_{t}\rangle}\nonumber\\
&&\quad\times 1_{\{\max Z_{t}-\sqrt{2}t>y\}}|\max Z_{t}-\sqrt{2}t>0\big]\nonumber\\
&=&\lim_{t\to +\infty}\pp_{\delta_{0}}\big[\mathrm{e}^{-\langle f,X_{t}-\max Z_{t}\rangle-\langle g,\sum_{u\in Z_{t}}I^{(u)}_{s}-\sqrt{2}s-\max Z_{t}\rangle
-\langle h, Z_{t}-\max Z_{t}\rangle}\big|\max Z_{t}-\sqrt{2}t>y\big]\nonumber\\
&&\quad\times \pp_{\delta_{0}}\left(\max Z_{t}-\sqrt{2}t>y|\max Z_{t}-\sqrt{2}t>0\right)\nonumber\\
&=&\mathrm{E}\left[\mathrm{e}^{-\langle f,\triangle^{X}\rangle-\langle g,\triangle^{I,s}\rangle-\langle h,\triangle^{Z}\rangle}\right]\times \mathrm{P}(Y>y).\nonumber
\end{eqnarray}
This yields the independence.\qed

\bigskip

\begin{remark}\rm
We have by \eqref{uf}  that for all $g\in C^{+}_{c}(\R)$ and $y\ge 0$,
\begin{eqnarray*}
\pp_{\cdot,\delta_{0}}\left[\mathrm{e}^{-\langle g,Z_{t}-\sqrt{2}t\rangle};\max Z_{t}-\sqrt{2}t>y\right]&=&u_{(1-\mathrm{e}^{-g})1_{(-\infty,y]}+1_{(y,+\infty)}}(t,-\sqrt{2t})-u_{1-\mathrm{e}^{-g}}(t,-\sqrt{2}t).
\end{eqnarray*}
Thus by Lemma \ref{lem4.1} with $f=0$ we have
\begin{eqnarray*}
&&\pp_{\delta_{0}}\left[\mathrm{e}^{-\langle g,Z_{t}-\sqrt{2}t\rangle};\max Z_{t}-\sqrt{2}t>y\,\big|\,\max Z_{t}-\sqrt{2}t>0\right]\\
&=&\frac{\mathrm{e}^{-u_{1-\mathrm{e}^{-g}}(t,-\sqrt{2}t)}-\mathrm{e}^{-u_{(1-\mathrm{e}^{-g})1_{(-\infty,y]}+1_{(y,+\infty)}}(t,-\sqrt{2t})}}{1-\mathrm{e}^{-u_{1_{(0,+\infty)}}(t,-\sqrt{2}t)}}\\
&\sim&\frac{u_{(1-\mathrm{e}^{-g})1_{(-\infty,y]}+1_{(y,+\infty)}}(t,-\sqrt{2t})-u_{1-\mathrm{e}^{-g}}(t,-\sqrt{2}t)}{u_{1_{(0,+\infty)}}(t,-\sqrt{2}t)}\\
&=&\pp_{\cdot,\delta_{0}}\left[\mathrm{e}^{-\langle g,Z_{t}-\sqrt{2}t\rangle};\max Z_{t}-\sqrt{2}t>y\,\big|\,\max Z_{t}-\sqrt{2}t>0\right],
\quad \mbox{ as }t\to+\infty.
\end{eqnarray*}
This implies that the limit of $(Z_{t}-\sqrt{2}t,\max Z_{t}-\sqrt{2}t)$ under $\pp_{\delta_{0}}\left(\cdot|\max Z_{t}-\sqrt{2}t>0\right)$ and that under $\pp_{\cdot,\delta_{0}}(\cdot|\max Z_{t}-\sqrt{2}t>0)$ are equal in distribution. Therefore the definition of $\triangle^{Z}$ given in Lemma \ref{lem4.4}
coincides with
that of \eqref{def-triangleZ}.
\end{remark}

\subsection{Proof of Theorem \ref{them3}}
In the following lemma, we establish an integral representation of $C(\phi)$ for $\phi\in\mathcal{H}_{1}$. It characterises the limiting extremal process of the skeleton BBM as a decorated Poisson point process.
\begin{lemma}\label{lem4.8}
For $\phi\in\mathcal{H}_{1}$,
\begin{equation}\label{lem4.8.0}
C(\phi)=c_{*}\int_{-\infty}^{+\infty}\sqrt{2}\mathrm{e}^{-\sqrt{2}y}\mathrm{E}\left[1-\mathrm{e}^{\langle \ln (1-\phi),\triangle^{Z}+y\rangle}\right]dy.
\end{equation}
Here $\triangle^{Z}$ is defined in Lemma \ref{lem4.4}, and $c_{*}=C(1_{(0,+\infty)})$.
\end{lemma}

\proof
First we suppose $\phi\in\mathcal{H}_{1}$ is a compactly supported continuous function with $\|\phi\|_{\infty}<1$. Let $g(x):=-\ln (1-\phi(x))$ for $x\in\R$.
The argument below Proposition \ref{lem3.5} implies that
\begin{equation}\label{lem4.8.1}
\mathrm{E}\left[\mathrm{e}^{-\langle g,\mathcal{E}^{Z}_{\infty}\rangle}\right]=\lim_{t\to+\infty}\pp_{\cdot,\delta_{0}}\left[\mathrm{e}^{-\langle g,\mathcal{E}^{Z}_{t}\rangle}\right]=\pp_{\cdot,\delta_{0}}\left[\mathrm{e}^{-C(\phi)\partial M_{\infty}}\right],
\end{equation}
where $\mathcal{E}^{Z}_{\infty}$ is the limit of $((\mathcal{E}^{Z}_{t})_{t\ge 0},\pp_{\cdot,\delta_{0}})$ in distribution. It is known that
$\mathcal{E}^{Z}_{\infty}$ is a DPPP($c_{*}\partial M_{\infty}\sqrt{2}\mathrm{e}^{\sqrt{2}y}dy,$
$\triangle^{Z})$ where $\triangle^{Z}$
is the distributional limit of $Z_{t}-\max Z_{t}$ under $\pp_{\cdot,\delta_{0}}\left(\cdot|\max Z_{t}-\sqrt{2}t>0\right)$ (and hence equal in law to $\triangle^{Z}$ defined in Lemma \ref{lem4.4}).  Using Poisson computations one has
$$\mathrm{E}\left[\mathrm{e}^{-\langle g,\mathcal{E}^{Z}_{\infty}\rangle}\right]=\pp_{\cdot,\delta_{0}}\left[\exp\{-c_{*}\partial M_{\infty}\int_{-\infty}^{+\infty}\sqrt{2}\mathrm{e}^{-\sqrt{2}y}\mathrm{E}\left[1-\mathrm{e}^{-\langle g,\triangle^{Z}+y\rangle}\right]dy\}\right].$$
We note that under our assumptions, $\pp_{\cdot,\delta_{0}}\left[\partial M_{\infty}>0\right]>0$. Hence \eqref{lem4.8.0} follows, otherwise there would be a contradiction.

For a general $\phi\in\mathcal{H}_{1}$, one can find a nondecreasing sequence of functions $\{\phi_{n}:n\ge 1\}\subset \mathcal{H}_{1}\cap C^{+}_{c}(\mathbb{R})$,
such that $\|\phi_{n}\|_{\infty}<1$ and $\phi_{n}(x)\uparrow \phi(x)$ for all $x\in \R$ as $n\to +\infty$.
Then by Corollary \ref{cor4},
$$C(\phi_{n})\le C(\phi)\le C(\phi_{n})+C(\phi-\phi_{n}).$$
Note that by the dominated convergence theorem $\int_{-\infty}^{+\infty}|x|\mathrm{e}^{\sqrt{2}x}\left(\phi(-x)-\phi_{n}(-x)\right)dx\to 0$ as $n\to +\infty$.
Thus by Lemma \ref{lem3.6} $C(\phi-\phi_{n})\to 0$ as $n\to +\infty$. So we get that $C(\phi_{n})\uparrow C(\phi)$ as $n\to+\infty$, and \eqref{lem4.8.0} follows immediately by the monotone convergence theorem.\qed

\bigskip

\noindent\textbf{Proof of Theorem \ref{them3}:}
It follows from Theorem \ref{them2} and Lemma \ref{lem4.8} that for all $g\in\mathcal{B}^{+}(\R)$ with $1-\mathrm{e}^{-g}\in\mathcal{H}$,
\begin{eqnarray}
\mathrm{E}\left[\mathrm{e}^{-\langle g,\mathcal{E}^{Z}_{\infty}\rangle}\right]&=&\pp_{\delta_{x}}\left[\mathrm{e}^{-C(1-\mathrm{e}^{-g})\partial M_{\infty}}\right]\label{them3.1}\\
&=&\pp_{\delta_{x}}\left[\exp\{-c_{*}\partial M_{\infty}\int_{-\infty}^{+\infty}\sqrt{2}\mathrm{e}^{-\sqrt{2}y}\mathrm{E}\left[1-\mathrm{e}^{-\langle g,\triangle^{Z}+y\rangle}\right]dy\}\right].\nonumber
\end{eqnarray}
 The above equations hold in particular for all $g\in C^{+}_{c}(\R)$. This implies that $\mathcal{E}^{Z}_{\infty}$ is a decorated Poisson point process with intensity $c_{*}\partial M_{\infty}\sqrt{2}\mathrm{e}^{-\sqrt{2}y}dy$ and decoration law $\triangle^{Z}$.

We have for all $\phi\in C^{+}_{c}(\R)$,
\begin{eqnarray}
\mathrm{E}\left[\mathrm{e}^{-\langle \phi,\sum_{i\ge 1}\mathcal{T}_{d_{i}}\triangle^{s}_{i}\rangle}\right]
&=&\mathrm{E}\left[\prod_{i\ge 1}\pp_{\cdot,\delta_{0}}\left[\mathrm{e}^{-\langle \mathcal{T}_{d_{i}}\phi,I_{s}-\sqrt{2}s\rangle}\right]\right]
=
\mathrm{E}\left[\prod_{i\ge 1}V_{\mathcal{T}_{d_{i}}\phi}(s,-\sqrt{2}s)\right]\nonumber\\
&=&\mathrm{E}\left[\mathrm{e}^{\langle \ln V_{\phi}(s,\cdot-\sqrt{2}s),\mathcal{E}^{Z}_{\infty}\rangle}\right]
=
\pp_{\delta_{x}}\left[\mathrm{e}^{-C(1-V_{\phi}(s,\cdot-\sqrt{2}s))\partial M_{\infty}}\right].\nonumber
\end{eqnarray}
The final equality follows from \eqref{them3.1}. Since $\lim_{s\to+\infty}C(1-V_{\phi}(s,\cdot-\sqrt{2}s))=C(\phi)$, we get by the above equality that
$$\lim_{s\to+\infty}\mathrm{E}\left[\mathrm{e}^{-\langle \phi,\sum_{i\ge 1}\mathcal{T}_{d_{i}}\triangle^{s}_{i}\rangle}\right]=\pp_{\delta_{x}}\left[\mathrm{e}^{-C(\phi)\partial M_{\infty}}\right]=\mathrm{E}\left[\mathrm{e}^{-\langle \phi,\mathcal{E}_{\infty}\rangle}\right]$$
for all $\phi\in C^{+}_{c}(\R)$. Hence we prove \eqref{them3.2}.\qed

\bigskip

\begin{remark}\label{rm5}\rm
We claim that for for each $i\ge 1$, $\triangle^{s}_{i}$ converges in distribution to the null measure as $s\to +\infty$.
This is because, by \eqref{decomforu1} for all $\phi\in C^{+}_{c}(\R)$,
$$\mathrm{E}\left[
\mathrm{e}^{-\langle\phi,\triangle^{s}_{i}\rangle}\right]
=\pp_{\cdot,\delta_{0}}\left[\mathrm{e}^{-\langle\phi,I_{s}-\sqrt{2}s\rangle}\right]=V_{\phi}(s,-\sqrt{2}s)
=1-\left(u_{\phi}(s,-\sqrt{2}s)-u^{*}_{\phi}(s,-\sqrt{2}s)\right).
$$
Noticing that  $0\le u^{*}_{\phi}(s,x-\sqrt{2}s)\le u_{\phi}(s,x-\sqrt{2}s)$, \eqref{lem3.8.1} implies that
$$\lim_{s\to+\infty}u_{\phi}(s,-\sqrt{2}s)=\lim_{s\to+\infty}u^*_{\phi}(s,-\sqrt{2}s)=0.$$ Then we have
$$\mathrm{E}\left[
\mathrm{e}^{-\langle\phi,\triangle^{s}_{i}\rangle}\right]
\to 1,\quad\mbox{ as }s\to +\infty.$$
\end{remark}

\subsection{Proof of Proposition \ref{prop5}}

\begin{lemma}\label{lem4.11} Suppose $\phi\in\mathcal{H}$. Then for all $s>0$ and $y\in\mathbb{R}$,
\begin{equation}\label{lem4.11.1}
\mathrm{E}\left[\mathrm{e}^{\langle \ln V_{\phi}(s,\cdot-\sqrt{2}s),\triangle^{Z}+y\rangle}\right]=\mathrm{E}\left[\mathrm{e}^{-\langle \phi,\triangle^{I,s}+y\rangle}\right].
\end{equation}
\end{lemma}

\proof Recall the definition of $(\bar{\mathcal{E}}^Z_\infty, Y)$ in Lemma \ref{lem4.3}.
We use $x_{j}\in \bar{\mathcal{E}}^{Z}_{\infty}$ to denote an atom of the random point measure $\bar{\mathcal{E}}^{Z}_{\infty}$.
For any $s>0$,
define random measure $\Theta_{s}$ by
$$\Theta_s:=\sum_{x_j\in \bar{\mathcal{E}}^Z_\infty}\mathcal{T}_{x_j}(I^j_s-\sqrt{2}s),$$
where $I^j$, $j\ge 1$ are i.i.d. copies of $(I,\mathbb{P}_{\cdot,\delta_0})$, and are independent of $(\bar{\mathcal{E}}^Z_\infty,Y)$.
Recall the Laplace functional of $I_{s}$ given in \eqref{def-V}.
It follows from Lemma \ref{lem4.3} that,
for all $f\in\mathcal{H}$, $s>0$ and $\lambda>0$
\begin{eqnarray*}
\mathrm{E}\left[e^{-\langle f , \Theta_s\rangle-\lambda Y}\right]&=&\mathrm{E}\left[e^{\langle \log V_{f}(s,\cdot-\sqrt{2}s),\bar{\mathcal{E}}^Z_\infty\rangle-\lambda Y}\right]\\
&=&\frac{1}{c_*}\left[ C\left(\left(1-V_{f}(s,\cdot-\sqrt{2}s)\right)1_{(-\infty,0]}+1_{(0,{+\infty})}\right)\right.\\
&&\left.-\int_0^{+\infty} e^{-y} C\left(\left(1-V_{f}(s,\cdot-\sqrt{2}s)\right)1_{(-\infty,y/\lambda]}+1_{(y/\lambda,{+\infty})}\right) dy\right]\\
&=&\mathrm{E}\Big[e^{-\langle{f}, \bar{\mathcal{E}}^{I,s}_\infty\rangle-\lambda Y}\Big].
\end{eqnarray*}
This implies that
$(\Theta_s,Y)\stackrel{\mbox{d}}{=}(\bar{\mathcal{E}}^{I,s}_\infty,Y)$.  Therefore, we obtain that
$$\triangle^{I,s}=\bar{\mathcal{E}}^{I,s}_\infty-Y\stackrel{\mbox{d}}{=}\Theta_s-Y=\sum_{x_j\in \triangle^Z}\mathcal{T}_{x_j}(I^j_s-\sqrt{2}s),$$
which implies that,
for all $\phi\in\mathcal{H}$ and $y\in\R$,
\begin{eqnarray*}
\mathrm{E}\left[\mathrm{e}^{-\langle \phi,\triangle^{I,s}+y\rangle}\right]=\mathrm{E}\left[\mathrm{e}^{-\langle \phi,\sum_{x_j\in \triangle^Z}\mathcal{T}_{x_j}(I^j_s-\sqrt{2}s)+y\rangle}\right]=\mathrm{E}\left[\mathrm{e}^{\langle \ln V_{\phi}(s,\cdot-\sqrt{2}s),\triangle^{Z}+y\rangle}\right].
\end{eqnarray*}
Now we finish the proof.\qed

\bigskip

\noindent\textbf{Proof of Proposition \ref{prop5}:}
It is easy to get by Poisson computations that
\begin{equation}\label{prop5.8.1}
\mathrm{E}\left[\mathrm{e}^{-\langle\phi,\sum_{i\ge 1}\mathcal{T}_{e_{i}}\triangle^{I,s}_{i}\rangle}\right]
=\pp_{\delta_{x}}\left[\exp\left\{-c_{*}\partial M_{\infty}\int_{-\infty}^{+\infty}\sqrt{2}\mathrm{e}^{-\sqrt{2}y}\mathrm{E}
\left[1-\mathrm{e}^{-\langle\phi,\triangle^{I,s}+y\rangle}\right]dy\right\}\right]
\end{equation}
for all $\phi\in C^{+}_{c}(\R)$.
It follows from Lemma \ref{lem3.8}, Lemma \ref{lem4.8} and Lemma \ref{lem4.11} that
for all $\phi\in C^{+}_{c}(\R)$,
\begin{eqnarray}
C(\phi)&=&\lim_{s\to+\infty}C(1-V_{\phi}(s,\cdot-\sqrt{2}s))\nonumber\\
&=&\lim_{s\to+\infty}c_{*}\int_{-\infty}^{+\infty}\sqrt{2}\mathrm{e}^{-\sqrt{2}y}\mathrm{E}\left[1-\mathrm{e}^{\langle \ln V_{\phi}(s,\cdot-\sqrt{2}s),\triangle^{Z}+y\rangle}\right]dy\nonumber\\
&=&\lim_{s\to+\infty}c_{*}\int_{-\infty}^{+\infty}\sqrt{2}\mathrm{e}^{-\sqrt{2}y}\mathrm{E}\left[1-\mathrm{e}^{-\langle \phi,\triangle^{I,s}+y\rangle}\right]dy.\label{lem4.12.2}
\end{eqnarray}
This together with \eqref{prop5.8.1} and Theorem \ref{them2} yields that
$$\lim_{s\to+\infty}\mathrm{E}\left[\mathrm{e}^{-\langle\phi,\sum_{i\ge 1}\mathcal{T}_{e_{i}}\triangle^{I,s}_{i}\rangle}\right]=\pp_{\delta_{x}}\left[\mathrm{e}^{-C(\phi)\partial M_{\infty}}\right]=\mathrm{E}\left[\mathrm{e}^{-\langle\phi,\mathcal{E}_{\infty}\rangle}\right],\quad\forall \phi\in C^{+}_{c}(\R).$$
Hence we complete the proof.\qed

\bigskip

\subsection{Proof of Theorem \ref{them5}}

We prove Theorem \ref{them5} in this section.
First we
relate
$C(\phi)$ to the Laplace functional of a certain random Radon measure. Then we observe that this random measure is infinitely divisible and thus get an expression for $C(\phi)$ which leads to the probabilistic interpretation presented in Theorem \ref{them5}.  Our observation on $C(\phi)$ is inspired by the work of \cite{M}.

\begin{lemma}\label{lem4.12}
Let $\{e_{i}:i\ge1\}$ be the atoms of a Poisson point process with intensity $c_{*}\sqrt{2}\mathrm{e}^{-\sqrt{2}y}dy$ and for every $s>0$, $\{\triangle^{I,s}_{i}:i\ge 1\}$ be an independent sequence of i.i.d. random measures with the same law as $\triangle^{I,s}$. Set
$$\mathcal{D}_{s}:=\sum_{i\ge 1}\mathcal{T}_{e_{i}}\triangle^{I,s}_{i}.$$
Then as $s\to +\infty$, the random measures $\mathcal{D}_{s}$ converges in distribution to a random Radon measure $\mathcal{D}_{\infty}$. Moreover,
we have for any $\phi\in \mathcal{H}$,
\begin{equation}\label{lem4.12.0}
\mathrm{E}\left[\mathrm{e}^{-\langle \phi,\mathcal{D}_{\infty}\rangle}\right]=\lim_{s\to+\infty}\mathrm{E}\left[\mathrm{e}^{-\langle \phi,\mathcal{D}_{s}\rangle}\right]=\mathrm{e}^{-C(\phi)}.
\end{equation}
\end{lemma}

\proof By the definition of $\mathcal{D}_{s}$, one can use simple Poisson computations to get that
\begin{equation}\label{lem4.12.1}
\mathrm{E}\left[\mathrm{e}^{-\langle f,\mathcal{D}_{s}\rangle}\right]=\exp\left\{-c_{*}\int_{-\infty}^{+\infty}\sqrt{2}\mathrm{e}^{-\sqrt{2}y}\mathrm{E}\left[1-\mathrm{e}^{-\langle f,\triangle^{I,s}+y\rangle}\right]dy\right\},\quad\forall f\in\mathcal{B}^{+}(\mathbb{R}).
\end{equation}
Combining \eqref{lem4.12.1} and \eqref{lem4.12.2} we get that
\begin{equation}\label{lem4.12.3}
\lim_{s\to+\infty}\mathrm{E}\left[\mathrm{e}^{-\langle \phi,\mathcal{D}_{s}\rangle}\right]=\mathrm{e}^{-C(\phi)},\quad\forall \phi\in\mathcal{H}.
\end{equation}
The above identity holds in particular for $\phi\in C^{+}_{c}(\mathbb{R})$. Moreover one has by Lemma \ref{lem3.6} that $\lim_{\lambda\to 0+}C(\lambda\phi)=0$ for $\phi\in C^{+}_{c}(\R)$. This implies the existence of the limit $\mathcal{D}_{\infty}$, and \eqref{lem4.12.0} follows immediately from \eqref{lem4.12.3}.\qed

\bigskip

A random measure $\mu$ is said to be exp-$\sqrt{2}$-stable if for any $a,b$ satisfying $\mathrm{e}^{\sqrt{2}a}+\mathrm{e}^{\sqrt{2}b}=1$, it holds that $\mathcal{T}_{a}\mu+\mathcal{T}_{b}\hat{\mu}\stackrel{\mbox{d}}{=}\mu$, where $\hat{\mu}$ is an independent copy of $\mu$. It is easy to see that an exp-$\sqrt{2}$-stable random measure is infinitely divisible. By \eqref{lem4.12.0} and the fact that $C(\mathcal{T}_{a}\phi)=\mathrm{e}^{\sqrt{2}a}C(\phi)$ for all $a\in\R$ and $\phi\in\mathcal{H}$, one can easily show that $\mathcal{D}_{\infty}$ is an exp-$\sqrt{2}$-stable random measure on $\R$.

\begin{lemma}\label{lem5.11}
There exist some constant $\iota\ge 0$ and some measure $\Lambda$ on $\mathcal{M}(\R)\setminus\{0\}$ satisfying that
$$\int_{-\infty}^{+\infty}\mathrm{e}^{-\sqrt{2}x}dx\int_{\mathcal{M}(\R)\setminus\{0\}}\left(1\wedge \mathcal{T}_{x}\mu(A)\right)\Lambda(d\mu)<+\infty,\quad\forall \mbox{ bounded Borel set }A\subset\R,$$
such that
\begin{equation}\label{lem5.11.1}
C(\phi)=\iota\int_{-\infty}^{+\infty}\phi(x)\mathrm{e}^{-\sqrt{2}x}dx+c_{*}\int_{-\infty}^{+\infty}\sqrt{2}\mathrm{e}^{-\sqrt{2}x}dx\int_{\mathcal{M}(\R)\setminus\{0\}}\left(1-\mathrm{e}^{-\langle \phi,\mathcal{T}_{x}\mu\rangle}\right)\Lambda(d\mu),\quad\forall \phi\in C^{+}_{c}(\R).
\end{equation}
Moreover, it holds that
\begin{equation}\label{lem5.11.2}
\iota=\frac{2c_{*}}{1-\mathrm{e}^{-\sqrt{2}}}\lim_{\epsilon\to 0+}\mathrm{limsi}_{s\to+\infty}\int_{-\infty}^{+\infty}\mathrm{e}^{-\sqrt{2}x}\mathrm{E}\left[\langle 1_{(0,1)},\mathcal{T}_{x}\triangle^{I,s}\rangle;\langle 1_{(0,1)},\mathcal{T}_{x}\triangle^{I,s}\rangle<\epsilon\right]dx,
\end{equation}
where ``$\mathrm{limsi}$" is supposed to hold with both $\liminf$ and $\limsup$, and
\begin{eqnarray}
\label{lem5.11.3}
&&\int_{-\infty}^{+\infty}\mathrm{e}^{-\sqrt{2}x}dx\int_{\mathcal{M}(\R)\setminus\{0\}}f\left(\mathcal{T}_{x}\mu(A_{1}),\cdots,\mathcal{T}_{x}\mu(A_{n})\right)\Lambda(d\mu)
\nonumber\\
&=&\lim_{s\to+\infty}\int_{-\infty}^{+\infty}\mathrm{e}^{-\sqrt{2}x}\mathrm{E}\left[f\left(\mathcal{T}_{x}\triangle^{I,s}(A_{1}),\cdots,\mathcal{T}_{x}\triangle^{I,s}(A_{n})\right)\right]dx
\end{eqnarray}
for any $n\ge 1$,
$f\in C_{c}\left(\bar{\R}^{n}\setminus\{0\}\right)$ and
any bounded open sets
$A_{1},\cdots,A_{n}\subset \R$.
\end{lemma}

\proof Since $\mathcal{D}_{\infty}$ is an exp-$\sqrt{2}$-stable random measure, \eqref{lem5.11.1} is a direct result of \cite[Theorem 3.1]{M}. We note that $\mathcal{D}^{s}$ converges in distribution to $\mathcal{D}_{\infty}$ and that
the exponent of the Laplace functional of $\mathcal{D}^{s}$ is given by
$$-\ln \mathrm{E}\left[\mathrm{e}^{-\langle \phi,\mathcal{D}^{s}\rangle}\right]=C\left(1-V_{\phi}(s,\cdot-\sqrt{2}s)\right)=c_{*}\int_{-\infty}^{+\infty}\sqrt{2}\mathrm{e}^{-\sqrt{2}x}
\mathrm{E}\left[1-\mathrm{e}^{-\langle\phi,\mathcal{T}_{x}\triangle^{I,s}\rangle}\right]dx$$
for all $\phi\in C^{+}_{c}(\R)$. \eqref{lem5.11.2} and \eqref{lem5.11.3} follow immediately from \cite[Exercise 6.4]{Kallenberg}.\qed

\bigskip

Lemma \ref{lem5.11} yields a short proof for Theorem \ref{them5}.

\noindent\textbf{Proof of Theorem \ref{them5}:}
We have that
$$\mathrm{E}\left[\mathrm{e}^{-\langle\phi,\mathcal{E}_{\infty}\rangle}\right]=\pp_{\delta_{x}}\left[\mathrm{e}^{-\partial M_{\infty}C(\phi)}\right],\quad\forall \phi\in C^{+}_{c}(\R)$$
where $C(\phi)$ can be represented by \eqref{lem5.11.1}. This yields the result of this theorem.\qed

\bigskip

\begin{remark}\label{rm6}\rm
It follows by \eqref{lem5.11.1} that for all $\lambda>0$,
\begin{equation}\label{r5.1}
\frac{C(\lambda\phi)}{\lambda}=\iota\int_{-\infty}^{+\infty}\phi(x)\mathrm{e}^{-\sqrt{2}x}dx+c_{*}\int_{-\infty}^{+\infty}\sqrt{2}\mathrm{e}^{-\sqrt{2}x}dx\int_{\mathcal{M}(\R)\setminus\{0\}}
\frac{1-\mathrm{e}^{-\lambda \langle\phi,\mathcal{T}_{x}\mu\rangle}}{\lambda}\Lambda(d\mu),\quad\forall \phi\in C^{+}_{c}(\R).
\end{equation}
Note that
$$\int_{-\infty}^{+\infty}\mathrm{e}^{-\sqrt{2}x}dx\int_{\mathcal{M}(\R)\setminus\{0\}}
\frac{1-\mathrm{e}^{-\lambda \langle\phi,\mathcal{T}_{x}\mu\rangle}}{\lambda}\Lambda(d\mu)\le\int_{-\infty}^{+\infty}\mathrm{e}^{-\sqrt{2}x}dx\int_{\mathcal{M}(\R)\setminus\{0\}}\frac{1}{\lambda}\wedge \langle\phi,\mathcal{T}_{x}\mu\rangle\Lambda(d\mu).$$
Thus by the dominated convergence theorem, the second term of \eqref{r5.1} converges to $0$ as $\lambda\to+\infty$. Consequently we get that
\begin{equation}
\iota\int_{-\infty}^{+\infty}\phi(x)\mathrm{e}^{-\sqrt{2}x}dx=\lim_{\lambda\to+\infty}\frac{C(\lambda\phi)}{\lambda},\quad\forall \phi\in C^{+}_{c}(\R).
\end{equation}
So a sufficient condition for the constant $\iota$ to be $0$ is that
$$\sup_{\lambda}C(\lambda\phi)<+\infty\quad\mbox{ for some }\phi\in C^{+}_{c}(\R).$$
This is true if the branching mechanism $\psi$ satisfies (A3), where one has $\sup_{\lambda}C(\lambda 1_{(0,+\infty)})\le \tilde{c}_{0}<+\infty$. In this case, $\mathcal{E}_{\infty}$ is equal in law to a Poisson random measure on $\mathcal{M}(\R)$ with intensity $c_{*}\partial M_{\infty}\int_{-\infty}^{+\infty}\sqrt{2}\mathrm{e}^{-\sqrt{2}x}\mathcal{T}_{x}\Lambda(d\mu)dx.$
We shall discuss this special case in detail in the next section.
\end{remark}

\subsection{Special case where (A3) is satisfied and proof of Theorem \ref{them4}}

In this section we assume in addition that (A3) is satisfied.
Recall the definition of
$\tilde{u}_{\phi}(t,x)$
given in \eqref{r0}. It is shown in the argument of Remark \ref{rm4.13} that when (A3) holds, $\tilde{u}_{\phi}(r,\cdot-\sqrt{2}r)\in \mathcal{H}$ for all $\phi\in\mathcal{H}$ and $r>0$. Applying Lemma \ref{lem3.8} to the function $u_{\tilde{u}_{\phi}(r,\cdot-\sqrt{2}r)}(t,x)=\tilde{u}_{\phi}(t+r,x-\sqrt{2}r)$, one gets that
for all $x\in\R$ and $\phi\in\mathcal{H}$,
\begin{equation}\nonumber
\lim_{t\to+\infty}\frac{t^{3/2}}{\frac{3}{2\sqrt{2}}\log t}\mathrm{e}^{-\sqrt{2}x}\tilde{u}_{\phi}(t,x-\sqrt{2}t)=\tilde{C}(\phi),
\end{equation}
where $\widetilde{C}(\phi)$ is defined in \eqref{defitildec}. In particular by taking $\phi=0$, one gets that
\begin{equation}\label{6.1}
\lim_{t\to+\infty}\frac{t^{3/2}}{\frac{3}{2\sqrt{2}}\log t}\mathrm{e}^{-\sqrt{2}x}\pp_{\delta_{x}}\left(\max X_{t}-\sqrt{2}t>0\right)=\tilde{c}_{0},
\end{equation}
where $\tilde{c}_{0}=\widetilde{C}(0)$.

\begin{lemma}\label{lem6.1}
Suppose $x,z\in \R$. Under $\pp_{\delta_{x}}$, conditioned on $\{\max X_{t}-\sqrt{2}t-z>0\}$, the random elements $\left(X_{t}-\sqrt{2}t-z,Z_{t}-\sqrt{2}t-z,\max X_{t}-\sqrt{2}t-z\right)$ converges, as $t\to+\infty$, in distribution to a limit
 $\left(\widetilde{\mathcal{E}}^{X}_{\infty},\widetilde{\mathcal{E}}^{Z}_{\infty},Y\right)$, where the limit is independent of $x$ and $z$, and $Y$ is an exponential random variable with mean $1/\sqrt{2}$. Moreover, given $\left(\widetilde{\mathcal{E}}^{X}_{\infty},Y\right)$, $\widetilde{\mathcal{E}}^{Z}_{\infty}$ is a Poisson random measure with intensity $\widetilde{\mathcal{E}}^{X}_{\infty}$.
\end{lemma}

\proof Fix $x,z\in \R$. It has been proved in \cite[Proposition 3.4]{RSZ} that conditioned on $\{\max X_{t}-\sqrt{2}t-z>0\}$, $\left(X_{t}-\sqrt{2}t-z,\max X_{t}-\sqrt{2}t-z\right)$ converges, as $t\to+\infty$, in distribution to a limit $\left(\widetilde{\mathcal{E}}^{X}_{\infty},Y\right)$, where the limit is independent of $x$ and $z$, and $Y$ is an exponential random variable with mean $1/\sqrt{2}$. We note that for all $f,g\in C^{+}_{c}(\R)$ and $\lambda_{i}\ge 0$, $i=1,2,3$,
\begin{eqnarray}
&&\pp_{\delta_{x}}\left[\mathrm{e}^{-\lambda_{1}\langle f,X_{t}-\sqrt{2}t-z\rangle-\lambda_{2}\langle g,Z_{t}-\sqrt{2}t-z\rangle-\lambda_{3}\left(\max X_{t}-\sqrt{2}t-z\right)}\,|\,\max X_{t}-\sqrt{2}t-z>0\right]\nonumber\\
&=&\pp_{\delta_{x}}\left[\mathrm{e}^{-\lambda_{1}\langle f,X_{t}-\sqrt{2}t-z\rangle-\lambda_{3}\left(\max X_{t}-\sqrt{2}t-z\right)}\pp_{\delta_{x}}\left[\mathrm{e}^{-\lambda_{2}\langle g,Z_{t}-\sqrt{2}t-z\rangle}\,|\,X_{t}\right]\,|\,\max X_{t}-\sqrt{2}t-z>0\right]\nonumber\\
&=&\pp_{\delta_{x}}\left[\mathrm{e}^{-\langle \lambda_{1}f+1-\mathrm{e}^{-\lambda_{2}g},X_{t}-\sqrt{2}t-z\rangle-\lambda_{3}\left(\max X_{t}-\sqrt{2}t-z\right)}\,|\,\max X_{t}-\sqrt{2}t-z>0\right]\nonumber\\
&\to&\mathrm{E}\left[\mathrm{e}^{-\langle\lambda_{1}f+1-\mathrm{e}^{-\lambda_{2}g},\tilde{\mathcal{E}}^{X}_{\infty}\rangle-\lambda_{3}Y}\right],\quad\mbox{ as }t\to +\infty.\label{lem6.1.1}
\end{eqnarray}
Obviously the right hand side of \eqref{lem6.1.1} converges to $1$ as $\lambda_{i}\to 0$, $i=1,2,3$. Hence conditioned on $\{\max X_{t}-\sqrt{2}t-z>0\}$, $\left(X_{t}-\sqrt{2}t-z,Z_{t}-\sqrt{2}t-z,\max X_{t}-\sqrt{2}t-z\right)$ converges in distribution to a limit $\left(\widetilde{\mathcal{E}}^{X}_{\infty},\widetilde{\mathcal{E}}^{Z}_{\infty},Y\right)$
with $Y$ being an exponential random variable with mean $1/\sqrt{2}$.
Moreover, it holds that
\begin{equation*}
\mathrm{E}\left[\mathrm{e}^{-\langle f,\widetilde{\mathcal{E}}^{X}_{\infty}\rangle-\langle g,\widetilde{\mathcal{E}}^{Z}_{\infty}\rangle-\lambda Y}\right]=\mathrm{E}\left[\mathrm{e}^{-\langle f+1-\mathrm{e}^{-g},\widetilde{\mathcal{E}}^{X}_{\infty}\rangle-\lambda Y}\right]
\end{equation*}
for all $f,g\in C^{+}_{c}(\R)$ and $\lambda\ge 0$. In particular one has
\begin{equation*}
\mathrm{E}\left[\mathrm{e}^{-\langle f,\widetilde{\mathcal{E}}^{X}_{\infty}\rangle-\lambda Y}\mathrm{E}\left[\mathrm{e}^{-\langle g,\widetilde{\mathcal{E}}^{Z}_{\infty}\rangle}\,|\,\widetilde{\mathcal{E}}^{X}_{\infty},Y\right]\right]=\mathrm{E}\left[\mathrm{e}^{-\langle f,\widetilde{\mathcal{E}}^{X}_{\infty}\rangle-\lambda Y}\cdot\mathrm{e}^{-\langle 1-\mathrm{e}^{-g},\widetilde{\mathcal{E}}^{X}_{\infty}\rangle}\right].
\end{equation*}
This implies that $\mathrm{E}\left[\mathrm{e}^{-\langle g,\widetilde{\mathcal{E}}^{Z}_{\infty}\rangle}\,|\,\widetilde{\mathcal{E}}^{X}_{\infty},Y\right]=\mathrm{e}^{-\langle 1-\mathrm{e}^{-g},\widetilde{\mathcal{E}}^{X}_{\infty}\rangle}$
$\mathrm{P}$-a.s.
So we prove the second conclusion of this lemma.\qed

\bigskip

\begin{lemma}\label{lem6.2}
Suppose $x,z\in \R$. Under $\pp_{\delta_{x}}$, conditioned on $\{\max X_{t}-\sqrt{2}t-z>0\}$, the random elements
$(X_{t}-\max X_{t},Z_{t}-\max X_{t},\max X_{t}-\sqrt{2}t-z)$
converges, as $t\to+\infty$, in distribution to a limit
$(\widetilde{\triangle}^{X},\widetilde{\triangle}^{Z},Y):=(\widetilde{\mathcal{E}}^{X}_{\infty}-Y,\widetilde{\mathcal{E}}^{Z}_{\infty}-Y,Y)$, where the limit is independent of $x$ and $z$, and $(\widetilde{\triangle}^{X},\widetilde{\triangle}^{Z})$ is independent of $Y$. Moreover, given $\widetilde{\triangle}^{X}$, $\widetilde{\triangle}^{Z}$ is a Poisson random measure with intensity $\widetilde{\triangle}^{X}$.
\end{lemma}

\proof The first conclusion follows from Lemma \ref{lem6.1} (in place of Lemma \ref{lem4.3} ) in the same way as Lemma \ref{lem4.4}. We only need to show the second conclusion.

Since $\left(\widetilde{\triangle}^{X},\widetilde{\triangle}^{Z}\right)$ is independent of $Y$, we have for all $f\in C^{+}_{c}(\R)$,
\begin{eqnarray}
\mathrm{E}\left[\mathrm{e}^{-\langle f,\widetilde{\triangle}^{Z}\rangle}\,|\,\widetilde{\triangle}^{X}\right]&=&\mathrm{E}\left[\mathrm{e}^{-\langle f,\widetilde{\triangle}^{Z}\rangle}\,|\,\widetilde{\triangle}^{X},Y\right]
=
\mathrm{E}\left[\mathrm{e}^{-\langle f,\widetilde{\mathcal{E}}^{Z}_{\infty}-Y\rangle}\,|\,\widetilde{\mathcal{E}}^{X}_{\infty},Y\right]\nonumber\\
&=&\mathrm{e}^{-\langle 1-\mathrm{e}^{-f},\widetilde{\mathcal{E}}^{X}_{\infty}-Y\rangle}
=
\mathrm{e}^{-\langle 1-\mathrm{e}^{-f},\widetilde{\triangle}^{X}\rangle}.\nonumber
\end{eqnarray}
The third equality follows from the second conclusion of Lemma \ref{lem6.1}. Hence we prove the second conclusion.\qed

\begin{lemma}\label{lem6.3}
For any $f,g\in C^{+}_{c}(\R)$,
\begin{equation*}
C\left(f+1-\mathrm{e}^{-g}\right)=\tilde{c}_{0}\int_{-\infty}^{+\infty}\sqrt{2}\mathrm{e}^{-\sqrt{2}y}\mathrm{E}\left[1-\mathrm{e}^{-\langle f,\widetilde{\triangle}^{X}+y\rangle-\langle g,\widetilde{\triangle}^{Z}+y\rangle}\right]dy.
\end{equation*}
\end{lemma}

\proof Fix arbitrary $f,g\in C^{+}_{c}(\R)$.
It follows by \cite[(3.12)]{RSZ} that
	 \begin{equation*}
	 	C\left(f+1-\mathrm{e}^{-g}\right)=\tilde{c}_{0}\int_{-\infty}^{+\infty}\sqrt{2}\mathrm{e}^{-\sqrt{2}y}\mathrm{E}\left[1-\mathrm{e}^{-\langle f+1-e^{-g},\widetilde{\triangle}^{X}+y\rangle}\right]dy.
	 \end{equation*}
 By Lemma \ref{lem6.2}, we have that
 $$\mathrm{E}\left[1-\mathrm{e}^{-\langle f,\widetilde{\triangle}^{X}+y\rangle-\langle g,\widetilde{\triangle}^{Z}+y\rangle}\right]=\mathrm{E}\left[1-\mathrm{e}^{-\langle f+1-e^{-g},\widetilde{\triangle}^{X}+y\rangle}\right].$$
 Now the desired result follows immediately.\qed
\bigskip

\noindent\textbf{Proof of Theorem \ref{them4}:}
Using  computation on Poisson point process,
we have for all $f,g\in C^{+}_{c}(\R)$,
\begin{eqnarray}
\mathrm{E}\left[\mathrm{e}^{-\left\langle f,\sum_{i\ge 1}\mathcal{T}_{\tilde{e}_{i}}\widetilde{\triangle}^{X}_{i}\right\rangle-\left\langle g,\sum_{i\ge 1}\mathcal{T}_{\tilde{e}_{i}}\widetilde{\triangle}^{Z}_{i}\right\rangle}\right]
&=&\mathrm{E}\left[\mathrm{e}^{\left\langle\ln \mathrm{E}\left[\mathrm{e}^{-\left\langle f,\widetilde{\triangle}^{X}+\cdot\right\rangle-\left\langle g,\widetilde{\triangle}^{Z}+\cdot\right\rangle}\right],\sum_{i\ge 1}\delta_{\tilde{e}_{i}}\right\rangle}\right]\nonumber\\
&=&\pp_{\delta_{x}}\left[\mathrm{e}^{-\tilde{c}_{0}\partial M_{\infty}\int_{-\infty}^{+\infty}\sqrt{2}\mathrm{e}^{-\sqrt{2}y}\mathrm{E}\left[1-\mathrm{e}^{-\left\langle f,\widetilde{\triangle}^{X}+y\right\rangle-\left\langle g,\widetilde{\triangle}^{Z}+y\right\rangle}\right]dy}\right].\nonumber
\end{eqnarray}
This together with Theorem \ref{them2} and Lemma \ref{lem6.3} yields that
$$\mathrm{E}\left[\mathrm{e}^{-\langle f,\mathcal{E}_{\infty}\rangle-\langle g,\mathcal{E}^{Z}_{\infty}\rangle}\right]=\mathrm{E}\left[\mathrm{e}^{-\left\langle f,\sum_{i\ge 1}\mathcal{T}_{\tilde{e}_{i}}\widetilde{\triangle}^{X}_{i}\right\rangle-\left\langle g,\sum_{i\ge 1}\mathcal{T}_{\tilde{e}_{i}}\widetilde{\triangle}^{Z}_{i}\right\rangle}\right],\quad\forall f,g\in C^{+}_{c}(\R).$$
Hence we complete the proof.\qed

\appendix
\section{Appendix}

\begin{lemma}\label{lemA1}
Let $L$ be the integer-valued random variable with distribution $\{p_{k}:k\ge 2\}$
as defined in Proposition \ref{prop1}.
 Suppose $f:[0,+\infty)\to [0,+\infty)$ satisfies the following conditions: There exist some constants $c,\kappa\ge0$ such that
\begin{description}
\item{(1)} $f$ is bounded in $[0,c)$ and convex on $[c,+\infty)$.
\item{(2)} $f(xy)\le \kappa f(x)f(y)$ for all $x,y\in [c,+\infty)$.
\end{description}
Then the following statements are equivalent.
\begin{description}
\item{(i)} $\mathrm{E}\left[f(L)\right]<+\infty$.
\item{(ii)} $\pp_{\cdot,k\delta_{0}}\left[f(\|Z_{t}\|)\right]<\infty$ for all $t>0$ and $k\in\mathbb{N}$.
\item{(iii)} $\pp_{\mu}\left[f(\|Z_{t}\|)\right]<\infty$ for all $t>0$ and $\mu\in\mc$.
\item{(iv)} $\pp_{\mu}\left[f(\|X_{t}\|)\right]<\infty$ for all $t>0$ and $\mu\in\mc$.
\item{(v)} $\int_{(1,+\infty)}f(x)\pi(dx)<+\infty.$
\end{description}
\end{lemma}

\proof
 Without loss of generality we may and do assume that the function $f:[0,+\infty)\mapsto [0,+\infty)$ satisfies that there is some $\kappa >0$ such that
\begin{description}
\item{(1') $f$ is convex on $[0,+\infty)$;}
\item{(2') $f(xy)\le \kappa f(x)f(y)$ for all $x,y\in [0,+\infty)$;}
\item{(3') $f$ is nondecreasing on $[0,+\infty)$ and $f(x)>1$ for all $x\ge 0$.}
\end{description}
In fact, \cite[Chapter IV, Lemma 1]{AN} shows that for any function $f$ which satisfies the original hypothesis, there is a function $\tilde{f}$ satisfying (1')-(3') such that for any probability measure $\mu$ on $[0,+\infty)$,
$\int_{[0,+\infty)}f(x)\mu(dx)$ is finite if and only if $\int_{[0,+\infty)}\tilde{f}(x)\mu(dx)$ is finite.

That $(i)\Leftrightarrow (ii)$ is established in \cite[Chapter III, Theorem 2]{AN}. We note that $(\|X_{t}\|)_{t\ge 0}$ is a continuous-state branching process. Thus $(iv)\Leftrightarrow (v)$ follows directly from \cite[Theorem 2.1]{JL}.

(ii) $\Leftrightarrow$ (iii): Since the branching rate and offspring distribution of $(Z_{t})_{t\ge 0}$ is spatially-independent, $(\|Z_{t}\|)_{t\ge 0}$ is a continuous-time Galton-Watson branching process. Thus we have for any nontrivial $\mu\in \mc$,
$$\pp_{\mu}\left[f(\|Z_{t}\|)\right]=\pp_{\mu}\left[\pp_{\mu}\left[f(\|Z_{t}\|)|Z_{0}\right]\right]=\pp_{\mu}\left[\pp_{\cdot,\|Z_{0}\|\delta_{0}}\left[f(\|Z_{t}\|)\right]\right].$$
Note that ($\|Z_{0}\|,\pp_{\mu})$ is a Poinsson random variable with parameter $\|\mu\|$. Thus we have
\begin{equation}\label{A1.1}
\pp_{\mu}\left[f(\|Z_{t}\|)\right]=\sum_{k=0}^{+\infty}\frac{\|\mu\|^{k}}{k!}\mathrm{e}^{-\|\mu\|}\pp_{\cdot,k\delta_{0}}\left[f(\|Z_{t}\|)\right].
\end{equation}
Hence  $\pp_{\mu}\left[f(\|Z_{t}\|)\right]<+\infty$ if and only if $\pp_{\cdot,k\delta_{0}}\left[f(\|Z_{t}\|)\right]<+\infty$ for all $k\in\mathbb{N}$.

(iii) $\Rightarrow$ (iv): Fix $t>0$ and a nontrivial $\mu\in\mc$. Suppose that $\pp_{\mu}\left[f(\|Z_{t}\|)\right]<+\infty$.
Recall that given $X_{t}$, $Z_{t}$ is a Poisson random measure with intensity $X_{t}(dx)$, and so $\|Z_{t}\|$ is a Poisson random variable with parameter $\|X_{t}\|$. Thus we get
$$\pp_{\mu}\left[\|Z_{t}\||X_{t}\right]=\|X_{t}\|\quad\pp_{\mu}\mbox{-a.s.}$$
Since $f$ in convex on $[0,+\infty)$, it follows by Jensen's inequality that
$$\pp_{\mu}\left[f(\|X_{t}\|)\right]=\pp_{\mu}\left[f\left(\pp_{\mu}\left[\|Z_{t}\||X_{t}\right]\right)\right]\le\pp_{\mu}\left[\pp_{\mu}\left[f(\|Z_{t}\|)|X_{t}\right]\right]
=\pp_{\mu}\left[f(\|Z_{t}\|)\right]<+\infty.$$

(iv) $\Rightarrow$ (iii): Fix $t,s>0$ and a nontrivial $\mu\in\mc$. Suppose $\pp_{\mu}\left[f(\|X_{t+s}\|)\right]<+\infty$. It follows by \eqref{skeletondecomp} that
$$\|X_{t+s}\|=\|X^{*}_{t+s}\|+\|I^{*,t}_{s}\|+\sum_{u\in Z_{t}}\|I^{(u)}_{s}\|.$$
We have
\begin{eqnarray}
\pp_{\mu}\left[f(\|X_{t+s}\|)\right]&\ge&\pp_{\mu}\left[f\left(\sum_{u\in Z_{t}}\|I^{(u)}_{s}\|\right)\right]
=
\pp_{\mu}\left[\pp_{\mu}\left[f\left(\sum_{u\in Z_{t}}\|I^{(u)}_{s}\|\right)|Z_{t}\right]\right]\nonumber\\
&\ge&\pp_{\mu}\left[f\left(\pp_{\mu}\left[\sum_{u\in Z_{t}}\|I^{(u)}_{s}\||Z_{t}\right]\right)\right]
=
\pp_{\mu}\left[f\left(\sum_{u\in Z_{t}}\pp_{\cdot,\delta_{z_{u}(t)}}\left[\|I_{s}\|\right]\right)\right].\nonumber
\end{eqnarray}
The first inequality follows from condition (3'),  the second inequality from Jensen's inequality.
We observe that the distribution of $\|I_{s}\|$ under $\pp_{\cdot,\delta_{x}}$ is independent of the starting location $x$. If we define $g(s):=\pp_{\cdot,\delta_{x}}\left[\|I_{s}\|\right]$ for all $s>0$, then we get from the above argument that
$$\pp_{\mu}\left[f\left(g(s)\|Z_{t}\|\right)\right]=\pp_{\mu}\left[f\left(\sum_{u\in Z_{t}}g(s)\right)\right]\le \pp_{\mu}\left[f(\|X_{t+s}\|)\right]<+\infty.$$
Thus by (2') we have
\begin{equation*}
\pp_{\mu}\left[f(\|Z_{t}\|)\right]=\pp_{\mu}\left[f(g(s)g(s)^{-1}\|Z_{t}\|)\right]\le \kappa f(g(s)^{-1})\pp_{\mu}\left[f(g(s)\|Z_{t}\|)\right]<+\infty.
\end{equation*}
Therefore we complete the proof.\qed

\bigskip

\begin{lemma}\label{lemA2}
Condition (A2) holds for some $0<\beta<1$ if and only if
\begin{equation}\label{A2.0}
\int_{0}^{1}\left(1+\psi'(s)\right)s^{-(1+\beta)}ds<+\infty.
\end{equation}
\end{lemma}

\proof Let $L$ be the integer-valued random variable with distribution $\{p_{k}:k\ge 2\}$. It follows from Lemma \ref{lemA1} that (A2) holds if and only if
\begin{equation}
\mathrm{E}\left[L^{1+\beta}\right]<+\infty.\label{A2.1}
\end{equation}
Let $\varphi(s)$ be the Laplace transform of $L$, that is, $\varphi(s):=\mathrm{E}\left[\mathrm{e}^{-s L}\right]$ for all $s\ge 0$. Then by \cite[Theorem B]{BD}, \eqref{A2.1} is equivalent to that
\begin{equation}\label{A2.2}
\int_{0}^{1}f_{1}(s)s^{-(2+\beta)ds}<+\infty,
\end{equation}
where $f_{1}(s):=\varphi(s)-1-\varphi'(0)s$. We use $F(s)$ to denote the generating function of $L$, i.e., $F(s)=\mathrm{E}\left[s^{L}\right]$ for $s\in [0,1]$. Since $\varphi(s)=F(\mathrm{e}^{-s})$ for all $s\ge 0$, setting $u=1-\mathrm{e}^{-s}$ we have
$$\varphi''(s)=F''(1-u)(1-u)^{2}+F'(1-u)(1-u).$$
Since $u\sim s$ as $s\to 0$, one has
\begin{equation}\label{A2.3}
f_{1}(s)\sim \frac{s^{2}}{2}\varphi''(s)\sim \frac{u^{2}}{2}\left[F''(1-u)(1-u)^{2}+F'(1-u)(1-u)\right]\sim c_{1}u^{2}F''(1-u)\quad \mbox{ as }s\to 0,
\end{equation}
for some constant $c_{1}>0$. On the other hand, in view of \eqref{eq2}, one has
$q F''(1-u)=\psi''(u)$. This together with \eqref{A2.3} implies that
$$f_{1}(s)\sim c_{2}u^{2}\psi''(u)\sim c_{2}u\left(\psi'(u)-\psi'(0)\right)=c_{2}u\left(\psi'(u)+1\right)\quad\mbox{ as }s\to 0, $$
for some constant $c_{2}>0$. Hence we have $f_{1}(s)s^{-(2+\beta)}\sim c_{2}u^{-(1+\beta)}(\psi'(u)+1)$ as $s\to 0$. So \eqref{A2.2} holds if and only if \eqref{A2.0} holds.\qed

\medskip

{\bf Yan-Xia Ren}

LMAM School of Mathematical Sciences \& Center for
Statistical Science, Peking
University,

Beijing, 100871, P.R. China.

E-mail: yxren@math.pku.edu.cn

\medskip

{\bf Ting Yang}

School of Mathematics and Statistics, Beijing Institute of Technology, Beijing, 100081, P.R.China;

Beijing Key Laboratory on MCAACI,

Beijing, 100081, P.R. China.

Email: yangt@bit.edu.cn

\medskip

{\bf Rui Zhang}

School of Mathematical Sciences \& Academy for Multidisciplinary Studies, Capital Normal
University,

Beijing, 100048, P.R. China.

Email: zhangrui27@cnu.edu.cn

\end{document}